\documentclass[a4paper,10pt,reqno]{amsart}
\usepackage{amscd,amssymb,graphicx,color}

\def\link{\operatorname{Link}}
\def\oM{\overline{M}}

\def\bM{\overline{\Bbb M}}

\def\bP{\Bbb P}

\def\Norm{\operatorname{Norm}}
\def\bF{\Bbb F}
\def\bD{\Bbb D}

\def\Aut{\operatorname{Aut}}
\def\Ann{\operatorname{Ann}}

\def\Sym{\operatorname{Sym}}
\def\Supp{\operatorname{Supp}}

\def\tY{\tilde{Y}}
\def\tsigma{\tilde{\sigma}}

\def\cO{\Cal O}

\def\codim{\operatorname{codim}}

\def\bZ{\Bbb Z}

\def\bQ{\Bbb Q}
\def\bG{\Bbb G}
\def\bR{\Bbb R}
\def\bA{\Bbb A}
\def\bB{\Bbb B}

\def\oY{\overline{Y}}

\def\cB{\Cal B}
\def\cA{\Cal A}
\def\cD{\Cal D}
\def\cS{\Cal S}
\def\cE{\Cal E}
\def\cF{\Cal F}
\def\tS{\tilde S}

\def\cA{\Cal A}
\def\tcA{\tilde{\cA}}
\def\tpi{\tilde{\pi}}
\def\tp{\tilde{p}}
\def\tA{\tilde{A}}

\def\oS{\overline{S}}
\def\cT{\Cal T}
\def\cP{\Cal P}
\def\cB{\Cal B}
\def\tS{\tilde{S}}

\def\oN{\overline{N}}

\def\tcA{\tilde{\cA}}
\def\cB{\Cal B}
\def\cY{\Cal Y}
\def\ocY{\overline{\cY}}
\def\cD{\Cal D}
\def\cX{\Cal X}
\def\cH{\Cal H}
\def\cM{\Cal M}
\def\cZ{\Cal Z}
\def\ocZ{\overline{\Cal Z}}

\def\Spec{\operatorname{Spec}}

\def\PGL{\operatorname{PGL}}
\def\tM{\tilde{M}}
\def\lc{\operatorname{lc}}
\def\sts{\operatorname{ss}}

\def\cE{\Cal E}

\def\bS{\Bbb S}
\def\bB{\Bbb B}

\def\ord{\operatorname{ord}}

\def\cU{\Cal U}

\def\cR{\Cal R}
\def\cO{\Cal O}

\def\oZ{\overline{Z}}

\def\oR{\overline{R}}

\def\Sp{\operatorname{Sp}}
\def\rk{\operatorname{rk}}

\def\Pic{\operatorname{Pic}}
\def\Spec{\operatorname{Spec}}
\def\Im{\operatorname{Im}}
\def\pr{\operatorname{pr}}
\def\Hom{\operatorname{Hom}}
\def\Cal{\mathcal}

\def\eps{\varepsilon}
\def\hY{\hat Y}
\def\hD{\hat D}
\def\hp{\hat p}

\def\W{\mathop{\circ}}

\def\arrow{\mathop{\longrightarrow}\limits}

\def\Efi#1#2#3#4#5{\displaystyle
#1\!\!-\!\!#2
\!\!-\!\!#3
\!\!-\!\!#4
\hskip-24.2pt\lower4.5pt\hbox{${\scriptstyle|}
\hskip-3.35pt\lower6pt\hbox{$#5$}$}}

\def\Evia#1#2#3#4#5{\displaystyle
#1\!\!-\!\!#2
\!\!-\!\!#3
\hskip-24.2pt\lower4.5pt\hbox{${\scriptstyle|}
\hskip-3.35pt\lower6pt\hbox{$#4\!\!-\!\!\!-\!\!\!-\!\!$}$\hskip2.3pt${\scriptstyle|}
\hskip-3.35pt\lower6pt\hbox{$#5$}$}}

\def\Ezia#1#2#3#4{\displaystyle
#1\!\!-\!\!#2
\hskip-14.8pt\lower4.5pt\hbox{${\scriptstyle|}
\hskip-3.35pt\lower6pt\hbox{$#3\!\!-\!\!$}$\hskip2.3pt${\scriptstyle|}
\hskip-3.35pt\lower6pt\hbox{$#4$}$}}

\def\Efia#1#2#3#4#5#6{\displaystyle
#1\!\!-\!\!#2
\!\!-\!\!#3
\!\!-\!\!#4
\hskip-24.2pt\lower4.5pt\hbox{${\scriptstyle|}
\hskip-3.35pt\lower6pt\hbox{$#5$}$\hskip5.7pt${\scriptstyle|}
\hskip-3.35pt\lower6pt\hbox{$#6$}$}}

\def\Esi#1#2#3#4#5#6{\displaystyle
#1\!\!-\!\!#2
\!\!-\!\!#3
\!\!-\!\!#4\!\!-\!\!#5
\hskip-24.2pt\lower4.5pt\hbox{${\scriptstyle|}
\hskip-3.35pt\lower6pt\hbox{$#6$
\lower3pt\hbox{\ }}$}}

\def\Esia#1#2#3#4#5#6#7{\displaystyle
#1\!\!-\!\!#2
\!\!-\!\!#3
\!\!-\!\!#4\!\!-\!\!#5
\hskip-24.2pt\lower4.5pt\hbox{${\scriptstyle|}
\hskip-3.35pt\lower6pt\hbox{$#6$\hskip-3.8pt\lower4.5pt\hbox{${\scriptstyle|}
\hskip-3.35pt\lower6pt\hbox{$#7$}$}}
\lower3pt\hbox{\ }$}}

\def\Ese#1#2#3#4#5#6#7{\displaystyle
#1\!\!-\!\!#2
\!\!-\!\!#3
\!\!-\!\!#4\!\!-\!\!#5\!\!-\!\!#6
\hskip-33.6pt\lower4.5pt\hbox{${\scriptstyle|}
\hskip-3.35pt\lower6pt\hbox{$#7$
\lower3pt\hbox{\ }
}$}}

\def\Esea#1#2#3#4#5#6#7#8{\displaystyle
#1\!\!-\!\!#2
\!\!-\!\!#3
\!\!-\!\!#4\!\!-\!\!#5\!\!-\!\!#6\!\!-\!\!#7
\hskip-33.6pt\lower4.5pt\hbox{${\scriptstyle|}
\hskip-3.35pt\lower6pt\hbox{$#8$
\lower3pt\hbox{\ }
}$}}

\def\Eei#1#2#3#4#5#6#7#8{\displaystyle
#1\!\!-\!\!#2
\!\!-\!\!#3
\!\!-\!\!#4\!\!-\!\!#5\!\!-\!\!#6\!\!-\!\!#7
\hskip-43.2pt\lower4.5pt\hbox{${\scriptstyle|}
\hskip-3.35pt\lower6pt\hbox{$#8$
\lower3pt\hbox{\ }
}$}}

\def\Eeia#1#2#3#4#5#6#7#8#9{{\displaystyle
#1\!\!-\!\!#2
\!\!-\!\!#3
\!\!-\!\!#4\!\!-\!\!#5\!\!-\!\!#6\!\!-\!\!#7\!\!-\!\!#8
\hskip-52.2pt\lower4.5pt\hbox{${\scriptstyle|}
\hskip-3.35pt\lower6pt\hbox{$#9$
\lower3pt\hbox{\ }
}$}}}

\def\arrow{\mathop{\longrightarrow}\limits}

\def\trop{\operatorname{trop}}
\def\oM{\overline{M}}
\def\os{\overline{S}}

\def\ocM{\overline{\Cal M}}
\def\ocH{\overline{\Cal H}}
\def\tcD{\tilde{\Cal D}}

\def\oH{\overline{H}}
\def\tN{\tilde{N}}
\def\tS{\tilde{S}}
\def\ts{\tS}
\def\ts7{\tilde{S}_7}
\def\tMM{\tilde{M}_3^{(2)}}
\def\tcMM{\tilde{\Cal M}_3^{(2)}}
\def\tcM{\tilde{\Cal M}}

\def\tcH{\tilde{\Cal H}}
\def\cN{\Cal N}

\def\cR{\Cal R}
\def\bP{\Bbb P}
\def\bF{\Bbb F}

\def\tdelta{\tilde{\delta}}
\def\cC{\Cal C}

\def\Sp{\operatorname{Sp}}

\def\Aut{\operatorname{Aut}}

\def\Sym{\operatorname{Sym}}
\def\Supp{\operatorname{Supp}}

\def\Hom{\operatorname{Hom}}

\def\tY{\tilde{Y}}
\def\tsigma{\tilde{\sigma}}

\def\cO{\Cal O}
\def\cS{\Cal S}
\def\cD{\Cal D}

\def\Pic{\operatorname{Pic}}

\def\codim{\operatorname{codim}}

\def\bZ{\Bbb Z}

\def\bQ{\Bbb Q}
\def\bG{\Bbb G}
\def\bR{\Bbb R}
\def\bA{\Bbb A}
\def\bB{\Bbb B}

\def\oY{\overline{Y}}

\def\cB{\Cal B}
\def\cA{\Cal A}
\def\cD{\Cal D}
\def\cS{\Cal S}
\def\cE{\Cal E}
\def\cF{\Cal F}
\def\tS{\tilde S}
\def\tF{\tilde F}
\def\tG{\tilde G}
\def\tcF{\tilde{\Cal F}}
\def\eY{\ddot Y}
\def\ep{\ddot p}

\def\cA{\Cal A}
\def\tcA{\tilde{\cA}}
\def\tA{\tilde{A}}

\def\oS{\overline{S}}

\def\cT{\Cal T}
\def\cP{\Cal P}
\def\cB{\Cal B}
\def\tS{\tilde{S}}

\def\oN{\overline{N}}

\def\tcA{\tilde{\cA}}
\def\cB{\Cal B}
\def\cY{\Cal Y}
\def\cX{\Cal X}
\def\cH{\Cal H}
\def\cM{\Cal M}

\def\Spec{\operatorname{Spec}}

\def\PGL{\operatorname{PGL}}
\def\tM{\tilde{M}}
\def\lc{\operatorname{lc}}

\def\cE{\Cal E}

\def\bS{\Bbb S}
\def\bB{\Bbb B}

\def\Ker{\operatorname{Ker}}
\def\Coker{\operatorname{Coker}}
\def\Im{\operatorname{Im
}}

\def\ord{\operatorname{ord}}

\def\cU{\Cal U}

\def\oZ{\overline{Z}}

\def\oR{\overline{R}}

\def\Cal{\mathcal}

\def\Pic{\operatorname{Pic}}
\def\val{\operatorname{val}}
\def\oN{\overline{N}}
\def\os7p{\oS_7'}
\def\os7{\oS_7}
\def\on6{\oN_6}
\def\n6{\oN_6}

\def\irr{\operatorname{irr}}

\def\cG{\Cal G}

\DeclareMathOperator{\ch}{char}

\newtheoremstyle{mystyle}{}{}{\itshape}{}{\scshape}{.}{ }{}
\theoremstyle{mystyle}

\swapnumbers
\newtheorem{Theorem}{Theorem}[section]

\newtheorem{Proposition}[Theorem]{Proposition}
\newtheorem{Lemma}[Theorem]{Lemma}

\newtheorem{Corollary}[Theorem]{Corollary}
\newtheorem{Claim}[Theorem]{Claim}

\newtheoremstyle{myreview}{}{}{}{}{\scshape}{.}{ }{}
\theoremstyle{myreview}

\newtheorem{Definition}[Theorem]{Definition}

\newtheorem{Remark}[Theorem]{Remark}
\newtheorem{Remarks}[Theorem]{Remarks}
\newtheorem{Notation}[Theorem]{Notation}
\newtheorem{Review}[Theorem]{}

\newcounter{et}[Theorem]
\def\cooltag{\tag{\arabic{section}.\arabic{Theorem}.\arabic{et}}\addtocounter{et}{1}}

\begin{document}
\title[Compact moduli of del Pezzo surfaces]{Stable pair, tropical, and log canonical\\
compact moduli of del Pezzo surfaces}
\author{Paul Hacking, Sean Keel, and Jenia Tevelev}

\begin{abstract}
We give a functorial normal crossing compactification of the moduli of smooth
cubic surfaces
entirely analogous to the Grothendieck-Knudsen
compactification $M_{0,n} \subset \oM_{0,n}$.
\end{abstract}

\maketitle

\section{Introduction and Statement of Results}\label{intro}
Throughout we work over an algebraically closed field $k$.
Let $Y^n$ be the 
moduli space 
of smooth marked del Pezzo surfaces $S$ of degree~$9-n$. We begin by observing
that each such surface comes with a natural boundary:
The marking of $S$ induces a labelling of its $(-1)$-curves, and we can
take their union $B \subset S$. This gives a natural way of 
compactifying the space:
Let 
$Y^n_{\times} \subset Y^n$ denotes the open locus where $B$ has normal
crossings. 
For $\ch k = 0$, let $\ocM$ denote the Koll\'ar--Shepherd-Barron--Alexeev 
stack of \nobreak{stable} surfaces with boundary \cite{KSB,Alexeev}
and $\oM$ its coarse moduli space. $\ocM$ is complete, so we can compactify by taking
the closure of $Y^n_{\times}$ in $\oM$. This closure turns out to have very
nice properties:

\begin{Theorem} \label{KSBcompact}
Assume $n \le 5$ or $n=6$ and $\ch k \neq 2$.
There is a compactification $Y^n \subset \oY^n_{\sts}$ and a 
family $p \colon (\cS,\cB) \rightarrow \oY^n_{\sts}$ of stable surfaces with boundary 
extending the universal family of smooth del Pezzo surfaces 
with normal crossing boundary over $Y^n_{\times}$. 
$\oY^n_{\sts}$ is a smooth projective variety, and 
$\oY^n_{\sts} \setminus Y^n$ is a union of smooth divisors with normal crossings.
The fibers of $p$ have stable toric singularities.
If $\ch k = 0$, $p$ defines a closed embedding $\oY^n_{\sts} \subset \oM$.
\end{Theorem}

(Recall that a pair $(S,B)$ has \emph{stable toric singularities} if it is 
locally isomorphic to a stable toric variety \cite{Alexeev2} 
together with its toric boundary. 
A stable toric variety is obtained by glueing normal toric varieties along toric boundary divisors, 
and its toric boundary is the union of the remaining boundary divisors.)

We call the fibres of $p \colon (\cS,\cB) \to \oY^n_{\sts}$
stable del Pezzo pairs. We will show that each
fibre $(S,B)$ comes with a canonical embedding in
a stable toric pair $(X,B)$, with $S \subset X$ transverse
to the toric strata (with the boundary of $S$ the
restriction of the toric boundary), see Theorems~\ref{univfams} and \ref{univfamE_6}.

Recall that a smooth variety $Y$ is \emph{log minimal} if for some (equivalently, 
any) smooth compactification $\tilde{Y}$ with normal crossing boundary $B$ the linear system 
$|N(K_{\tilde{Y}}+B)|$ defines an embedding on $Y$ for $N$ sufficiently divisible. 
Such a $Y$ is expected to come with a natural log canonical compactification $\oY_{\lc}$, see \cite[2.5]{BCHM}. 

\begin{Theorem} \label{lccompact}
The variety $Y^n$ is log minimal. For $n \leq 6$ or $n=7$ and $\ch k \neq 2$, its log canonical
compactification $\oY^n_{\lc}$ is smooth and $\oY^n_{\lc}\setminus Y^n$ is a union of smooth divisors with normal crossings.

Let $\pi:\,Y^{n+1}\to Y^n$ be the natural functorial morphism.
Assume $n \leq 5$ or $n=6$ and $\ch k \neq 2$. Then $\pi$ extends to a
morphism $\pi:\,\oY^{n+1}_{\lc} \to \oY^n_{\lc}$
and we have a commutative diagram
$$
\begin{CD}
\cS @>>> \oY^{n+1}_{\lc}\\
@VpVV                 @VV\pi V \\
\oY^n_{\sts}  @>>> \oY^n_{\lc},
\end{CD}
$$
where the horizontal arrows are isomorphisms for $n \le 5$ and log crepant birational morphisms for $n=6$.
Moreover, $\oY^6_{\sts} \to \oY^6_{\lc}$ is a composition of blowups with smooth centers.
In particular $\oY^n_{\sts}$ is a log minimal
model for $Y^n$ and its universal family is a log minimal
model for $Y^{n+1}$.
\end{Theorem}

Note the close analogy with $M_{0,n} \subset \oM_{0,n}$:
This is the log canonical compactification as well as the
closure in the moduli space of pointed stable curves, and
$\oM_{0,n+1} \to \oM_{0,n}$ is the universal family.
The two cases are instances of
a single construction which we explain below.

\begin{Remarks}\noindent
\begin{enumerate}
\item In general boundary strata of irreducible components of $\ocM$ have arbitrary singularities, 
even for (degenerations of) lines in $\bP^2$. Already for $9$ lines the compactification
has non log canonical boundary. See \cite{ChQII}. The del Pezzo
case is thus highly exceptional. 
\item The Weyl group $W(E_n)$ acts on $Y^n$ by changing the marking. Our constructions are
equivariant, so we can quotient by $W(E_n)$ to get moduli spaces of unmarked
del Pezzo surfaces, and orbifold versions of the above theorems.
\item We expect that $\oY^n_{\sts}$ is a connected component of $\oM$ and is a fine moduli space.
We have verified this for $n=5$ by explicit deformation theory calculations.
\item Following the philosophy of \cite{H}, and \cite{Alexeev2} 
we believe that $\oY^n_{\lc}$ ($n \geq 6$) itself has a natural functorial
meaning -- namely the connected component of $\ocM$ corresponding to 
$(S,eB)$ where we choose $e$ minimal so that $K_S + eB > 0$. 
\item The universal family $(\cS,\cB) \to \oY^n_{\sts}$ is highly non-trivial.
E.g. in the first non-constant case 
$(\cS,\cB) \to \oY^5_{\sts}$, the base is a copy of
$\oM_{0,5}$. $\cB$ has $16$ irreducible
components, each a copy of
the universal family $\oM_{0,6} \to \oM_{0,5}$. 
The simplest degenerate
fibre is a surface with $6$ (smooth) components (across which the $16$ $(-1)$-curves are distributed). 
\item The tautological map $\pi:Y_{n+1} \to Y_n$
(given by blowing down a $(-1)$-curve) induces a map of tropical varieties 
$\trop(\pi):\trop Y_{n+1} \to \trop Y_n$ --- a map of fans induced by a linear map of ambient
real vector spaces (we review tropical geometry below and give a new, and self contained,
treatment of the theory in \ref{Coxco}). 
$\trop(\pi)$ is canonically associated with the inclusion of root systems
$E_n \subset E_{n+1}$ as we explain below. 
Our compactifications, and the universal families, together
with all the combinatorics of the boundary and fibres, can be read off from $\pi$ and
this map of fans: together they determine 
the functorial compactification, and the functor that it represents. 
In particular we get the complete picture of possible degenerations, 
without considering degenerations (or even surfaces) at all. 
\end{enumerate}
\end{Remarks}

We will identify $\oY^n_{\lc}$ with known spaces:
$\oY^5_{\lc}\simeq\oM_{0,5}$,
$\oY^6_{\lc}$ is isomorphic to the Naruki space $\oY^6$ of cubic
surfaces~\cite{Naruki},
and $\oY^7_{\lc}$ is isomorphic to the Sekiguchi space~$\oY^7$ of
del Pezzo surfaces of degree~$2$~\cite{Seki1,Sekiguchi}.
Sekiguchi constructed
$\oM_{0,n}$, $\oY^6$, and $\oY^7$ in a unified way,
using ``cross-ratio maps'' associated to root subsystems of type $D_4$
in root systems $D_n$, $E_6$, and $E_7$. 

\begin{Definition}\label{typeIdef}
Let $S$ be a marked del Pezzo surface.
By a {\em KSBA cross-ratio} we mean the cross-ratio of $4$ points on
a $(-1)$-curve $L$ cut out by $(-1)$-curves $L_1,\ldots,L_4$
(assuming that these points are distinct and that $L_i \cdot L = 1$).
If, in addition, $L_i \cdot L_j = 0$ for $1\le i < j\le 4$ then
we have a {\em cross-ratio of type $I$}.
\end{Definition}

We show that KSBA cross-ratios of type $I$ are the same as
Sekiguchi's $D_4$ cross-ratios.
Sekiguchi defines $\oY^n$ as the closure of the image of $Y^n$
in the product $(\bP^1)^N$ given by all cross-ratios of type $I$.
Naruki proved that $\oY^6$ is smooth
and its boundary has simple normal crossings. Sekiguchi
conjectured the analogous result for $n=7$. We prove:

\begin{Theorem}[Sekiguchi Conjecture]\label{sekconj}
For $n\le 6$ or $n=7$ and $\ch k \neq 2$, $\oY^n\simeq\oY^n_{\lc}$. In particular,
$\oY^7$ is smooth  and the boundary $\oY^7\setminus Y^7$ has
simple normal crossings.
\end{Theorem}

It is rather unnatural to take only cross-ratios of type $I$:

\begin{Theorem}\label{crKSB}
For $n \le 5$ or $n=6$ and $\ch k \neq 2$, the closure of the image of $Y^n_\times$ in the
product $(\bP^1)^M$ given by all KSBA cross-ratios
is isomorphic to $\oY^n_{\sts}$.
\end{Theorem}

We obtain $\oY^n_{\lc}$, $\oY^n_{\sts}$, and its universal
family $(\cS,\cB)$, together with an ambient family
of stable toric pairs $(\cX,\cB)$,
canonically from the interior~$Y^n$ by applying elementary
ideas from tropical algebraic geometry~\cite{Tevelev}.
The same construction
applied to $M_{0,n}$ yields $\oM_{0,n}$. Next we explain
the procedure.

Let $Y$ be a variety defined over an algebraically closed field $k$.
The group of units $M:=\cO^*(Y)/k^*$ is a free Abelian group of finite
rank.
Let $N=M^*$. The tropical variety $\cA$ of $Y$ is an intrinsic subset
of $N_\bQ$.
We refer to \cite{Tevelev,EKL,SS} for its definition.
In \ref{Coxco} we will prove the following new characterization
of~$\cA$:
\begin{equation}\label{newcharact}
\cA(Y)=\bigl\{[v]\,\bigm|\,v\ \hbox{\rm is a discrete valuation of}\
k(Y)\bigr\}\subset N_\bQ,\cooltag
\end{equation}
where $[v](u)=v(u)$ for any unit $u\in M$. Note
$[v] =0$ if $v$ has center on $Y$, thus by
\eqref{newcharact} $\cA(Y)$ ``sees'' the boundary divisors
from all compactifications of $Y$. We show it can be
computed from a single normal crossing compactification,
see Cor.~\ref{constcasencc}.

\begin{Definition}[\cite{Tevelev}] 
Let $T:=\Hom(M,\bG_m)$ be the {\em intrinsic} algebraic torus.
Choosing a splitting of the exact sequence
$$0 \rightarrow k^{*} \rightarrow \cO^*(Y) \rightarrow M \rightarrow 0$$
defines an evaluation map $Y \to T$ (any two such are related by a translation by an element of $T$).
Suppose that $Y\to T$ is a closed embedding
(in which case $Y$ is called {\em very affine}).
Let $\cF\subset N_\bQ$ be a fan,
$X(\cF)$ the corresponding toric variety, and $\oY(\cF)$
the closure of $Y$ in $X(\cF)$. We call $\cF$ and $\oY(\cF)$ {\em tropical}
if $\oY(\cF)$ is complete and the multiplication map
$\oY(\cF) \times T \to X(\cF)$ is flat and surjective.
\end{Definition}

\begin{Theorem}[\cite{Tevelev}] Tropical fans exist, and if
$\cF$ is tropical then $|\cF| = \cA$. Any refinement
of a tropical fan is tropical.
\end{Theorem}

It is natural to wonder if there is a canonical fan structure
on $\cA$. The answer is no in general, but there is a natural
Mori-theoretic sufficient condition:

\begin{Definition} \label{defhubsch}
We say that a very affine variety $Y$ is {\em Sch\"on}\/ if the multiplication map of one
(and hence of any \cite[1.4]{Tevelev}) tropical compactification
is smooth. A Sch\"on variety $Y$
is called {\em H\"ubsch}\/ if it is log minimal  and its log canonical
compactification $\oY_{\lc}$
is tropical, i.e., $\oY_{\lc}\simeq\oY(\cF)$ for some fan $\cF$ (called the {\em log canonical fan}).
\end{Definition}

In fact we prove that a Sch\"on subvariety of an algebraic torus is either log minimal or preserved by a subtorus, see \ref{ReidQ}.
This answers a question of M.~Reid.

\begin{Theorem} \label{minfan}
Let $Y$ be a H\"ubsch very affine variety.
Then every fan structure on $\cA$ is tropical, and
a refinement of the log canonical fan.
\end{Theorem}

Thus in the H\"ubsch case
$\cA \subset N_{\bQ}$ has a canonical fan structure, which
yields the log canonical compactification embedded in a
canonical ambient toric variety, and transverse to the
toric boundary. The construction works in the other direction as
well --- $\cA$ together with its canonical fan structure
can be read off from the boundary
stratification of the log canonical compactification, see
\S 2.

Here are some examples of H\"ubsch varieties:

\begin{Theorem} $M_{0,n}$ and $Y^n$ (for $n \leq 6$ or $n=7$ and $\ch k \neq 2$) are H\"ubsch.
The corresponding tropical ($=$ log canonical) compactifications are
$\oM_{0,n}$ and $\oY^n_{\lc}$.
\end{Theorem}

$M_{0,n}$ is isomorphic to the complement of the braid arrangement of hyperplanes in $\bP^{n-3}$. More generally, recall that a hyperplane arrangement in $\bP^m$ is \emph{connected} if 
the subgroup of $\PGL_{m+1}$ preserving the arrangment is finite (or equivalently, if
the complement is log minimal).
By \cite{HackingKeelTevelev}:

\begin{Theorem}\label{hpch} Complements of arbitrary connected hyperplane arrangments
are H\"ubsch and their tropical ($=$ log canonical)
compactifications are Kapranov's visible contours~\cite{CQ1}.
\end{Theorem}

We write $Y(E_n) := Y^n$ and $Y(D_n) := M_{0,n}$. We
will treat the two cases in a unified way.
$\Delta$ will indicate
a root system of type $D_n$ or $E_n$.
Next we describe the log canonical fans for $Y(\Delta)$
using combinatorics of root subsystems introduced by
Sekiguchi~\cite{Sekiguchi}.
Let
$W(\Delta)$ be the Weyl group.
Let $\Delta_+\subset\Delta$ be the positive roots.
$W(D_n)$ acts on $M_{0,n}$ through its quotient $S_n$
and $W(E_n)$ acts on $Y^n$ by changing markings, see \cite{DO} or
\ref{intrtorus} for details.

\begin{Theorem}
Let $\bZ^{\Delta_+}$ be the lattice with basis vectors $[\alpha]$ for
$\alpha\in\Delta_+$.
Let
\begin{equation}\label{funnyMdef}
M(\Delta)=\left\{\sum n_{\alpha}[\alpha]\,\Big|\,
\sum n_{\alpha}\alpha^2=0\right\}.\cooltag
\end{equation}
Then $M(\Delta)$ is an irreducible $W(\Delta)$-module
of rank equal to the number of roots in $\Delta_+$ with  three-legged
support.
We have an isomorphism of $W(\Delta)$-modules
$$\cO^*(Y(\Delta))/k^*=M(\Delta).$$
\end{Theorem}

\begin{Review}\label{defrn}
Let $\psi:\,\bZ^{\Delta_+}\to N(\Delta)$ be the dual map.
For any root subsystem $J\subset\Delta$, let
$$\psi(J):=\sum_{\alpha\in J\cap\Delta_+}\psi(\alpha)$$
and let $\Cal J$ be the set of root subsystems of $\Delta$ of type $J$.
For example, $\cA_1$ is the set of $A_1$'s in $\Delta$.
Let $\cR(\Delta)$ be a simplicial complex defined as follows.
As a set,
$$\cR(D_{2n+1})=\cD_2\ \sqcup\ \cD_3\ \sqcup\ \ldots\ \sqcup \cD_n;$$
$$\cR(D_{2n})=\cD_2\ \sqcup\ \cD_3\ \sqcup\ \ldots
\ \sqcup\ \cD_{n-1}\ \sqcup\ (\cD_{n}\times\cD_n);$$
$$\cR(E_6)=\cA_1\ \sqcup\ (\cA_2\times \cA_2\times \cA_2);$$
$$\cR(E_7)=\cA_1\ \sqcup\ \cA_2\ \sqcup\ (\cA_3\times \cA_3)\ \sqcup\
\cA_7.$$

Subsystems $\Theta_1,\ldots,\Theta_k\in\cR$ form a simplex if
and only if
$\Theta_i\perp\Theta_j$, or
$\Theta_i\subset\Theta_j$,
or $\Theta_j\subset\Theta_i$ for any pair $i$, $j$ with the following exception:
For $\Delta=E_7$, we exclude ``Fano simplices''
formed by $7$-tuples of pairwise orthogonal~$A_1$'s.

Let $\cF(\Delta)\subset N(\Delta)_\bQ$ be the collection
of cones determined by the rays $\psi(\Theta)$
for each $\Theta\in\cR(\Delta)$, where rays span a cone in $\cF(\Delta)$
iff the corresponding subsystems form a simplex in~$\cR(\Delta)$.
\end{Review}

\begin{Definition}
We say that a collection $\{\sigma_i\}$ of convex subsets of
$\bR^n$ is {\em convexly disjoint}\/
if any convex subset of $\bigcup_i\sigma_i$ is contained in one of them.
\end{Definition}

For $N$ a free abelian group of finite rank and $\sigma$ a cone in $N \otimes \bQ$, we say $\sigma$ is \emph{strictly simplicial}
if it is generated by a subset of a basis of $N$. We say a fan $\cF$ is strictly simplicial if each cone of $\cF$ is so. Equivalently, the associated
toric variety $X(\cF)$ is smooth.

\begin{Theorem} ($n \le 6$ or $n=7$ and $\ch k \neq 2$ for $E_n$)
$\cF(\Delta)$ is the log canonical fan of $Y(\Delta)$. This
fan is convexly disjoint and strictly simplicial.
\end{Theorem}

In what follows we frequently use

\begin{Theorem}[{\cite[3.1]{Tevelev}}]\label{dominanta}
A dominant map of very affine varieties $Y'\to Y$ induces a
surjective homomorphism of intrinsic tori $T_{Y'}\to T_Y$ and a
surjective map
of tropical varieties  $\cA(Y') \to\cA(Y)$.
\end{Theorem}

Note that if $Y$ is H\"ubsch and the log canonical fan is
strictly convex, then $\cA(Y') \to \cA(Y)$ is a map of
fans, for any fan structure on $\cA(Y')$
(and the log canonical fan structure on $\cA(Y)$).

\begin{Theorem} \label{univfams}
Let $(\Delta_{n+1},\Delta_n)$ be either
$(D_{n+1},D_n)$ or $(E_{n+1},E_n)$ for $n \le 5$ or $n=6$ and $\ch k \neq 2$.
The natural morphism $Y(\Delta_{n+1}) \to Y(\Delta_n)$ induces
a surjective map of fans $\cF(\Delta_{n+1}) \to \cF(\Delta_{n})$.
This induces a commutative
diagram:
$$
\begin{CD}
\oY_{\lc}(\Delta_{n+1}) @>>> X(\cF(\Delta_{n+1})) \\
@VpVV  @V{\pi}VV \\
\oY_{\lc}(\Delta_n) @>>> X(\cF(\Delta_n))
\end{CD}
$$
where the horizontal maps are closed embeddings.
For $\Delta_n = D_n$, or $E_n$, $n \leq 5$,
$\pi$ is a flat family of stable toric varieties,
and $p$ is a flat
family of stable pairs which induces an isomorphism $\oY_{\lc}(\Delta_n) \stackrel{\sim}{\rightarrow} \oM_{0,n}$ for $D_n$
and a closed embedding of $\oY_{\lc}(\Delta_n)$ into the moduli space $\oM$
of stable surfaces with boundary for $E_n$ (assuming $\ch k =0$).
Each fibre of $p$ is transverse to the toric strata
of the corresponding fibre of $\pi$.
\end{Theorem}

If $\Delta_n=E_6$ then $\pi$ is not flat.
However there is a canonical combinatorial
procedure for flattening a toric morphism.
We write $\cF^n$ for $\cF(E_n)$. 

\begin{Theorem} \label{univfamE_6}
Assume $\ch k \neq 2$. 
There is a unique minimal refinement $\tcF^6$ of $\cF^6$ such that the
closure $X(\tcF^7)$ of the torus $T(Y^7)$ in the fibre product $X(\cF^7) \times_{X(\cF^6)} X(\tcF^6)$
is flat over $X(\tcF^6)$. 
$X(\tcF^7)$ is a normal toric variety, let $\tcF^7 \subset N(E_7)$ be the corresponding fan.
There is an induced commutative diagram
$$
\begin{CD}
\oY(\tcF^{7})  @>{i}>> X(\tcF^{7}) @>>> X(\cF^{7}) \\
@V\tp VV       @V{\tpi}VV @V{\pi}VV \\
\oY(\tcF^6) @>{j}>> X(\tcF^{6}) @>>>   X(\cF^6)
\end{CD}
$$
where $i$ and $j$ are closed embeddings.
$\tp$ and $\tpi$ are flat with reduced fibers.
$X(\tcF^6)$ and $\oY(\tcF^6)$ are smooth with simple normal crossing boundary.
Away from Eckhart points, each fibre of $\tp$ is transverse to the toric strata
of the corresponding fibre of $\tpi$. 
Blowing up the Eckhart points of $\tp$ yields a family of stable
surfaces with boundary which induces a closed embedding of $\oY(\tcF^6)$ in $\oM$ for $\ch k = 0$.
\end{Theorem}
Recall that an \emph{Eckhart point} on a smooth cubic surface $S$ is an ordinary triple point of the union $B$ of the $(-1)$-curves.
The Eckhart points of the family $\tilde{p}$ are by definition the closure of the locus of Eckhart points on the smooth fibres. 
If $s$ is an Eckhart point and $(S,B)$ is the fiber through $s$, then $S$ is smooth near $s$ and $B$ has an ordinary triple point at $s$. 
See \ref{eckpoints}.

$Y^n$ has received considerable attention from the point of view of hyperplane 
arrangements on ball quotients,
by Kond\=o, Dolgachev, Heckman, Looijenga and others. See \cite{Looijenga07} and the references
there. In particular $Y^n$ has a Baily-Borel compactification, and 
Looijenga gives desingularisations which are naturally determined by the space together with its arrangement of hyperplanes. 
$\oY^n_{\lc}$, $n \leq 7$ can be constructed by Looijenga's procedure. 
None of these authors (to our knowledge) consider the modular meaning of the compactification (i.e. whether the boundary
parameterizes any sort of geometric object). 

\begin{Review}\textsc{Outline of the proof.}\label{outlineoftheproof}
Let $Y=Y(\Delta)$. 
By Cor.~\ref{constcasencc}, to
find $\cA(Y)$, expressed as the union of a collection of cones $\cF$ (which may or
may not form a fan), it is
enough to find an orbifold normal crossing compactification
and describe its boundary.
We use the compactification $\oY_{\lc}$ (which we show has orbifold normal crossings).
We identify $\oY_{\lc}$ for $D_n$ with $\oM_{0,n}$ and
$\oY_{\lc}$ for $E_6$ with the Naruki space $\oY^6$,
see Cor.~\ref{n6lcm}.
For $E_7$ there is a straightforward way to compute~$\oY_{\lc}$:
The anti-canonical map for a smooth del Pezzo of degree $2$ is
a double cover of $\bP^2$ branched over a smooth quartic. This
expresses $Y$ as a finite Galois cover
$Y \to M_3 \setminus H$ of the moduli space of smooth non-hyperelliptic
curves of genus~$3$. The log canonical compactification of $Y$ is
then the normalization in this field extension of the log
canonical compactification of $M_3 \setminus H$.
The latter turns out to be not $\oM_3$ but
the blowup $\tM_3 \to \oM_3$ along the locus of curves
whose generic point corresponds to two genus~$1$ curves glued
at two points.

Using root systems we can
describe $\cF$ and prove that it is a strictly simplicial fan and
convexly disjoint by purely combinatorial means. 
This in particular implies that $\oY_{\lc}$ is smooth, the boundary has simple normal crossings,
and that $\oY_{\lc}$ embeds in the toric variety with the fan $\cF$.
This is carried out in sections 4--8.
Our main technical tool is a
(it seems to us rather miraculous) larger more symmetric
fan $\cG$ coming from the real points of $Y$, see \ref{realcomp}.
The dominant tautological map
$Y^{n+1} \to Y^n$ induces a surjection of sets
$\cA(Y^{n+1}) \to \cA(Y^n)$, which is a map of
fans  $\cF^{n+1}\to \cF^n$, as $\cF^n$ is convexly disjoint.
Therefore, there is an induced morphism of toric varieties, and thus
of the log canonical compactifications $\oY^{n+1}_{\lc} \to \oY^n_{\lc}$.
For $n \leq 5$ this is the universal family of stable
surfaces. For $n =6$ the map of toric varieties is not
flat, but the canonical procedure for flattening such a map
leads to the universal family. See~\ref{sds}.
\end{Review}

\begin{Review}\textsc{Acknowledgements.}
Daniel Allcock helped us a great deal with the branch
cover constructions of \ref{branchcover}. We
thank A.-M. Castravet, I. Dolgachev, B. Fantechi, G. Farkas, S.
Grushevskii, B. Hassett,
G. Heckmann, R. Heitmann, M. Luxton, M.~Olsson, Z.~Qu, M.~Reid, B.~Sturmfels,
and L.~Williams
for many helpful discussions. We are particularly
grateful to Professor Jiro Sekiguchi for sending us copies of his
(quite remarkable) papers,
from which we got a great deal of combinatorial inspiration.
The first author was partially supported by NSF grant DMS-0650052,
the second author by NSF grant DMS-0353994, and the
third author by NSF grant DMS-0701191 and a Sloan research fellowship. 
\end{Review}

\tableofcontents

\section{Geometric Tropicalization and Cox Coordinates} \label{geotrop}\label{Coxco}

Let $K$ be the field of Puiseux series over $k$ and $R$ its valuation ring.
The valuation $\deg$ induces an isomorphism $\,K^*/R^*\simeq\bQ$.
Let $Y$ be an affine variety over~$K$,
$M=\cO^*(Y)/K^*$, and $\tM := \cO^*(Y)/R^*$.
We have an exact
sequence
$$
0 \to K^*/R^* \to \tM \to M \to 0.
$$

\begin{Definition} \label{defBG}
For a field extension $K \subset L$, a map
$f: \Spec(L) \to Y$, and a
valuation $v: L^* \to \bQ$ such that $v(R^*)=0$ and $v(R)\ge0$,
let $[f,v] \in \tN_\bQ:= \Hom(\tM,\bQ)$ be the
element induced by the composition
$\cO^*(Y)\arrow^{f^*} L^* \arrow^v \bQ$.
Let $\tcA \subset \tN_{\bQ}$ be the union of
such $[f,v]$.
The {\em Bieri-Groves set}\/ \cite{BG,EKL} $\cA \subset \tcA$ is
 the union of $[f,v]$
such that the restriction of $v$ to $R$ is
the standard valuation.
\end{Definition}

\begin{Definition} \label{bfan} Let $\oY \to \Spec R$ be a normal $\bQ$-factorial
variety of finite type over $\Spec R$. Let $Y \subset \oY$
be an open subset such that the boundary
$\partial \oY := \oY \setminus Y$ is divisorial. Let
$\cF(\partial \oY)$ denote the collection of cones in
$\tN(Y)_{\bQ}$ defined as follows: For each irreducible
boundary divisor $D \subset \partial \oY$
let $\ord_D$ be
the corresponding
valuation of the function field $K(Y)$. This is trivial
on $R^*$ (since elements of $R^*$ pull back to units on $\oY$)
and so induces an element $[D] \in \tN_{\bQ}$.
For any collection $\sigma=\{D_1,\ldots,D_k\}\subset\partial \oY$ such
that
$\bigcap D_i\ne\emptyset$,
let $\cF(\sigma)$ be the cone in $\tN_\bQ$ spanned by
$[D_1],\ldots,[D_k]$. 

If $\oZ$ is a normal $\bQ$-factorial variety over $k$ and $Z \subset \oZ$ an open subset
with divisorial boundary $\partial \oZ := \oZ \setminus Z$, we define a collection of cones
$\cF(\partial \oZ) \subset N(Z)_{\bQ}$, where $N(Z)=\Hom(M(Z),\bZ)$ and $M(Z)=\cO^{*}(Z)/k^{*}$,
in the same way as above. That is, the rays of $\cF(\partial \oZ)$ are generated by the valuations $\ord_D$
for $D$ a component of $\partial \oZ$, and a collection of rays spans a cone of $\cF(\partial \oZ)$ if the corresponding divisors intersect.
\end{Definition}

A version of the following theorem 
was proved by W.~Gubler \cite[7.10]{G}.
We observed the result independently in joint work with Z.~Qu \cite{Quthesis}.

\begin{Theorem} \label{tropforncc}
Suppose we have a diagram
$$\begin{CD}
\cY @>>> \ocY \\
@VVV  @VVV \\
\Spec K @>>> \Spec R
\end{CD}$$
where $\cY \subset \ocY$ is an open substack of a
smooth Deligne-Mumford stack, proper over $\Spec(R)$,
such that the boundary $\partial \ocY:=\ocY \setminus \cY$ is a divisor with
simple normal crossings. Let
$$\begin{CD}
Y @>>> \oY \\
@VVV  @VVV \\
\Spec K @>>> \Spec R
\end{CD}$$
be the induced diagram of coarse moduli spaces --- so in particular
the boundary has orbifold simple normal crossings.

Then $\tcA$ is the union of the collection of cones $\cF(\partial \oY)$.
In particular, $\tcA$ is the underlying set of a fan, and the
Bieri-Groves set $\cA$
is the underlying set of the fibre of the map of fans induced by
$\tN_{\bQ} \to \bQ=\Hom(K^*/R^*,\bQ)$.
\end{Theorem}

\begin{proof}
Let $\cG \subset \tN_\bQ$ be the union of cones in the statement.
Since $R^* \subset \cO^*(\oY)$ and $R\subset\cO(\oY)$,
divisors of $\oY$ have trivial
valuation on $R^*$ and nonnegative valuation on~$R$. Therefore, $[D] \in \tcA$.
Moreover, let $\sigma=\{\cD_1,\ldots,\cD_k\}\subset\partial \ocY$ with $\cap
D_i\ne\emptyset$.
Then any linear combination
$\sum n_i[D_i]$ with nonnegative rational coefficients belongs to
$\tcA$. Indeed, the exceptional divisors of weighted blowups of the regular sequence 
$\sigma$ give arbitrary nonnegative combinations (up to scalar multiples). 
Thus $\cG \subset \tcA$.

Let $S \subset L$ be the valuation ring for $L,v$ as in
Definition~\ref{defBG}. As $v|_{R}$ is non-negative, we have the morphism $\Spec(S)
\to \Spec(R)$,
and so by properness  there is a unique
extension $f: \Spec(S) \to \oY$.
Let $\sigma$ be the collection of all boundary divisors of $\oY$
that contain the image of the closed point of $\Spec(S)$. Let $u \in
\cO^*(Y)$
be a unit, and let $B \subset \partial \oY$ be its locus of poles, the
union of boundary divisors with negative valuation on $u$.
Then $u \in \cO(\oY \setminus B)$, so if
$[D](u) \geq 0$ for all $D \in \sigma$, $f^{-1}(B) = \emptyset$ so
$f^*(u) \in S$ and
$[v,f](u) = v(f^*(u)) \geq 0$.
Thus $[v,f] \in \cF(\sigma)$. \end{proof}

\begin{Corollary}\label{constcasencc} Let
$\cZ \subset \ocZ$ be an open substack of a
Deligne-Mumford stack proper and smooth over $\Spec k$,
with simple normal crossing
boundary, and $Z \subset \oZ$ the corresponding orbifold
simple normal crossing compactification of coarse moduli
spaces.  Then $\cA(Z)$ is the underlying set of a fan,
and is the union of the collection of cones $\cF(\partial \oZ)$.
\end{Corollary}

\begin{proof}
If $\cY = \cZ \times_k K$ for a $k$-variety $\cY'$ then $\ocY:= \ocZ \times_k
R$
satisfies all conditions of Theorem~\ref{tropforncc}.
It follows that $\tcA_{Y}$ is the product $\cA_{Z}\times \bQ_{\geq 0}$.
The statement follows easily.
\end{proof}

In what follows we frequently use

\begin{Theorem}[{\cite[3.1]{Tevelev}}]\label{dominanta}
A dominant map of very affine varieties $Y'\to Y$ induces a
surjective homomorphism of intrinsic tori $T_{Y'}\to T_Y$ and a
surjective map
of tropical varieties  $\cA(Y') \to\cA(Y)$.
\end{Theorem}

Note if $Y$ is H\"ubsch, and the log canonical fan is
strictly convex, then $\cA(Y') \to \cA(Y)$ is a map of
fans, for any fan structure on $\cA(Y')$
(and the log canonical fan structure on $\cA(Y)$).

\begin{proof}[Proof of \eqref{newcharact}]
Since discrete valuations lift to field extensions it is clear that the
set defined in \eqref{newcharact} satisfies the analog of Th. ~\ref{dominanta}.
Now by de Jong's theorem \cite[4.1]{Jong} we may assume $Y$ is smooth and admits
a normal crossing compactification. Now the result follows
from the proofs of Th.~\ref{tropforncc} and Cor.~\ref{constcasencc}.
\end{proof}

\begin{proof}[Proof of Theorem~\ref{minfan}]
Let $k=\dim Y$. Recall that top-dimensional cones of $\cF$
are $k$-dimensional.
By \cite[2.5]{Tevelev}, any refinement of
a tropical fan is tropical, so it is enough to show that an
arbitrary fan $\cF''$ with $|\cF''| = \cA_{Y}$ refines
the log canonical fan $\cF$.
Suppose on the contrary that $\cF''$ does not refine $\cF$. Then
there exists a $(k-1)$-dimensional cone $\alpha \in\cF$ which meets the
interior of a $k$-dimensional cone $\sigma'' \in\cF''$.
Let $\cF'$ be a strictly
simplicial common refinement of $\cF$ and $\cF''$ (with the same support).
There exists a $(k-1)$-dimensional cone $\alpha' \in\cF'$,
$\alpha' \subset \alpha \cap \sigma''$, and meeting
the interiors of each. Let $X'\arrow^p X''$, $X'\arrow^\pi X$
be the corresponding proper birational toric morphisms.
Let $Z$, $Z'$, $Z''$ be closed toric
strata corresponding to $\alpha$, $\alpha'$, $\sigma''$. Then
\begin{itemize}
\item $Z''=T''$ is a quotient torus of $T$.
\item $p:\,Z' \to T''$ is a surjective proper
toric map, thus $Z' = \bP^1\times T''$.
\item $\pi: Z' \to Z$ is birational.
\end{itemize}

By \cite[1.4]{Tevelev}, the structure map $\oY'\times T\to X'$ is smooth and surjective,
and therefore $G = \oY' \cap Z' \subset \bP^1\times T''$ is $1$-dimensional, smooth, reduced, and proper.
Since $G$ is proper and $T''$ is affine,
it follows that $G=\bP^1 \times z$ for a reduced $0$-dimensional subscheme $z \subset T''$.
So $K_G + B_G$ is trivial. $G$ is
also a closed stratum of $\oY'$ and thus (since $\oY'$
has smooth structure map)
$(K_{\oY'} + B)|_G = K_G + B_G$ is trivial. Thus $G$
is contracted by $\pi$, since $\pi: \oY' \to \oY$ is
the log canonical model, and log crepant by \cite[1.4]{Tevelev}. By
equivariance of the map $\pi:\,Z' \to Z$, it contracts all fibres of
$Z' \to Z''$, i.e., it is not birational, a contradiction.
\end{proof}

Suppose that $\oY$ is a compactification of
the smooth variety $Y$ with simple normal crossing boundary $\partial \oY$.
We assume until the end of this section that $\Pic(\oY)$ is a lattice generated by the classes of
boundary divisors.
Let $A:=\bZ^{\partial \oY}$ be the free Abelian group with a basis
given by irreducible boundary divisors.
We have a canonical exact sequence
$0 \to \Pic(\oY)^\vee  \to
A \to N_Y \to 0$,
where $N_Y$ is dual to $M_Y=\cO^*(Y)/k^*$.
Let $T_{\Pic}$, $\bG_m^{\partial\oY}$, $T_Y$ be the algebraic tori with characters
$\Pic(\oY)$, $A^\vee$, $M_Y$.
We have an exact sequence
$\{e\} \to T_{\Pic}\to \bG_m^{\partial\oY}\to T_Y\to\{e\}$.

Let $\bA^{\partial \oY}$ be the affine space and let
$\cU \subset \bA^{\partial \oY}$ be the open $\bG_m^{\partial\oY}$-toric subvariety
defined as follows: a collection of coordinates vanish
simultaneously at some point of~$\cU$ iff the intersection
of the corresponding divisors of $\oY$ is non-empty.
Let $\tcF \subset A_\bQ$ be the fan of $\cU$,
i.e., for each open stratum $S\subset\oY$, we take the convex
hull $\tsigma_S$ of the rays associated with each of the boundary
divisors
that contains $S$.
Let $\cU_S \subset \cU$ be
the open subset associated to $\tsigma_{S}$
and $Z_S\subset \cU_S$ the closed orbit.

\begin{Lemma} \label{equivs} The following conditions are equivalent:
\begin{enumerate}
\item $T_{\Pic}$ acts freely on $\cU_S$.
\item $\tsigma_S$ maps isomorphically onto a
strictly simplicial cone in $N_Y$.
\item Boundary divisors that don't contain $S$ generate
$\Pic(\oY)$.
\item For each boundary divisor $D\supset S$, we can find
a unit $u \in \cO^*(Y)$ with valuation one on $D$, and
valuation zero on other boundary divisors containing~$S$.
\end{enumerate}
\end{Lemma}
\begin{proof} Immediate from the definitions. \end{proof}

\begin{Proposition}\label{nonseptoric}
If the (equivalent) conditions of Lemma~\ref{equivs} hold for all strata~$S$, then $T_{\Pic}$ acts freely on $\cU$,
with quotient a smooth (not necessarily separated) $T_Y$-toric variety,
the union of affine toric varieties $X(\sigma_S)$ where $\sigma_S \subset N_Y$
is the (isomorphic) image of $\tsigma_S \subset A_{\bQ}$. Denote this
toric variety by $X(\cF)$.
\end{Proposition}

\begin{proof}
$X(\cF)$ can be glued from affine charts $X(\sigma_S)=\cU_S/T_{\Pic}$ as
in the standard definition of a toric variety \cite{Oda}.
\end{proof}

Consider the sheaf of algebras $B := \bigoplus\limits_{D \in \Pic(\oY)} \cO_{\oY}(D)$
with multiplication given by
tensor product (more precisely, to define the multiplication we have to fix line bundles whose classes form a basis of $\Pic(\oY)$).
Let $W: =\Spec(B) \arrow^{q} \oY$.

\begin{Theorem} \label{coxmap}
$T_{\Pic}$ acts freely on $W$ with quotient $q:W \to \oY$.
There is a natural $T_{\Pic}$ equivariant map $f:\,W \to \cU$.
Suppose that $T_{\Pic}$ acts freely on $\cU$ with quotient $X(\cF)$ as in Prop.~\ref{nonseptoric}.
Then $f$ induces a natural map of quotients $\oY \to X(\cF)$.
The scheme-theoretic inverse image
of a toric stratum is a stratum of~$\oY$, and this establishes
a one-to-one correspondence between toric strata of $X(\cF)$
and boundary strata of $\oY$. The map
$\oY \times T_Y \to X(\cF)$ is smooth and surjective.
\end{Theorem}

\begin{proof}
We can check the first statement locally on $\oY$ and so may assume that all
line bundles are trivial. Then $W = \oY \times T_{\Pic}$ and the
statement is obvious.

The sections $1_D \in H^0(\oY,\cO(D))$ induce an
equivariant map
$f:W \to \bA^{\partial \oY}$ and a scheme-theoretic equality
$f^{-1}(Z_S) = q^{-1}(S)$. In particular, $f(W)\subset\cU$.
Next we show that $W \times \bG_m^{\partial\oY} \to \cU$ is smooth and surjective.
Since both spaces are smooth, it's enough
to check that all fibres are smooth (and non-empty) of the same
dimension.
If we pullback the map to the torus orbit $Z_S$ we obtain
$$
F:\,q^{-1}(S) \times \bG_m^{\partial\oY} \to Z_S.
$$
Since $\bG_m^{\partial\oY} \to Z_S$ is a surjective homomorphism,
$dF$ is everywhere surjective.
It follows that $F$ is a smooth surjection, of
relative dimension independent of $S$.

The statements for $\oY$ follow by taking the quotient by the
free action of $T_{\Pic}$.
\end{proof}

\begin{Definition}\label{M^Snotation}
Let $M^S_Y \subset M_Y$ be the sublattice generated
by units with zero valuation on all boundary
divisors containing $S$.
Note we have a canonical restriction map
$M^S_Y \to M_S = \cO^*(S)/k^*$ and thus a canonical map
\begin{equation}\label{anothermap}
S \to T^S_Y:=\Hom(M_Y^S,\bG_m).\cooltag
\end{equation}
\end{Definition}

\begin{Theorem} \label{hubsch}
Suppose that the conditions of Lemma~\ref{equivs} hold for all strata~$S$.
The map $\oY \to X(\cF)$ of Th.~\ref{coxmap}
is an immersion iff \eqref{anothermap} is an immersion for any~$S$.

If $S$ is very affine  and $M^S_Y  \to  M_S$ is surjective
then \eqref{anothermap} is an immersion.

If \eqref{anothermap} is an immersion for any $S$
and the image cones $\sigma_S \subset N_Y$ form a fan~$\cF$,
then $X(\cF)$ is the associated toric variety, $Y$ is
Sch\"on, and $\cF$ is tropical.
\end{Theorem}

\begin{proof}
$\oY \to X(\cF)$ is an immersion
iff each fiber is a single (reduced) point. We can
check this stratum by stratum on $X(\cF)$. But one checks immediately
from the definitions that the maps in the statement are exactly
the maps over the various torus orbits of $X(\cF)$. If $M^S_Y$
generates the units and $S$ is very affine, it is clear that the
map for this stratum is a closed embedding.

Finally, if the cones $\sigma_S$ form a fan, it is clear
that $X(\cF)$ is the associated
toric variety. Since the structure map is smooth by
Th.~\ref{coxmap}, the fan is tropical.
\end{proof}

\section{Iitaka fibration for subvarieties of algebraic tori} \label{ReidQ}

If $X \subset A$ is a smooth subvariety of an abelian variety then either $K_X$ is ample or $X$ is preserved by a positive dimensional subgroup 
$A' \subset A$. Here we prove an analogous statement for subvarieties of algebraic tori.
This answers a question posed to the second author by Miles Reid.

\begin{Theorem} 
Let $Y \subset T$ be a Sch\"on closed subvariety of an algebraic torus.
Exactly one of the following holds:
\begin{enumerate}
\item $Y$ is log minimal, or 
\item $Y$ is preserved by a nontrivial subtorus $S\subset T$.
\end{enumerate}
\end{Theorem}

\begin{proof}
Let $X$ be a smooth $T$-toric variety such that the closure $\oY$ of $Y$ in $X$
is smooth, projective, and transverse to the toric boundary. So $\oY$ is smooth
and the restriction $B$ of the toric boundary is a reduced simple normal crossing divisor.

By \cite[1.4]{Tevelev}, the line bundle $K_{\oY}+B$ is globally generated.
So $Y$ is log minimal iff ~$K_{\oY}+B$ is positive on any curve $C\subset\oY$
which is not contained in the boundary.

Suppose that $Y$ is not log minimal, and let $C$ be a curve as above such that 
$(K_{\oY}+B).C=0$. Then $C$ intersects $B$ 
(because $T$ is affine and so cannot contain a proper curve).
In particular, $l:=-K_{\oY}\cdot C=B\cdot C>0$.
Using translation by an element of $T$, we can assume that $C$ 
contains the unit element $e\in T$ and that $C$ is smooth at $e$.

By the bend-and-break argument, we can assume that $C$ is a rational curve 
and $l\le \dim Y+1$ (use \cite[II.5.14]{Ratcurves}, with $M:=K_{\oY}+B$).
Let $\nu:\,\bP^1\to C$ be the normalization. 

\begin{Claim} $l>1$.
\end{Claim}

\begin{proof}
Suppose $l=B\cdot C=1$.  Then $\nu^{-1}(B)$ is a singleton
and therefore its complement is $\bA^1$. But there are no non-constant maps $\bA^1\to T$.
\end{proof}

We choose $C$ as above of minimal degree with respect to some polarisation of $\oY$.
Let $P=\nu^{-1}(e)$.
By the bend-and-break argument \cite[1.11]{Wisniewski},
$\nu$ can be included in an irreducible 
$l$-dimensional family $Z$ of morphisms $\bP^1\to\oY$ such that $P \mapsto e$, and 
the locus $\oS\subset \oY$ of points spanned by these rational curves has dimension~$l-1$
(Wi\'sniewski assumes that the class of $C$ generates an extremal ray,
but we only need that the curve $C$ has minimal degree, see, e.g.,~\cite[IV.2.6]{Ratcurves}).

\begin{Claim}
$S:=\oS\cap Y$ is an algebraic subtorus of $T$. 
\end{Claim}

\begin{proof}
Let $C'$ be a rational curve from the family $Z$. 
Let $\nu':\,\bP^1\to C'$ be the normalization. Let $C'_0:=\bP^1\setminus\nu'^{-1}(B)$.
Then $B\cdot C'=-K_{X}\cdot C'=l$, and therefore $\nu'^{-1}(B)$ 
(set-theoretically) is a union of at most $l$ points.
It follows that the intrinsic torus $T'=\Hom(\cO^*(C'_0)/k^*,\bG_m)$
has dimension at most $l-1$.
By the universal property of the intrinsic torus, 
the morphism $C'_0\mathop{\longrightarrow}\limits^{\nu'} C'\cap Y\hookrightarrow T$ factors through~$T'$.
It follows that $C'\cap Y$
is contained in a subtorus of dimension at most $l-1$.
Since this is true for any curve in the family $Z$, and since subtori of an algebraic torus
do not deform, it follows that the locus $S$ spanned by them is contained in a subtorus of dimension at most $l-1$.
But $S$ itself is $(l-1)$-dimensional. Therefore $S$ is an algebraic subtorus.
\end{proof}

We claim that $Y$ is preserved by $S$. It suffices to prove that $Y$ is preserved 
by any fixed one-parameter subgroup $R\subset S$. 
The factorization morphism $T\to T/R$ extends to the equivariant $\bP^1$-bundle $\pi:\,A\to T/R$,
where the fan $\cF$ of $A$ consists of the ray spanned by $R$ and the opposite ray.
Let $\oR\simeq\bP^1\subset A$ be the closure of~$R$.
Since $R\subset Y$, $\cF$ belongs to the tropicalization $\cT(Y)$ of $Y$.
It follows from \cite[2.5]{Tevelev} that the closure $\tilde Y$ of $Y$ in $A$ is smooth
with normal crossing boundary $\tilde B$. 
Note that $\tilde Y$ is not complete and that $\oR\subset\tilde Y$.
We have a morphism $f \colon W \rightarrow \oS$ where $W$ is the total space of a $\bP^1$-bundle $p \colon W \rightarrow V$ 
(the universal family of smooth rational curves),
such that $f$ contracts a section of $p$ to $e \in \oS$, $f$ is finite over $\oS \setminus \{e\}$, and $f^*(K_{\oY}+B)$ is zero on fibers of $p$.
It follows that $K_{\oY}+B$ is numerically trivial on $\oS$. 
By \cite[1.4]{Tevelev}, the morphism $\pi \colon \tilde Y\to\oY$ is log crepant, i.e.,
$K_{\tilde{Y}}+\tilde{B}= \pi^*(K_{\oY}+B)$. Hence $(K_{\tilde{Y}}+\tilde{B})|_{\oR} = 0$.
Since the log canonical line bundle of $\oR$ is trivial, it follows that the determinant of the normal bundle $\cN_{\oR/\tilde Y}$ is trivial by adjunction.
The bundle $\cN_{\oR/\tilde Y}$ is a subbundle of $\cN_{\oR/A}$, which is trivial 
because it is the normal bundle of a fiber of a fibration. So 
$\cN_{\oR/\tilde Y}\simeq\oplus_{i=1}^k\cO(a_i)$ with $a_i\le 0$ for each $i$ and $\sum a_i =0$, i.e., $\cN_{\oR/\tilde Y}$ is the trivial bundle. 
Therefore,  by standard deformation theory, $\tilde Y$ is covered by deformations of $\oR$.
It follows that $Y$ is preserved by $R$.
\end{proof}

\section{Genus $3$ Curves: Moduli Stack $\tcM_3^{(2)}$} \label{branchcover}

For this section we assume $\ch k  \neq 2$.
We construct a Deligne--Mumford stack $\tcM_3^{(2)}$
that, as we will later prove, has $\oY^7_{\lc}$ as its coarse moduli space.
Throughout we use caligraphic font to indicate stacks, and
ordinary font for coarse moduli spaces.

Let $\Gamma$ be the dual graph of a stable curve,
with each vertex $v$  labeled by $g(v)$, the
genus of the corresponding connected component of the normalization.
Define
$$
\ocM_{\Gamma} = \prod_{v \in \Gamma^0} \ocM_{g(v),E_v}
$$
where $E_v$ is the set of incident edges locally about $v$ (so
a loop counts twice).

\begin{Theorem}[{\cite[Appendix]{GP}}]\label{strataofmg}
The closed strata of
$\ocM_g$ are in bijection with dual graphs of stable curves.
$\Aut(\Gamma)$ acts naturally on $\ocM_{\Gamma}$ and
$\ocM_{\Gamma}/\Aut(\Gamma)$ is the normalization of the stratum
corresponding to~$\Gamma$.
\end{Theorem}

\begin{Lemma} \label{smoothstratum} The closed stratum corresponding
to the dual graph $\Gamma$ fails to be smooth iff there is
a dual graph $A$ and two distinct sets of edges $N_1, N_2$ of $A$
such that the curve corresponding to $\Gamma$ can be
obtained from the curve corresponding to $A$
by smoothing nodes corresponding to $N_1$ (resp.~corresponding to
$N_2$).
\end{Lemma}

\begin{proof} We consider the versal deformation space of
the singularities. It is smooth, with one variable for each
node, and coordinate hyperplanes, each corresponding to
smoothing all but a single node. To each coordinate
subspace (intersection of coordinate hyperplanes) is
attached a dual graph, and two subspaces are local analytic
branches of the same global stratum iff they have the same
dual graph. Inclusion of strata corresponds to smoothing
of nodes. Now the result is clear.
\end{proof}

\begin{Review}
Strata of $\ocM_3$ are listed in~\cite[pg. 340--347]{Fabersthesis}.
We draw some of them in Fig.~\ref{m3strata}.

\begin{figure}[htbp]
  \includegraphics[width=4in]{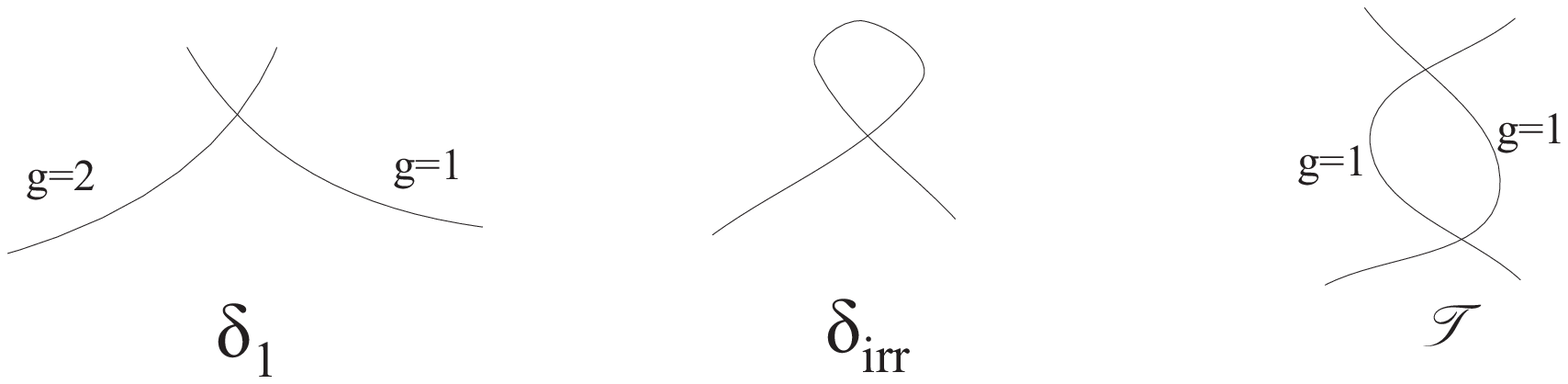}\\
  \caption{\small Some Strata of $\ocM_3$}\label{m3strata}
\end{figure}

Strata $\delta_1$ and~$\delta_{irr}$ have codimension $1$ and
$\cT$ has codimension $2$.
These strata are closed.
Checking the criterion from Lemma~\ref{smoothstratum} for each substratum of $\cT$ gives
\end{Review}

\begin{Proposition} \label{Tsmooth} The substack $\cT$ is smooth and
isomorphic to
$$
(\ocM_{1,2} \times \ocM_{1,2})/(\bZ/2\bZ)^{2},
$$
where one generator of $(\bZ/2\bZ)^2$ interchanges the two
copies of $\ocM_{1,2}$ and the other acts diagonally on
the two copies, interchanging the two marked points.
\end{Proposition}

\begin{Definition} \label{tcM3}
Let $\cH\subset\cM_3$ be the hyperelliptic divisor
with closure $\ocH\subset\ocM_3$. By
the boundary of $\ocM_3$ we mean the complement to $\cM_3 \setminus
\cH$
(the Deligne--Mumford boundary $\delta_1$, $\delta_{irr}$ plus $\ocH$).
Let $\tcM_3$ be the blowup of $\ocM_3$ along $\cT$.
For a divisor $\cD$ on $\ocM_3$ we write $\tcD$ for
its strict transform on~$\tcM_3$.
Let $\cE$ be the exceptional divisor.
\end{Definition}

\begin{Proposition}\label{tm3nc}\label{delta1t}
$(\tcM_3,\cB)$ is smooth with
normal crossings. $\tcH \to \ocH$ is an isomorphism.
The normalization of $\tilde{\Delta}_1$ is
$\oM_{1,1} \times \oM_{0,7}/S_6$.
\end{Proposition}

\begin{proof}
$\ocH \subset \ocM_3$ is a smooth substack with coarse moduli space isomorphic
to $\oM_{0,8}/S_{8}$.
Let $\cB_i$ ($i=2,3,4$) be the boundary divisor of
$\ocH$ that corresponds to stable rational curves with two components
and with $i$ marked points on one of the components.
By \cite{CornalbaHarris},
$\delta_1$~meets $\ocH$ transversally in~$\cB_3$;
$\cT$ is contained in $\ocH$ and equal to
the substack $\cB_4$; $\delta_{\irr}$ meets~$\ocH$
transversally along $\cB_2$, but $\delta_{\irr}$
has two transverse analytic branches along $\cT$,
each of which meets $\ocH$ transversally. In terms of Weil
divisors
$\delta_{\irr}|_{\ocH} = \cB_2 + 2 \cB_4$.
It~follows that $(\tcM_3,\cB)$ is smooth with
normal crossings and $\tcH \simeq \ocH$.

The normalization of $\delta_1$ is
$\ocM_{1,1} \times \ocM_{2,1}$ by Th.~\ref{strataofmg}.
Since $\cT$ is transverse to $\delta_1$, the normalization
of $\tdelta_1$ will be the blowup of the inverse
image of $\cT$ on the normalization of $\delta_1$. This
inverse image is the product of $\ocM_{1,1}$ with
the codimension two stratum~$S$ of $\ocM_{2,1}$ corresponding
to a pointed rational curve meeting a genus~$1$ curve in
two points. By \cite[\S 10]{Rullasthesis} the coarse moduli space of the blowup of
$\ocM_{2,1}$ along $S$ is $\oM_{0,7}/S_6$. \end{proof}

\begin{Definition}
Consider the locally constant sheaf $R_{et}^1c_*(\bZ/2\bZ)$, where
$c:\cC \to \cM_3$ is the universal curve.
The associated locally constant sheaf of symplectic bases
gives a Galois \'etale cover $\cM_3^{(2)} \to \cM_3$ with group $G = \Sp(2,6)$
\cite[\S5]{DM}.
For any $\cX$ birational to $\cM_3$ we write
$\cX^{(2)}$ for the normal closure of $\cX$ in the function field
of $\cM_3^{(2)}$, and for any closed substack $\cZ \subset \cX$ we write $\cZ^{(2)}$
for its inverse image in $\cX^{(2)}$.
\end{Definition}

\begin{Theorem}\label{almoststacksmooth}\label{incidencies}\label{stacksmooth}
The pair $(\tcM_3^{(2)},\cB)$ is smooth with simple normal crossing boundary, in particular,
irreducible components of boundary divisors are smooth.
They are in $W(E_7)$-equivariant bijection with certain elements of
$\cR(E_7)$:
\begin{table}[htbp]
\begin{tabular}{|c|c|c|c|c|c|}
  \hline
  $\Theta$                      & $A_1$                  & $A_2$                                             & $A_3\times A_3$& $A_7$        \\
  \hline
  type of $\cD(\Theta)$           & $\tdelta_{\irr}^{(2)}$ & $\tdelta_1^{(2)}$                                 & $\tilde\cE^{(2)}$& $\tcH^{(2)}$ \\
  \hline
  $D(\Theta)$ is isomorphic to  &                 & \vphantom{$F^{F^F}$}$\oM_{0,4} \times \oM_{0,7}$  & $\oM_{0,5}\times\oM_{0,5} \times \oM_{0,4}$& $\oM_{0,8}$  \\
  \hline
\end{tabular}
\end{table}
\end{Theorem}

\begin{Remark}\label{rootsEn}\label{Sp(2,6)}
Let us recall how $\Sp(2,6)$ is related to $E_7$.
Let $\Pic_n$ be the Picard group of the blow-up
of $\bP^2$ in $n$ points with the intersection pairing
$h\cdot h=1$, $e_i\cdot e_j=-\delta_{i,j}$, $h\cdot e_i=0$.
Let $K_n=-3h+e_1+\ldots+e_n$ be the canonical class. Then
$E_n:=\{\alpha\in (K_n)^\perp\ |\ \alpha^2=-2\}$.
For example, positive roots of $E_7$ are given by
\begin{equation}
\alpha_{ij}=e_i-e_j;\quad
\alpha_{ijk}=h-e_i-e_j-e_k;\quad
\beta_{i}=2h-e_1-\ldots-\hat e_i-\ldots -e_7.
\cooltag
\end{equation}
Sometimes we skip $\alpha$ and $\beta$ and denote roots simply as $ij$, $ijk$, and $i$.

For any $v\in\Pic_7$, let $\bar v=v\otimes1\in\Pic_7\otimes\bF_2$.
Then $(\bar K_7)^2=0$ and the induced bilinear form on
$V:=(\bar K_7)^\perp/(\bar K_7)\simeq\bF_2^6$
is nondegenerate and symplectic. We have
$$V\setminus\{0\}=\{\bar\alpha\,|\,\alpha\in E_7^+\}.$$
The Weyl group $W(E_7)$
is generated by reflections $x\mapsto x+(\alpha,x)\alpha$ for
$\alpha\in E_7$.
The induced transformation of $V$ is a symplectic transvection
$x\mapsto x+(\bar\alpha,x)\bar\alpha$.
The corresponding homomorphism
$W(E_7)\to G$ is surjective  with kernel $\{\pm E\}$.
\end{Remark}

\begin{Remark}\label{e7calcs}\label{DynkinTh}
We frequently use well-known facts about conjugacy
classes of root subsystems in an irreducible root system $\Delta$.
Most of them can be summarized in the
following observation of Eugene Dynkin:
any irreducible root subsystem of $\Delta$ is $W(\Delta)$-equivalent
to a subsystem such that its Dynkin diagram is contained in the affine Dynkin diagram of $\Delta$.
For example, the affine Dynkin diagram of $E_7$ is
$$\Esea\W\W\W\W\W\W\W\W$$
It is easy to see that
all subsystems of $E_7$ of types
$A_1$, $A_2$, $A_3\times A_3$, $A_7$ (vertices of $\cR(E_7)$) are conjugate.
Moreover, the action of $W(E_7)$ on $\cA_7(E_7)$ is $2$-transitive and the intersection of any $2$ $A_7$'s
is $A_3\times A_3$. Notice also that $A_2^\perp$ is isomorphic to $A_5$ but not all $A_5$'s are conjugate:
there is another $W(E_7)$-orbit.
\end{Remark}

{\em Proof of Th.~\ref{stacksmooth}.}\/
Let  $(\ocM, \cB)$ be the formal neighborhood
of a point on $\tcM_3$ or $\ocM_3$, let $\cM := \ocM \setminus \cB$.
By \cite{GM}, the tame algebraic fundamental group of $\cM$ is $\hat\bZ^m$, with one
generator~$e_i$ for each component of $\cB$.
The cover $\pi:\,\ocM^{(2)}\to\ocM$ gives rise
to a homomorphism $h: \bZ^m \to G$ that factors
through a product of cyclic groups $\tilde h:\,\Gamma_1\times\ldots\Gamma_m\to G$
(here we take each $|\Gamma_i|$ minimal). We use a variation of Abhyankar's lemma:

\begin{Lemma}[\cite{GM}]\label{bc}
The cover $\pi:\,\ocM^{(2)}\to\ocM$ is a disconnected union of $|G|/|\Gamma|$
generalized Kummer coverings, where $\Gamma$ is the image of $h$.
More precisely, $\pi$ is isomorphic
to the restriction of\/ $\coprod(\bA',B') \to (\bA,B)$ to the germ of $(\bA,B)$
at the origin, where $\bA = \bA^m$ (resp.~$\bA'$) is the toric variety
associated to the first octant in the lattice $\bZ^m$
(resp. the intersection of the octant with $\Ker h$).
$(\ocM^{(2)},\cB^{(2)})$ is smooth with normal crossings if and only if
$\tilde h$ is an isomorphism. In this case, if $\cD \subset \ocM$ is a smooth divisor
transverse to the boundary, then ${\pi}^{-1}(\cD)$ is smooth.
\end{Lemma}

The monodromy around $\delta_1$ and $\delta_{irr}$ is well-known, see e.g. \cite{GeemenOort}.
The analytic picture is as follows:
a small disc transverse to a boundary divisor corresponds
to a degeneration in which a smooth curve $C$ degenerates
to a singular curve $C'$. The special fibre is
a deformation retract of the total family (over the
disc) and this induces a surjective map on homology
$H_1(C,\bZ) \to  H_1(C',\bZ)$
the kernel of which is generated by one class, the
so called vanishing cycle.
We use the same term for
the generator of $H^1(C',\bZ)^{\perp}\subset H^1(C,\bZ)$ and the
analogous classes for other coefficient rings, and
refer to \cite[2.10]{GeemenOort} for an algebraic description of vanishing cycles.
The corresponding monodromy action on $V=H^1_{et}(C,\bZ/2\bZ)$ is by
symplectic transvection by the vanishing cycle $\alpha$.

The monodromy action is trivial for loops around $\delta_1$,
as the vanishing cycle is trivial.
The monodromy around $\cE$ is a composition of monodromies
around two branches of $\delta_{\irr}$ meeting at $\cT$.
A generic point of $\cT$ is a union of two elliptic curves
glued at two points and the vanishing cycles for these two nodes
are non-trivial but equal modulo~$2$.
It follows that the monodromy around $\cE$ is also trivial.
Thus $\ocM_3^{(2)}\to\ocM_3$
(resp.~$\tcM_3^{(2)}\to\tcM_3$) is \'etale outside of $\delta_{\irr}$
(resp.~$\tdelta_{\irr}$).

Consider the formal neighborhood in $\tcMM$ as in Lemma~\ref{bc}.
We claim that $\tilde h$ is an isomorphism and therefore
$(\tcMM,\cB)$ is smooth with normal crossings.
By Prop.~\ref{Tsmooth}, $\cT$ has no self-intersection.
It follows that strata contained in~$\cT$ are characterized by the
existence of a unique pair of nodes which, when simultaneously
removed, disconnect the curve into two connected curves.
The equivalent condition on the vanishing cycles
is that they are non-trivial, but equal modulo~$2$. Call this class
$\alpha$ if it exists. Note that such a stratum
corresponds to two disjoint strata in $\tcM_3$ (that both belong to~$\cE$),
with the same local monodromy generated by transvections by $\alpha$, and the
other non-trivial vanishing cycles (for the nodes of type~$\delta_{\irr}$ other than the distinguished pair).
We get a collection of different and pairwise orthogonal
vanishing cycles. Since their number is less than~$7$ (because $(\tcM_3,\cB)$ has normal crossings),
the following lemma implies that $\tilde h$ is an isomorphism.

\begin{Lemma}\label{stupido}\label{monodromy}
Let $S \subset V\setminus\{0\}$ be a set of pairwise
orthogonal vectors. Let $H \subset G$ be the subgroup
generated by transvections by elements of $S$. Then
$H = (\bZ/2\bZ)^{|S|}$ unless $S\cup\{0\}$ is the maximal isotropic
subspace
of $2^3=8$ vectors.

More generally, a subset $S\subset V\setminus\{0\}$ is a reduction of a root
subsystem of $E_7$
modulo~$2$ if and only if it is closed under transvections by its
elements.
In this case the group generated by these transvections
is isomorphic to the Weyl group of the root subsystem or to its
quotient
modulo $\{\pm E\}$ if the latter is contained in it.
\end{Lemma}

\begin{proof}
The second statement follows from the simple fact that  a subset of a root system
is a root subsystem if and only if it is closed under reflections
by its elements.

For the first statement, $S$ is obviously closed under transvections.
The corresponding root subsystem is $A_1^{|S|}$ because any other
subsystem
contains non-perpendicular roots.
The corresponding Weyl group is $(\bZ/2\bZ)^{|S|}$.
It is obvious that a group generated by reflections in perpendicular
vectors
in $\bR^7$ contains $-E$ if and only if these vectors form a basis,
i.e., there are $7$ of them.
\end{proof}

It remains to prove that the boundary divisors are smooth (so that the boundary has \emph{simple} normal crossings), 
to identify them, and to enumerate them by elements of $\cR(E_7)$.

\begin{Lemma} \label{delta1}
Irreducible components
of $\tdelta_1^{(2)}$
are smooth and in one-to-one $G$-equivariant
correspondence with symplectic $2$-planes $W \subset V$.
The latter  $G$-set is canonically identified with $\cA_2(E_7)$.
The Galois group of the branched cover
$\pi: \tilde{\Delta}_{1,W} \to \Norm(\tilde{\Delta}_1)$
 is isomorphic to $S_3 \times S_6$
(we write $\Norm(X)$ for the normalization of the reduced space).
The branched cover~$\pi$ is the natural quotient map
$$
\oM_{0,4} \times \oM_{0,7} \to \oM_{0,4}/S_3 \times \oM_{0,7}/S_6
$$
\end{Lemma}

\begin{proof}
We begin by identifying the Galois group --- which by the
general theory is the subgroup of $\Sp(2,6)$ that preserves
a given irreducible component of $\delta_1^{(2)}$. We note
that the local system $R^1\pi_*(\bZ/2\bZ)$ is canonically
split over the open stratum $\delta_1^0$, corresponding
to the splitting of the cohomology of a stable curve
with a single node of type $\delta_1$, $V = W \oplus W^{\perp}$
where $W$ is the cohomology of the genus~$1$ component. It
follows that the stabilizer is contained in the subgroup
preserving~$W$, which is equal to $\Sp(2,W)\times\Sp(2,W^\perp)$.
Notice that $W$ is nondegenerate and
the set of subspaces of given signature forms a single $\Sp(2,6)$-orbit by the Witt theorem.
Any subspace is obviously closed under transvections
by its elements and therefore corresponds to a root subsystem by Lemma~\ref{stupido}.
It easy to see that this subsystem is  $A_2$.
Therefore the set of $2$-dimensional
symplectic subspaces $W\subset \bF_2^6$ is $\Sp(2,6)$-equivariantly
identified with $\cA_2(E_7)$ and
$\Sp(W)$ is equal to $W(A_2)=S_3$ (note that the symplectic group is generated by symplectic transvections).
By a similar argument $W^\perp$ is identified with $A_2^\perp=A_5$
and $\Sp(W^\perp)=W(A_5)=S_6$.

Therefore  $\Sp(2,W)\times\Sp(2,W^\perp)$ is isomorphic to
$S_3\times S_6$ and is generated
by transvections by classes in either $W$ or~$W^{\perp}$.
Such transvections are given by monodromy around boundary
divisors of $\delta_1$, corresponding to further degenerations
of the curve where we shrink to a point a loop on either the
elliptic, or genus $2$ component. Thus it follows that
the Galois group is precisely $S_3 \times S_6$. It is clear from this
discussion that the cover of $\delta_1$
is (at least generically) the product of the covers $\ocM_{1,1}^{(2)} \to \ocM_{1,1}$ and $\ocM_{2,1}^{(2)} \to \ocM_{2,1}$.
By \cite[VIII.3]{DO} and Prop.~\ref{delta1t}, after pulling back to $\tcM_3$ and passing to the coarse moduli spaces,
these covers are the quotients $\oM_{0,4} \to \oM_{0,4}/S_3$ and $\oM_{0,7} \to \oM_{0,7}/S_6$. So it is enough
to prove that components of $\tdelta_{\irr}^{(2)}$ are normal.

Since they are Cohen-Macaulay (because the ambient stack is smooth)
it's enough to check there is no self intersection in codimension
one. Such a self intersection would lie over a codimension
one self intersection of $\delta_1$. There is a unique such
stratum in~$\ocM_3$, $\Delta_{1,1}$ of \cite[pg. 340]{Fabersthesis}.
This corresponds to a curve $C$ with two nodes, each of type
$\delta_1$. It's clear that smoothing the two nodes gives rise
to distinct and orthogonal $W_1,W_2 \subset H^1(C,\bZ/2\bZ)$,
and thus, on the cover, local analytic branches that belong
to distinct irreducible components. \end{proof}

\begin{Lemma} \label{H} Irreducible
components of $\oH^{(2)}$ are smooth, isomorphic to
$\oM_{0,8}$ and in one-to-one $G$-equivariant correspondence
with
quadratic forms $Q$ of plus type inducing the symplectic form of
$\bF_2^6$.
This  $G$-set is canonically identified with $\cA_7(E_7)$.
Each component $\ocH_{Q}$ is isomorphic to its strict transform on~$\tcMM$. Irreducible components of $\tcH^{(2)}$
are pairwise disjoint. The Galois group of the cover
$\oH_{Q} \to \oH$
is the subgroup of $G$ preserving $Q$. It is isomorphic
to~$S_8$ naturally acting on $\oM_{0,8}$.
\end{Lemma}

\begin{proof}
A nonsingular quadratic form on $\bF_2^{2m}$ has \emph{plus type} if there exists a totally isotropic subspace
of dimension $m$ \cite[p.~xi]{Atlas}.
The symplectic group acts transitively on the quadratic forms of plus type inducing the symplectic form.
The parametrization of irreducible components of $\oH^{(2)}$ by quadratic forms of plus type and description
of the Galois group is given in \cite[VIII.3]{DO}, where in addition
the cover $\oH_{Q} \to \oH$ is generically identified with $M_{0,8} \to M_{0,8}/S_8$.

Irreducible components of $\tcH^{(2)}$ are smooth by  \eqref{bc}
since $\tcH$ is smooth and transverse to the branch locus of $\pi: \tcMM \to \tcM_3$.
In particular they are normal, so, since $\oH=\oM_{0,8}/S_8$, the cover
$\oH_Q \to \oH$ is the quotient $\oM_{0,8} \rightarrow \oM_{0,8}/S_8$.
To show that each component of $\tcH^{(2)}$ maps isomorphically onto its image
in $\ocM_3^{(2)}$ it suffices to show that each component of $\ocH^{(2)}$ is normal.
It is enough to check singular points of the cover, in codimension
one over $\ocH$, and there is only one such stratum,~$\cT$, of
$\ocM_3$. A generic point of $\cT$ is a union of two elliptic curves
glued at two points and the vanishing cycles for these two nodes
are non-trivial but equal modulo~$2$.
It follows by Th.~\ref{bc} that the cover
$\ocM_3^{(2)}\to \ocM_3$ is given
locally by
(the product with $\bA^4$ of)
$(z^2 = xy) \subset \bA^3 \to \bA^2$ by dropping the $z$
coordinate. Here $x=0$, $y=0$ are the
two analytic branches of $\delta_{\irr}$, each of which meets
$\ocH$ transversally along $\cT$. Thus the inverse image
of $\ocH$ has two smooth analytic branches meeting
transversally. It is enough to argue that the two branches belong
to different irreducible components of $\ocH^{(2)}$.
Fix one and let
$Q$ be the corresponding quadratic form. Let $\alpha$
be the vanishing cycle in $V$ (for either of the two nodes).
We compute $Q(\alpha)$.

Consider a smooth hyperelliptic curve degenerating to a curve
of type $\cT$. We obtain this by taking two genus $1$ curves, and
on each choosing a disc and its image under a hyperelliptic
involution, and now joining the two curves, to obtain a hyperelliptic
curve with two thin collars, interchanged by the involution.
In the degeneration the collars contract to a pair of conjugate
points.  $\bF_2^6=H^1(C,\bZ/2\bZ)$ is identified with even cardinality
subsets
of the branch locus, with addition exclusive or, modulo identifying
a subset with its complement. $Q$ assigns to
a set half its cardinality, modulo $2$, see \cite[VIII.3]{DO}. The
symplectic form assigns to a pair of subsets the cardinality
of their intersection, modulo $2$. Under
this identification, the collar circle $\alpha$ is represented
by the sum of the $4$ branch points on one of the two elliptic
curves. Thus $Q(\alpha) =0$.
Transvection by $\alpha$ interchanges the two sheets of the
local analytic model of the cover of $\ocM_3$ above. It does not preserve
$Q$: Let $A$ be the subset representing $\alpha$ and $Z$
a pair of points with $|A \cap Z| =1$. One checks
that $Q(Z) =1$ but $Q(Z + (Z,A)A) = 0$.
Thus the other analytic branch of the inverse image of
$\ocH$ belongs to a different irreducible component.

It remains to construct an $A_7\subset E_7$ using the form $Q$.
Recall the description of the root system $E_7$ from \eqref{rootsEn}.
As above, we identify $\bF_2^6$ with
even cardinality subsets of the set $\{1,\ldots,8\}$
modulo identifying a subset with its complement.
Then $\bar\alpha_{ij}$ corresponds to $\{i,j\}$, $\bar\alpha_{ijk}$
corresponds to $\{1,\ldots,7\}\setminus\{i,j,k\}$,
and $\bar\beta_{i8}$ corresponds to $\{1,\ldots,7\}\setminus\{i\}$.
The set $\{Q(x)=1\}$ then corresponds to a subsystem
\begin{equation}\label{standA7}
\{\alpha_{ij},\ \beta_{i8}\}.\cooltag
\end{equation}
It is easy to see that this is a root subsystem of type $A_7$.
\end{proof}

There is a unique $W(E_7)$-orbit of root systems of type $A_1^7$ in $E_7$.
We call such a root system a \emph{Fano simplex} (because it corresponds to the set of
$(-2)$-curves on the surface over a field $k$ of characteristic $2$ obtained by blowing up the seven points of the Fano plane 
$\bP^2(\bF_2)$ in $\bP^2_k$).

\begin{Lemma} \label{deltairr}  Irreducible components
of $\tilde\delta_{\irr}^{(2)} \subset \tcMM$ are smooth.
The $G$-set of irreducible components is canonically identified with
$\cA_1(E_7)\simeq V\setminus\{0\}$.
A collection of $\tdelta_{\irr,\alpha}$ divisors
has non-trivial intersection iff the roots $\alpha\in E_7$ are
pairwise orthogonal and don't form a Fano simplex.
\end{Lemma}

\begin{proof} Let $W \subset \ocM_3$ be the complement
of all the codimension two boundary strata.
To compute the set $I$ of irreducible components of $\delta_{\irr}^{(2)}$
we may replace $\ocM_3$ by the formal completion $U$ of $W$ along
$D:=\delta_{\irr} \cap W$.
Then $I=G/H$ for $H \subset G$ the image of the monodromy map
$\pi^D_1(U) \to G$.
We have an exact sequence \cite[7.3]{GM}
$$\{e\}\to K\to \pi^D_1(U)\to \pi_1(D)\to\{e\},$$
with $K$ central and cyclic, generated by monodromy
around $\delta_{\irr}$ (this is the algebraic analog
of the fibration sequence for the unit sphere bundle in
the normal bundle to a submanifold). The image of $K$ in $H$
is generated by transvection by a vanishing cycle $\alpha$.
As $K$ is central,
$H \subset G_\alpha$. It is easy to check by direct calculation that
$G_\alpha$ is generated by transvections by elements in $\alpha^{\perp}$,
and as in the proof of
Lemma~\ref{delta1} we can realize such a transvection by monodromy
in $U$ (by considering degenerations contracting two disjoint
cycles to points). Thus $H = G_\alpha$. So as a $G$-set
$$I = G/G_{\alpha} = V\setminus\{0\}=\cA_1(E_7).$$

Now we argue that each component $\tdelta_{\alpha}$  of
$\tdelta_{\irr}^{(2)}$
is smooth. As in the proof of \eqref{delta1} we consider
a divisor of self intersection. There are two such
divisors for $\delta_{\irr} \subset \ocM_3$, one of which,
$\cT$, is removed by the blowup. So a self intersection
divisor of $\tdelta_{\alpha}$ must lie over the strict
transform of the stratum $\Delta_{0,0}$ of \cite[pg.
340]{Fabersthesis}.
This corresponds to a curve with two $\delta_{\irr}$ nodes, for
which there are two distinct vanishing cycles. The cover is smooth
over this stratum, and $\tdelta^{(2)}$ has two analytic braches.
They belong to the two
irreducible components corresponding to the two vanishing cycles.

Finally an intersection point of $k$ $\delta_{\alpha_i}$ divisors
lies over a $k$-fold self intersection of $\delta_{\irr}$. Such
strata are listed in \cite[pg 340--347]{Fabersthesis}. Clearly
the collection of $k$ vanishing cycles are pairwise orthogonal.
In the other direction, it's easy to check
that any two $k$ element subsets of pairwise orthogonal elements
of $V$ are conjugate under $\Sp(2,V)$ unless $k=3$ or $4$ in which case
there are two conjugacy classes. It follows that all
such intersections occur.
\end{proof}

\begin{Lemma} \label{E}
Irreducible
components of $T^{(2)}$ (resp $E^{(2)}$) are smooth, pairwise disjoint, and isomorphic to
$\oM_{0,5}\times\oM_{0,5}$  (resp. $\oM_{0,5}\times\oM_{0,5}\times
\oM_{0,4}$).
They are in one-to-one correspondence with root subsystems
of type $A_3 \times A_3$.
\end{Lemma}

\begin{proof}
The description of $T^{(2)}$ follows
from Lemma ~\ref{H} since $T \subset \oH=\oM_{0,8}/S_8$ is the
boundary divisor $B_4$. The correspondence with root
systems of type $A_3 \times A_3$ now follows from \eqref{e7calcs},
where it is noted that $A_3 \times A_3 \subset E_7$ are in one to
one correspondence with pairs $A_7 \neq A_7' \subset E_7$, under
$A_7 \cap A_7' = A_3 \times A_3$. The local analytic description of the cover $\ocM_3^{(2)} \rightarrow \ocM_3$
in the proof of Lemma~\ref{H} shows the components of $T^{(2)}$ are pairwise disjoint.

In order to describe $E^{(2)}$, we need to analyze the stacky structure at $\cT^{(2)} \subset \ocM^{(2)}_3$.
Recall that the substack $\cT \subset \ocM_3$ is the intersection of two
branches of $\delta_{\irr}$, and that $\ocH \subset \ocM_3$ also contains
$\cT$ and intersects each branch of $\delta_{\irr}$ transversely.
Locally over a point of $\cT$, a connected component of the cover
$\ocM^{(2)}_3 \rightarrow \ocM_3$ is the double cover branched over
$\delta_{\irr}$. Thus on $\ocM^{(2)}_3$, a slice transverse to
$\cT^{(2)}$ is an $A_1$-singularity, and each component of $\cT^{(2)}$ is
contained in two branches of $\delta^{(2)}_{\irr}$ and two branches of
$\ocH^{(2)}$, any two of which intersect tranversely (cutting out $\cT^{(2)}$).
The exceptional divisor $\cE^{(2)}$ of the blowup
$\tcM_3^{(2)} \rightarrow \ocM_3^{(2)}$ is a $\bP^1$-bundle over $\cT^{(2)}$,
and on each component the strict transforms of the branches of $\delta_{\irr}$
and $\ocH$ give 4 disjoint sections.

The irreducible components of the stack $\ocH^{(2)}$ have coarse
moduli space $\oM_{0,8}$. For a smooth hyperelliptic
curve $C$, the hyperelliptic involution acts trivially on $H^1(C,\bZ/2\bZ)$,
so preserves a given 2-level structure. It follows that $\ocH^{(2)}$ has
stabiliser $\bZ/2\bZ$ at a general point given by the hyperelliptic involution.
Since $\ocH^{(2)}$ and its coarse moduli space are both smooth, the
automorphism group of an arbitrary point is generated by the specialisations
of automorphisms in codimension one. We find that the automorphisms are generated by
the hyperelliptic involution and the involution along $\cB_3$ given by the hyperelliptic
involution of the elliptic component, which is the restriction of the involution along
$\delta_1$.

The components of $\cT^{(2)}$ are the closed substacks of $\ocH^{(2)}$ with
coarse moduli spaces the boundary divisors of type $\oM_{0,5} \times \oM_{0,5}$
in $\oM_{0,8}$. The components of $\cE^{(2)}$ are $\bP^1$-bundles
over these stacks. Note that the stabiliser $\bZ/2\bZ$ of a general point of
$\cT^{(2)}$ acts
nontrivially on the fibre of the bundle: it fixes the sections given by $\tcH^{(2)}$ pointwise
and interchanges the sections given by $\tilde{\delta}_{\irr}^{(2)}$
(because the hyperelliptic involution of a general point $[C]$ of $\cT$ interchanges the
two nodes of $C$). The remaining generators of the stabiliser of an arbitrary point of $\cT^{(2)}$
are given by the involution along $\delta_1^{(2)}$ so act trivially on the fibre.
Hence a component of the coarse
moduli space $E^{(2)}$ is a $\bP^1$-bundle over $\oM_{0,5} \times \oM_{0,5}$
with 3 disjoint sections given by one component of
$\tilde{\delta}^{(2)}_{\irr}$
and two components of $\tcH^{(2)}$. It follows that the bundle is trivial,
so each component of $E^{(2)}$ is isomorphic to $\oM_{0,5} \times \oM_{0,5} \times \oM_{0,4}$, as required.
\end{proof}

\section{Irreducible Representation $M(\Delta)$ of the Weyl Group}\label{weylg}

\begin{Definition}
Let $\Delta=\{\alpha\in\Lambda\,|\,\alpha^2=-2\}$ be an A-D-E root system,
where $\Lambda$ is a negative-definite $\bZ$-lattice spanned by  $\Delta$.
Consider the linear map
$$\phi:\,\Sym^2\Lambda^\vee\to\bZ^{\cA_1(\Delta)},
\quad f\mapsto \sum_{A_1\in\cA_1(\Delta)}f(A_1)[A_1],$$
where $\Sym_2(\Lambda^{\vee})$ is the
space of quadratic forms on $\Lambda$ and
$f(A_1)$ is equal to the value of $f$ on one of the two (opposite) roots of $A_1$.
We define $N(\Delta):=\Coker\phi$.
\end{Definition}

We choose positive roots $\Delta_+\subset\Delta$
and identify $\cA_1(\Delta)$ with $\Delta_+$.
Let $\{\alpha_i\}$ be the associated simple roots, which we identify with vertices of
the Dynkin diagram~$\Gamma$.
Recall that any root $\alpha\in\Delta_+$ is a nonnegative linear combination
$\sum n_i\alpha_i$ of simple roots and the support $\Supp\alpha=\{\alpha_i\,|\,n_i\ne0\}\subset\Gamma$
is a connected subgraph.

\begin{Proposition}\label{exsequn}
We have an exact $W$-equivariant sequence of free $\bZ$-modules
\begin{equation}\label{seqcoolformular}
0\to\Sym^2\Lambda^\vee\arrow^\phi\bZ^{\Delta_+}\arrow^\psi
N(\Delta)\to0.\cooltag
\end{equation}
Let $M(\Delta)$ be the dual lattice with the induced embedding
$\psi^\vee:\,M(\Delta)\hookrightarrow\bZ^{\Delta_+}$. We have \eqref{funnyMdef}.
Let $\Delta_+^T\subset\Delta_+$ be the set of roots with $3$-legged support.
Then $\rk N(\Delta)=|\Delta_+^T|$
and the restrictions $n_\alpha|_{\Im\psi^\vee}$ for $\alpha\in\Delta_+^T$
are coordinates on $M(\Delta)$.
\end{Proposition}

\begin{proof}
Since $\Gamma$ is a tree, there exists a (not unique) total ordering
$\prec$
on the set of pairs of vertices
$S=\{(i,j)\,|\,1\le i\le j\le\rk\Delta\}$ that satisfies the following
property.
If $(i,j)\prec(i',j')$
then either $i'$ or~$j'$ is not contained in the minimal substring
(a~tree with two legs)
of $\Gamma$ that contains both $i$ and~$j$.
We denote this substring by $[i,j]$.
It is well-known that $\alpha_{i,j}:=\sum\limits_{k\in[i,j]}\alpha_k$
is a root of $\Delta$, moreover, $\alpha_{i,j}$ is the unique root with support $[i,j]$.
Let $\{\omega_i\}$ be the basis of $\Lambda^\vee$ dual
to~$\{\alpha_i\}$.
It follows from the property of the ordering and the definition of
$\alpha_{i,j}$ that
$$\omega_i\omega_j(\alpha_{i',j'})=\begin{cases}
1& (i',j')=(i,j)\cr
0& (i',j')\prec(i,j)\cr
\end{cases}$$
It follows that $\phi$ is injective and that $N(\Delta)$ is
torsion-free.
Since the $\alpha_{i,j}$'s are the only roots in $\Delta_+\setminus\Delta_+^T$,
$\rk N(\Delta)=|\Delta_+^T|$.
Finally, \eqref{funnyMdef} follows by dualizing \eqref{seqcoolformular}.
\end{proof}

\begin{Corollary}\label{torsionfree}
$M(\Delta)$ is an irreducible $W$-module
of rank given in the table
\begin{table}[htbp]
\begin{tabular}{|c||c|c|c|c|c|c|c|}
  \hline
  $\Delta$        &$A_n$          &    $D_n$          & $E_6$ & $E_7$ & $E_8$ \\
  \hline
  $\rk M(\Delta)=|\Delta^T_+|$ &$0$             &    $n(n-3)\over2$  & $15$ & $35$ & $84$ \\
  \hline
\end{tabular}
\end{table}
\end{Corollary}

\begin{proof}
Irreducibility follows from the tables of characters in \cite{Atlas}.
\end{proof}

\begin{Lemma}\label{subs}
If $\Delta'\subset\Delta$ then we have a commutative diagram
\begin{equation}\label{basicdiagram}
\begin{CD}
0 @>>> \Sym^2\Lambda^\vee(\Delta) @>>> \bZ^{\Delta_+} @>>> N(\Delta)
@>>> 0\\
  && @VVV @VVV @VVV \\
0 @>>> \Sym^2\Lambda^\vee(\Delta') @>>> \bZ^{\Delta'^+} @>\psi>>
N(\Delta') @>>> 0\\
\end{CD}
\cooltag\end{equation}
where the first two vertical arrows are restrictions of functions.
The induced map $N(\Delta)\to N(\Delta')$ is surjective. In particular,
$M(\Delta')\to M(\Delta)$ is injective.
\end{Lemma}

\begin{proof}
Clear from definitions.
\end{proof}

\begin{Lemma}\label{stand}
Consider the standard partial order on $\Delta_+$:
$\alpha\ge\beta$ iff\/ $\alpha-\beta=\sum n_i\alpha_i$, where $n_i\ge 0$.
If\/ $\alpha>\beta$ then there exists a sequence
of simple roots $\alpha_{i_1},\ldots,\alpha_{i_p}$ such that
$\alpha=\beta+\alpha_{i_1}+\ldots+\alpha_{i_p}$
and $\beta+\alpha_{i_1}+\ldots+\alpha_{i_k}\in\Delta_+$ for any
$k=1,\ldots,p$.
\end{Lemma}

\begin{proof}
Let $\alpha-\beta=\sum_{i\in I}n_i\alpha_i$, where $n_i>0$ for any $i\in I$.
Recall that the inner product on $\Lambda$ is negative-definite.
Arguing by induction on $\sum n_i$, it suffices to prove that there exists $i\in I$
such that $\beta+\alpha_i$ or $\alpha-\alpha_i$ is a root. Suppose that the index
with this property does not exist. Then $(\beta,\alpha_i)\le0$ and $(\alpha,\alpha_i)\ge0$
for any $i\in I$.
Therefore, $(\alpha-\beta,\alpha-\beta)=(\alpha-\beta,\sum n_i\alpha_i)\ge0$. This is a contradiction.
\end{proof}

\begin{Lemma}\label{U(D4)}
There exists a unique decomposition
$D_4^+=F_1\coprod F_2\coprod F_3$, where each $F_i$ is a fourtuple
of pairwise orthogonal roots.
$M(D_4)\subset\bZ^{D_4^+}$ consists of functions
that are constant on each fourtuple $F_i$ and add up to~$0$.
\end{Lemma}

\begin{proof}
By \eqref{exsequn},
$M(D_4)\subset\bZ^{\cA_1(D_4)}$ consists of linear combinations
$$\sum\limits_{1\le i<j\le 4}a_{ij}[e_i-e_j]+b_{ij}[e_i+e_j]$$
such that $$\sum_{1\le i<j\le 4}a_{ij}(e_i-e_j)^2+b_{ij}(e_i+e_j)^2=0.$$
Looking at the coefficient at $e_ie_j$, we see that $a_{ij}=b_{ij}$ and the
remaining condition is that
$\sum\limits_{1\le i<j\le 4}a_{ij}(e_i^2+e_j^2)=0$. A short
calculation gives the claim.
\end{proof}

\begin{Theorem}\label{messy}
The map $\bigoplus\limits_{D_4\subset\Delta} M(D_4)\to M(\Delta)$ induced
by~\eqref{basicdiagram} is surjective.
In particular, the map $N(\Delta)\to\bigoplus\limits_{D_4\subset\Delta}N(D_4)$ is injective.
\end{Theorem}

\begin{proof}
We will prove by induction on the partial order $\ge$ of Lemma~\ref{stand} that
\begin{Claim}
For any $\alpha\in\Delta^T_+$ there exists a $D_4$ such that
$D_4^+=F_1\coprod F_2\coprod F_3$ as in Lemma~\ref{U(D4)}, $\alpha\in F_1$,
and all other roots of $F_1$ and $F_2$ are smaller than $\alpha$.
\end{Claim}

Assuming the Claim, let $I \subset M(\Delta)$ be the image of
$\bigoplus\limits_{D_4\subset\Delta} M(D_4)$.
By the Claim, for any $\alpha\in\Delta_+^T$, $I$ contains a function
$\alpha+\ldots$,
where all other
terms are less than~$\alpha$. Using induction
and Prop.~\ref{exsequn}, we conclude
that $I = M(\Delta)$.

Now we prove the Claim. We again argue by induction on $\ge$.
Let $\bD_4\subset\Delta$ be the unique $D_4$ with
$\Gamma(D_4)\subset\Gamma(\Delta)$.
Let $\alpha_0 \in \bD_4$ be the sum of the simple roots.
$\alpha_0$ is the smallest root in $\Delta_+^T$ and
one checks that $\alpha_0$ is greater than other positive roots of $\bD_4$
except for one root, which we include in $F_3$.
Thus the
Claim holds for $\alpha_0$.

Let $\alpha\in\Delta^T_+$.
Since $\alpha\ge\alpha_0$, by Lemma~\ref{stand} there
exists a simple root~$\gamma$  such that
$\alpha-\gamma$ is a root and $\alpha-\gamma\ge\alpha_0$, i.e., $\alpha-\gamma\in\Delta_+^T$.
By the induction hypothesis, the Claim is true for
$\alpha-\gamma$.
Let $D_4^+=F_1\coprod F_2\coprod F_3$ be the corresponding~$D_4$.
Let $r\in W(\Delta)$ be the reflection w.r.t.~$\gamma$. Then
$r(\alpha-\gamma)=\alpha$.
Consider $D_4'=r(D_4)$ and the fourtuples $F_i'$ obtained from $r(F_i)$ by
replacing
any negative root with its opposite.
Then $\alpha\in F_1'$ and we claim that any other
root $\pm r(\delta)$ in $F_1'$ or $F_2'$ is smaller than~$\alpha$.
Indeed, if $\delta\ne\gamma$ then $r(\delta)\in\Delta_+$
(because $r$ permutes roots in $\Delta_+\setminus\{\gamma\}$)
and is equal to $\delta+\gamma$, $\delta$, or $\delta-\gamma$.
In any case, $r(\delta)<\alpha$ because $\delta<\alpha-\gamma$.
If $\delta=\gamma$ then $r(\delta)=-\gamma$ and the corresponding
positive root in a fourtuple
is $\gamma$. And $\gamma<\alpha$.
\end{proof}

\begin{Remark}[\cite{Sekiguchi}]\label{D4class}
We will need the list of all possible $D_4\subset E_7$.
For any partition $\{ijkl\,|\,a\,|\,uv\}$,
$\{ij\,|\,kl\,|\,mn\,|\,b\}$, or $\{ab\,|\,cd\,|\,ijk\}$ of $\{1\dots7\}$,
let $D(ijkl,a)$, $D(ij,kl,mn)$, or $D(ab,cd)$ respectively,
be the $D_4$-subsystem  with fourtuples
$$D(ijkl,a)  =\left\{\{ij,kl,aij,akl\};\quad \{ik,jl,aik,ajl\};\quad \{il,jk,ail,ajk\}\right\}$$
$$D(ij,kl,mn)=\left\{\{ij, kl ,mn,b\};\quad\{ikm,iln,jkn,jlm\};\quad\{ikn,ilm,jkm,jln\}\right\}$$
$$D(ab,cd)   =\left\{\{i,jk, abi,cdi\};\quad\{j,ik, abj,cdj\};\quad\{k,ij, abk,cdk\}\right\}.$$
\end{Remark}

\begin{Theorem}\label{PsiCalculation}
We define $\psi(\Theta)$ as in \eqref{defrn}.
Let $\Theta=D_k\subset D_n$. Then $\psi(\Theta)=\psi(\Theta^\perp)=2\psi(A_{k-1})$ for any $A_{k-1}\subset D_k$.
In Table~\ref{perpsE7} we list
conjugacy classes of irreducible root subsystems in $E_6$ and $E_7$
and the relation between their $\psi$-images.

In a few cases $\psi(\Theta)$ is divisible in the lattice $N(\Delta)$:\quad
${1\over 2}\psi(D_k)\in N(D_n)$,\break
${1\over 4}\psi(D_k\times D_k)\in N(D_{2k})$,
${1\over 3}\psi(A_2^{\times 3})\in N(E_6)$,
${1\over 2}\psi(A_3^{\times 2})\in N(E_7)$,
${1\over 4}\psi(A_7)\in N(E_7)$.

\begin{table}[h]
\small\begin{tabular}{|c|c|c|c|c|c|}
  \hline
  $\Theta$ & $\Theta^\perp\subset E_6$  & $\psi$ & $f$ \cr
  \hline
  $A_1$& $A_5$&$\psi(A_5)=3\psi(A_1)$&$\sum\limits_{i=1}^6x_i(d+x_i)$              \cr
  $A_2$& $A_2^{\times2}$  & $2\psi(A_2)=\psi(A_2^{\times2})$ & $\sum\limits_{1\le i<j\le 3}(x_ix_j+x_{i+3}x_{j+3}) - d^2$\cr
  $A_3$& $A_1^{\times2}$  &$\psi(A_3)=\psi(A_1^{\times2})$&
  $x_1x_2-x_3x_4+x_5x_6-d(x_3+x_4)-d^2$\cr
  $A_4$& $A_1$ &$\psi(A_4)=2\psi(A_1)$&$x_1(d+x_1)-\sum\limits_{i=2}^6x_i(d+x_i)$              \cr
  $A_5$&$A_1$ &&\cr
  $D_4$& $\emptyset$ &$\psi(D_4)=0$&$\sum\limits_{i=2}^5x_i^2-x_6^2-(d+x_1)(x_1+3d+2x_6)$        \cr
  $D_5$& $\emptyset$ &$\psi(D_5)=0$&$3d^2+x_1^2+4x_1d-\sum\limits_{i=2}^6x_i^2$        \cr
  \hline
%
%
  \hline
  $\Theta$ & $\Theta^\perp\subset E_7$& $\psi$ & $f$ \cr
  \hline
  $A_1$  & $D_6$ &$\psi(D_6)=3\psi(A_1)$                  &$3d^2+x_1^2+4x_1d-\sum\limits_{i=2}^7x_i^2$\cr
  $A_2$& $A_5^-$  &$\psi(A_5^-)=2\psi(A_2)$&$3d^2+2d(x_1+x_2)-2x_1x_2+\sum\limits_{i=1}^2x_i^2-\sum\limits_{i=3}^7x_i^2$              \cr
  $A_3$  & \!$A_3'\!\times\!\! A_1$\! &$\psi(A_3)=\psi(A_3')$&$\sum\limits_{i=1}^3 x_i(d+x_i)-\sum\limits_{i=4}^7 x_i(d+x_i)$\cr
  $A_5^-$ &$A_2$&&\cr
  $D_4$  & $A_1^{\times3}$& $\psi(D_4)=\psi(A_1^{\times3})$&$\sum\limits_{i=2}^5x_i^2-\sum\limits_{i=6}^7x_i^2-(d+x_1)(x_1+3d+2\sum\limits_{i=6}^7x_i)$\cr
  $D_5$  & $A_1$ &$\psi(D_5)=2\psi(A_1)$                  &$3d^2+x_1^2+4x_1d-\sum\limits_{i=2}^7x_i^2-4dx_7-2x_7\sum\limits_{i=2}^6x_i$\cr
  $D_6$  & $A_1$&&\cr
  $E_6$  & $\emptyset$ & $\psi(E_6)=0$                    &$d^2-\sum\limits_{i=1}^6x_i^2+x_7^2$\cr
  $A_7$  & $\emptyset$ &&\cr
\hline
  $A_4$  & $A_2$ &$4\psi(A_4)\!=\!4\psi(A_2)\!+\!\psi(A_7)$&$5\sum\limits_{i=1}^2x_i(d+x_i)-3\sum\limits_{i=3}^7x_i(d+x_i)$\cr
  $A_5^+$& $A_1$ &$2\psi(A_5^+)\!=\!2\psi(A_1)\!+\!\psi(A_7)$&$3x_1(d+x_1)-\sum\limits_{i=2}^7x_i(d+x_i)$             \cr
  $A_6$  & $\emptyset$ & $4\psi(A_6)=3\psi(A_7)$ &$\sum\limits_{i=1}^7 x_i(d+x_i)$\cr
\hline
\end{tabular}
\smallskip
\caption{\small Root subsystems of $E_6$ and $E_7$. In the last three rows, $A_7$
is a unique subsystem of this type that contains both $\Theta$ and $\Theta^\perp$.}\label{perpsE7}
\end{table}
\end{Theorem}

\begin{proof}
The classification of irreducible subsystems is well-known \cite{Seki1,Sekiguchi}.
In each case it is easy to compute $\Theta^\perp$.
We will check only three equalities involving $\psi(\Theta)$,
leaving other cases to the reader.
Our proof is a routine calculation using Prop.~\ref{exsequn},
which identifies elements of $\bZ^{\Delta_+}$ trivial in $N(\Delta)$
with functions of the form $f(\alpha)$, where $f\in \Sym^2\Lambda^\vee(\Delta)$.
We provide $f$ in each case in Table~\ref{perpsE7}.

Let $\Theta=D_k\subset\Delta=D_n$. Write $D_I = \{ \pm \eps_i \pm \eps_j \ | \ i,j \in I \}$ for $I \subset \{1,\ldots,n\}$.
We can assume that $\Theta=D_{1,\ldots,k}$ and $\Theta^\perp=D_{k+1,\ldots,n}$.
Let $A_{k-1}\subset D_k$ be the standard subsystem of Lemma~\ref{anrem}.
Let $x_1,\ldots,x_{n}$ be standard coordinates on $\Lambda(D_{n})\subset\bZ^n$.
Consider $f=x_1^2+\ldots+x_k^2-x_{k+1}^2-\ldots-x_n^2$.
Then $f(\alpha)=2$ if $\alpha\in \Theta$, $f(\alpha)=-2$ if $\alpha\in\Theta^\perp$,
and $f(\alpha)=0$ otherwise. It follows that $\psi(\Theta)=\psi(\Theta^\perp)$.
Now consider $f=\sum\limits_{1\le i<j\le k}x_ix_j$.
Then $f(\alpha)=-1$ if $\alpha\in A_{k-1}$, $f(\alpha)=1$ if $\alpha\in D_k\setminus A_{k-1}$, and $f(\alpha)=0$ otherwise.
So  $\psi(D_{k})=2\psi(A_{k-1})$.

Let $\Theta=A_2\subset\Delta=E_6$. Then $A_2^\perp=A_2'\times A_2''$, where we can assume that
$A_2=\{123, 7, 456\}$,
$A_2'=\{45, 56, 46\}$, and
$A_2''=\{12, 23, 13\}$.
We denote by $\{d,x_1,\ldots,x_6\}$ the standard coordinates on $\Pic_6$
and their restrictions on $\Lambda(E_6)$.
Consider the function $f\in\Sym^2\Lambda^\vee(E_6)$
from the $A_2\subset E_6$ row of Table~\ref{perpsE7}.
Then $f(\alpha)=-1$ if $\alpha\in A_2'\times A_2''$,
$f(\alpha)=2$ if $\alpha\in A_2$, and $f(\alpha)=0$ otherwise.
It follows that $-\psi(A_2'\times A_2'')+2\psi(A_2)=0$.

Take the standard $A_7\subset E_7$ \eqref{standA7}.
It contains an $A_6$ with  roots $\alpha_{ij}$.
We use coordinates on $\Lambda(E_7)$ as above.
Take $f\in\Sym^2\Lambda^\vee(E_7)$ from the $A_6\subset E_7$ row of Table~\ref{perpsE7}.
Then $f(\alpha)=2$ for $\alpha\in A_6$,
$f(\alpha)=-6$ for $\alpha\in A_7\setminus A_6$, and $f(\alpha)=0$ otherwise.
It follows that $8\psi(A_6)-6\psi(A_7)=0$.

It remains to check divisibility of $\psi(\Theta)$.
In $D_n$, $\psi(D_k)=2\psi(A_{k-1})$.
In~$D_{2k}$, $\psi(D_k\times D_k)=\psi(D_k\times D_k^\perp)=2\psi(D_k)=4\psi(A_{k-1})$.
In $E_6$, $\psi(A_2^{\times 3})=\psi(A_2)+\psi(A_2^\perp)=3\psi(A_2)$.
In $E_7$, $\psi(A_3^{\times 2})=\psi(A_3)+\psi(A_3')=2\psi(A_3)$.
Finally, we have $4\psi(A_6)=3\psi(A_7)$, and therefore $\psi(A_7)$  is divisible by $4$.
\end{proof}

\section{Intrinsic Torus of $Y(\Delta)$ and KSBA cross-ratios}\label{intrtorus}

Now we identify $M(\Delta)$, combinatorially defined in \ref{weylg}, with units of $Y(\Delta)$.

\begin{Lemma}\label{haha} Let $Z_1,\dots,Z_k \subset \bG_m^l$ be
distinct irreducible Weil divisors defined by equations
$F_1,\dots,F_k \in \cO(\bG_m^l)$.
Then $Y:=\bG_m^l\setminus\cup Z_i$ is very affine,
with intrinsic torus $\bG_m^{k+l}$, and $\cO^*(Y)/k^*$ is generated by the $F_i$ and the coordinates
of~$\bG_m^l$.
\end{Lemma}

\begin{proof}
Let $W$ be the graph of \raise2.7pt\hbox{$\begin{CD}\bG_m^l@>{F_1,\ldots,F_k}>>\bA^k\end{CD}$}.
Then $W\cap \bG_m^{k+l}$ is closed in  $\bG_m^{k+l}$ and isomorphic to $Y$ via the first projection.
Therefore $Y$ is very affine.
Factoriality of $\bG_m^l$ implies that units of~$Y$ have form
$\chi F_1^{p_1}\ldots F_k^{p_k}$, where $\chi$ is a unit of $\bG_m^l$.
\end{proof}

\begin{Review}
We define $Y(D_n):=M_{0,n}$. Here $D_3=A_3$.
The Weyl group $W(D_n)$ acts on $M_{0,n}$ through its quotient $S_n$ (permuting marked points).
Fixing the first $3$ points, we identify $M_{0,n}$ with an open subset of~$\bA^{n-3}$:
a point $(z_1,\ldots,z_{n-3})\in\bA^{n-3}$ corresponds to
$n$ points in $\bP^1$ given by the columns of the matrix
$$
\left[\begin{matrix}
1& 0& 1& z_1  &\ldots&z_{n-3}\cr
0& 1& 1& 1&\ldots&1\cr
\end{matrix}\right].
$$
The boundary of $M_{0,n}$ in $\bA^{n-3}$
consists of the hyperplanes $z_i=0$, $z_i=1$, and $z_i=z_j$.
Therefore $M_{0,n}$ is an open subvariety of $\bG_m^{n-3}\subset \bA^{n-3}$
with divisorial boundary and so
is very affine, with intrinsic torus of rank $n(n-3)\over 2$
by Lemma~\ref{haha}.
\end{Review}

\begin{Review}
$Y(E_n)= Y^n$ is the moduli space of marked del Pezzo surfaces, $n=4,5,6,7,8$.
Here $E_4:=A_4$ and $E_5:=D_5$.
A \emph{marking} of $S$ is an
isometry $m:\,\Pic_n\to\Pic S$ such that
$m(K_n)=K_S$.
We denote $H=m(h)$ and $E_i=m(e_i)$.
Two marked surfaces $(S,m)$ and
$(S',m')$ are isomorphic if there exists an isomorphism
$f:\,S\to
S'$ such that $m=f^*\circ m'$.
$W(E_n)$ acts on~$Y^n$ by twisting markings: $w\cdot(S,m)=(S,m\circ w)$.

Let $X^n$ be the quotient of the open subset of
$(\bP^2)^n$ of points in linearly general position
by the free action of $\PGL_3$.
As above, we identify $X^n$
with an open subset of $\bA^{2n-8}$ by
taking $n$ points in $\bP^2$ given by columns of the matrix
\begin{equation}\label{matrix}
\left[\begin{matrix}
1& 0& 0 & 1& x_1&\ldots&x_{n-4}\cr
0& 1& 0 & 1& y_1&\ldots&y_{n-4}\cr
0& 0& 1 & 1& 1  &\ldots&1\cr
\end{matrix}\right].\cooltag
\end{equation}
For any $(S,m)\in Y^n$,
$\cO_S(H)$ gives a morphism $S\to\bP^2$ that blows down
$E_1,\ldots,E_n$ to $n$ points $p_1,\ldots,p_n\in\bP^2$
such that no $3$ points lie on a line, no $6$ points lie on a conic,
and no $8$ points lie on a  cubic curve singular at one of them~\cite{Dem}.
For any collection of points with these conditions,
their blow up is a del Pezzo surface \cite{Dem}.
Therefore, $Y^n$ is isomorphic to an open subset of $\bA^{2n-8}$ with a divisorial boundary~\cite[\S 4]{Seki1}
that includes the hyperplanes $x_i=0$ and $y_i=0$.
Therefore $Y^n$ is very affine by Lemma~\ref{haha}.

Blowing up the $n$ sections \eqref{matrix}
of the trivial bundle $Y^n\times\bP^2$
gives a universal family of del Pezzo surfaces $\cS\to Y^n$.
The natural morphism
\begin{equation}\label{mapP}
p:\,Y^{n+1}\to Y^n\cooltag
\end{equation}
given by blowing down~$E_{n+1}$
is induced by the projection $\bA^{2n-6}\to\bA^{2n-8}$ given by
truncating
the matrix~\eqref{matrix}. Therefore it factors
through an open embedding $i:\,Y^{n+1}\hookrightarrow\cS$
with a divisorial boundary
that we denote by $\cB$.
\begin{equation}\label{diagramEnEn+1}
\begin{matrix}
Y^{n+1}&\hookrightarrow&(\cS,\cB)\cr
&\searrow&\downarrow\cr
&&Y^n\cr
\end{matrix}\cooltag
\end{equation}
The fiber of $\cB$ over $[S]\in Y^n$
is identified with the union of $(-1)$-curves in~$S$
and (if $n=7$) the ramification curve of the double cover $S\to\bP^2$
given by $-K_S$. It is well-known (see e.g.~\cite{Dem}) that the number of irreducible
components in $\cB$ is equal to $10$ if $n=4$, $16$ if $n=5$, $27$ if $n=6$, and $57=56+1$ if $n=7$.
\end{Review}

\begin{Lemma}\label{grot}
\eqref{diagramEnEn+1} is $W(E_n)$-equivariant,
in particular $W(E_n)$ permutes irreducible components of $\cB_n$.
In the Grothendieck ring of $W(E_n)$-modules, one has
\begin{equation}\label{grotenen+1}
[\cO^*(Y^{n+1})/k^*)]=[\cO^*(Y^n)/k^*]+[\bZ^{\cB_n}]-[\Pic_n].\cooltag
\end{equation}
In particular, $\rk\cO^*(Y^n)/k^*=\rk M(E_n)$.
\end{Lemma}

\begin{proof}
We have an exact sequence of sheaves on ${\cS^n}$
$$0\to \cO^*_{\cS^n}\to i_*\left[\cO^*_{Y^{n+1}}\right]\arrow^{\val}\bigoplus_{k=1}^{\#\cB_n}F_k\to0,$$
where $F_k$ is a push-forward of the constant sheaf $\bZ$
on the $k$-th irreducible component of $\cB_n$ and
$\val$ is the order of vanishing along components of~$\cB_n$.
The long exact sequence of cohomology gives
$$0\to \cO^*(Y^{n})\to \cO^*(Y^{n+1})\to
\bZ^{\#\cB_n}\to\Pic_n\to0,$$
which implies \eqref{grotenen+1}.
The equality of ranks follows by counting and~\eqref{grotenen+1}.
\end{proof}

\begin{Definition} \label{mapfromcoxeterDn}\label{mapfromcoxeterEn}
Let $\cX(\Delta)$ be the complement of the Coxeter arrangement 
$$\cH =\{(\alpha=0) \ | \ \alpha \in \Delta_+ \} \subset \Lambda_k^{\vee}.$$
We define a map $\Psi:\,\cX(\Delta)\to Y(\Delta)$ as follows.

Let $\eps_1,\ldots,\eps_n$ \eqref{rootsdn} be the coordinates  on $\Lambda^\vee_k(D_n)$.
Consider the map
$$\Psi:\,\cX(D_n)\to M_{0,n},\quad \Psi(\eps_1,\ldots,\eps_n)=(\bP^1;\
\eps^2_1,\ldots,\eps^2_n).$$
The image is $M_{0,n}$ because $\eps_i=\pm \eps_j$ are precisely the root hyperplanes.

Let $q_1=h-3e_1,\ldots,q_n=h-3e_n$ be coordinates on  $\Lambda_k^\vee(E_n)$ defined using \eqref{rootsEn}.
We define $\Psi:\,\cX(E_n)\to Y^n$ as follows:
$\Psi(q_1,\ldots,q_n)$ is the blow up of $\bP^2$ in points
$$p_1=(q_1:q_1^3:1),\quad\ldots,\quad p_n=(q_n:q_n^3:1).$$
Our points lie on a cuspidal cubic curve (with a cusp at the infinity).
The blow-up is a del Pezzo surface \cite{Dem}:
the points are distinct $\Leftrightarrow$ $q_i\ne q_j$,
no three points lie on a line $\Leftrightarrow$ $q_i+q_j+q_k\ne0$,
no six points lie on a conic $\Leftrightarrow$
$q_{i_1}+\ldots+q_{i_6}\ne0$, and
no eight points lie on a cubic singular at one of them
$\Leftrightarrow$ $2q_{i_1}+q_{i_2}\ldots+q_{i_8}\ne0$.
One can check that these equations are precisely the root hyperplanes.
\end{Definition}

\begin{Theorem}\label{pullbackofunits}
There exists a $W(\Delta)$-equivariant commutative diagram
$$
\begin{CD}
M(\Delta) @>{\psi^{\vee}} >> \bZ^{\Delta_+} \\
@V{f}VV                              @V{h}VV \\
\cO^*(Y(\Delta))/k^* @>{\Psi^*}>> \cO^*(\cX(\Delta))/k^*
\end{CD}
$$
where $h$ and $f$ are isomorphisms, and $h([\alpha])$ is equal to $\alpha$ as a
function on $\Lambda^{\vee}_k$.
\end{Theorem}

\begin{proof}
It is clear that $\cO^*(\cX)/k^*=\bZ^{\Delta_+}$.
Since the complement to a subarrangement given by simple roots
is an algebraic torus $\bG_m^n$, it follows from Lemma~\ref{haha} that
$\cX$ is very affine with the intrinsic torus
$T_\cX=\Hom(\bZ^{\Delta_+},\bG_m)$.

Let $\Psi^*:\,\cO^*(Y(\Delta))/k^*\to\cO^*(\cX)/k^*=\bZ^{\Delta_+}$ be the pull-back of
units.
By Lemma~\ref{grot}, $\rk M(\Delta)=\rk \cO^*(Y(\Delta))/k^*$.
Therefore, it suffices to prove that any
$m\in M(\Delta)\subset\bZ^{\Delta_+}$ is a pull-back of a unit of~$Y(\Delta)$.
By Th.~\ref{messy} it suffices to prove that, for any $D_4\subset\Delta$, any
$m\in M(D_4)\subset M(\Delta)$ is a pull-back of a unit.
By~\ref{DynkinTh}, it suffices to prove this for any fixed $D_4\subset\Delta$.
Using the action of $W(D_4)$ and Lemma~\ref{U(D4)}, it suffices to
construct as a pull-back of a unit the function in $\bZ^{D_4^+}\subset\bZ^{\Delta_+}$ that is equal to $0$, $1$,
and $-1$
on three fourtuples of orthogonal $A_1$'s.

Suppose $\Delta=D_n$.
Consider a forgetful map $M_{0,n}\to M_{0,4}$.
A pull-back of a unit of $M_{0,4}$ (i.e.~a cross-ratio of points $p_1,p_2,p_3,p_4\in\bP^1$)
to $\cX(D_n)$ is equal to
$${(\eps_1^2-\eps_2^2)(\eps_3^2-\eps_4^2)\over(\eps_1^2-\eps_3^2)(\eps_2^2-\eps_4^2)}=
{(\eps_1-\eps_2)(\eps_1+\eps_2)(\eps_3-\eps_4)(\eps_3+\eps_4)\over
(\eps_1-\eps_3)(\eps_1+\eps_3)(\eps_2-\eps_4)(\eps_2+\eps_4)}.$$
In the additive notation, this function is equal to $1$ on the
four orthogonal roots of the numerator and to $-1$ on
the four orthogonal roots of the denominator.
By \eqref{rootsdn}, these fourtuples belong to a $D_4$.

Suppose $\Delta=E_n$. A cross-ratio of projections of points
$p_1,p_2,p_3,p_4$ from $p_5$
is a unit of $Y(\Delta)$.
This is a geometric cross-ratio of type I, where (in the notation of Def.~\ref{typeIdef}),
$L_5$ is an exceptional divisor over $p_5$ and $L_1,\ldots,L_4$ are strict transforms
of lines connecting $p_1,\ldots,p_4$ with $p_5$.
The pull-back of this unit to $\cX(E_n)$ is
\begin{equation}\label{cool}
\small\frac{
\left|\begin{matrix}
q_1&q_2&q_5\cr
q_1^3&q_2^3&q_5^3\cr
1&1&1\cr
\end{matrix}\right|
\cdot
\left|\begin{matrix}
q_3&q_4&q_5\cr
q_3^3&q_4^3&q_5^3\cr
1&1&1\cr
\end{matrix}\right|
}{
\left|\begin{matrix}
q_1&q_3&q_5\cr
q_1^3&q_3^3&q_5^3\cr
1&1&1\cr
\end{matrix}\right|
\cdot
\left|\begin{matrix}
q_2&q_4&q_5\cr
q_2^3&q_4^3&q_5^3\cr
1&1&1\cr
\end{matrix}\right|
}.
\cooltag\end{equation}
Using the identity
$
\left|\small\begin{matrix}
a&b&c\cr
a^3&b^3&c^3\cr
1&1&1\cr
\end{matrix}\right|=(a-b)(b-c)(c-a)(a+b+c)
$,
\eqref{cool} is equal to
$$\frac{(q_3-q_4)(q_1-q_2)(q_3+q_4+q_5)(q_1+q_2+q_5)
}{(q_1-q_3)(q_2-q_4)(q_1+q_3+q_5)(q_2+q_4+q_5)}=
\frac{\alpha_{34}\alpha_{12}\alpha_{345}\alpha_{125}
}{\alpha_{13}\alpha_{24}\alpha_{135}\alpha_{245}}.$$
By \eqref{D4class}, these fourtuples of orthogonal roots belong to a
$D_4$.
\end{proof}

\begin{Corollary}\label{D4inclpullback}
For $D_4 \subset \Delta$, the composition of morphisms
$$\cO^*(M_{0,4})/k^*=M(D_4) \mathop{\hookrightarrow}^i M(\Delta)=\cO^*(Y(\Delta))/k^*$$
is the pull-back of units
induced either by a forgetful map $M_{0,n}\to M_{0,4}$ (if\/ $\Delta=D_n$)
or by a KSBA cross-ratio map of type~I (if\/ $\Delta=E_n$).
\end{Corollary}

\section{Simplicial Complex $\cR(\Delta)$ of Boundary Divisors}

\begin{Notation}
Let $\oY(D_n): = \oM_{0,n}$ and let $\oY(E_6)$ be
Naruki's space of cubic surfaces.
For $\ch k \neq 2$, let $\oY(E_7)$ be the coarse moduli
space for $\tcM_3^{(2)}$ .
Note that by \cite{DO}, the coarse moduli space $\tM_3^{(2)}\setminus B$ is smooth and
isomorphic to $Y(E_7)$. In what follows, whenever we refer to $\oY(E_7)$ we implicitly assume $\ch k \neq 2$.
\end{Notation}

\begin{Lemma}\label{RDn}\label{Dnclass}\label{anrem}
Let $\bZ^n$ be the standard $\bZ$-lattice
with the
product $\eps_i\cdot\eps_j=-\delta_{ij}$. Then
\begin{equation}\label{rootsdn}
D_n=\{\alpha\in \bZ^n\ |\ \alpha^2=-2\}=
\{\pm\eps_i\pm\eps_j\ |\ 1\le i<j\le n\}.\cooltag
\end{equation}
The Weyl group $W(D_n)\simeq S_n\ltimes (\bZ_2)^{n-1}$, where $(\bZ_2)^{n-1}$ acts by even sign changes.
Let $N_n := \{1,2,\dots n\}$.
$\cR(D_n)$ is equivariantly isomorphic ($\bZ_2^{n-1}$ acts trivially)  to
\begin{itemize}
\item the set of boundary divisors of $\oM_{0,n}=\oY(D_n)$;
\item the set of bipartitions $I\coprod I^c=N_n$, where $2\le|I|\le|I^c|$;
\item the set of subsystems of type $A_k$ for $n>2k+2$ and $A_l\times A_l$ if $n=2l+2$
of the `standard' $A_{n-1}\subset D_n$ with roots $\eps_i-\eps_j$, $i\neq j\in N_n$.
\end{itemize}
Simplices are subsets of divisors with non-empty intersection,
pairwise compatible bipartitions,
and pairwise orthogonal or nested subsystems, respectively.
Note also that if we fix $N_{n-1} \subset N_{n}$ then under the action of
$S_{n-1} \subset S_{n}$, $\partial \oM_{0,N_{n}}$ is equivariantly
isomorphic to the collection of subsets of $N_{n-1}$ of cardinality
between $2$ and $n-2$.
\end{Lemma}

\begin{proof}
To any  bipartition  we assign
a boundary divisor $\delta_{I,I^c}\subset\oM_{0,n}$ (see e.g.~\cite{Keel}),
a subsystem $D_I=\{\pm\eps_i\pm\eps_j\,|\,i,j\in I\}\subset D_n$ of type $D_{|I|}$
(if $|I|=|I^c|$ then we assign $D_I\times D_{I^c}$), and
the intersection of this subsystem of $D_n$ with $A_{n-1}$.
\end{proof}

\begin{figure}[htbp]
  \includegraphics[width=4in]{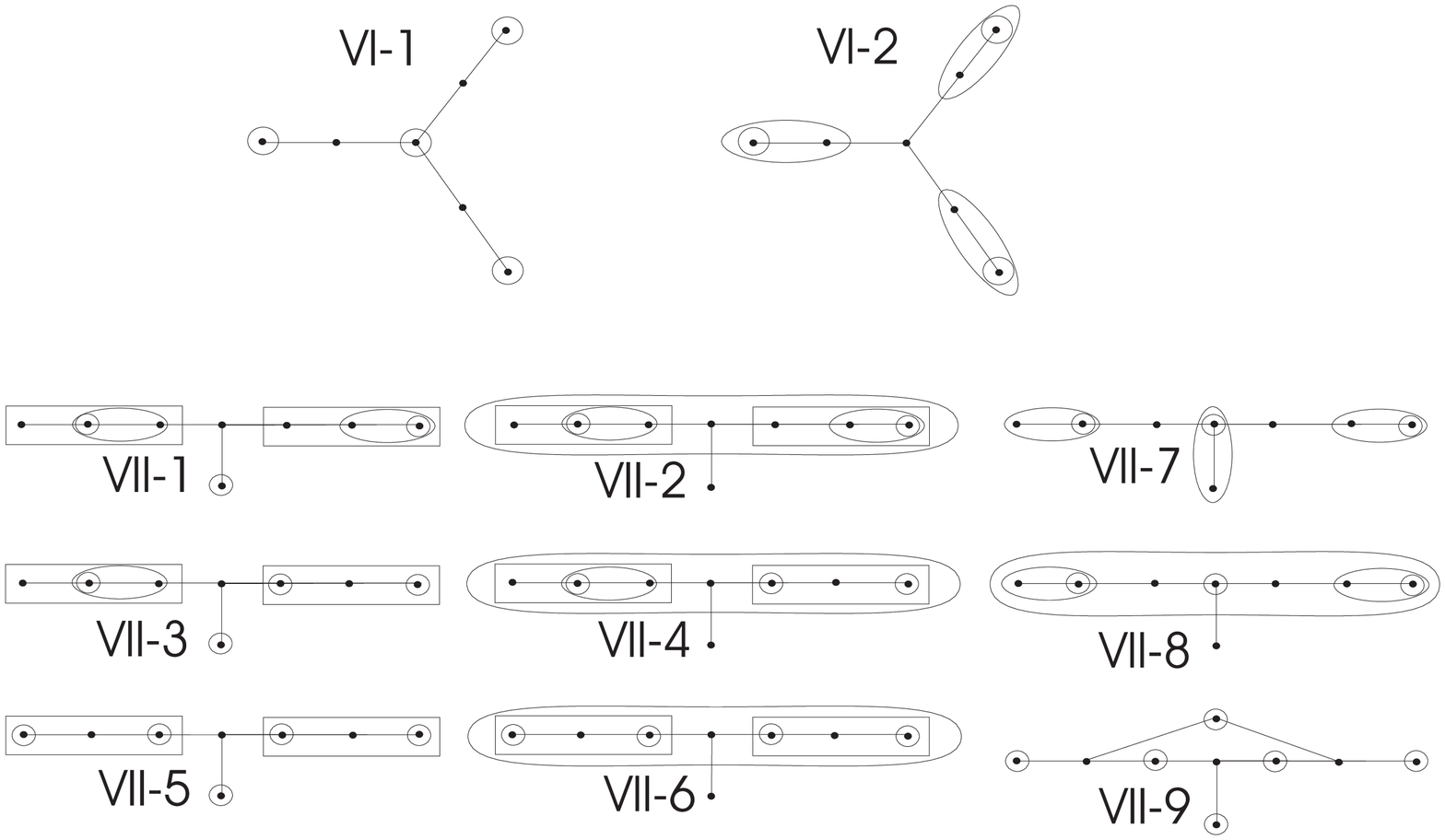}\\
  \caption{\small Maximal Simplices of $\cR(E_6)$ and $\cR(E_7)$}\label{conesofF(Delta)}
\end{figure}

\begin{Review}\textsc{Complex $\cR(E_n)$.}
Suppose that $\sigma=\{\Theta_1,\ldots,\Theta_k\}\in\cR(E_n)$
has the property that (for some choice of simple roots) $\Gamma(\Theta_i)\subset\tilde\Gamma(E_n)$
for any $i$, where $\tilde\Gamma(E_n)$ is the affine Dynkin diagram.
Then we can `tube' $\sigma$ as in \cite{Williams}
by drawing $k$ `tubes' encircling $\Gamma(\Theta_1),\ldots,\Gamma(\Theta_k)$.
To avoid confusion, we box reducible subsystems
$A_2\times A_2\times A_2$ and $A_3\times A_3$ instead of circling.
We define a few simplices on Fig.~\ref{conesofF(Delta)}, where
$VII_9$ is tubed on the diagram $\Gamma'(E_7)$ obtained by adding
the root $\alpha_{167}$ to~$\tilde\Gamma(E_7)$.
\end{Review}

\begin{Proposition}\label{classification}
Points of $\cR(E_n)$ are $W(E_n)$-equivariantly identified with boundary divisors of
$\oY(E_n)$. Its simplices correspond to subsets
of intersecting divisors.
Any simplex of $\cR(E_n)$ is $W(E_n)$-equivalent to a face
of a simplex in Figure~\ref{conesofF(Delta)}.
\end{Proposition}

\begin{proof}
The identification of vertices of $\cR(E_6)$
with divisors of $\oY(E_6)$ and the description
of maximal simplices are explicit in \cite{Naruki}.

In Th.~\ref{incidencies} we
parametrized boundary divisors of~$\tcM_3^{(2)}$
by points of~$\cR(E_7)$.
We also identified $D(A_2)$, $D(A_3\times A_3)$, and $D(A_7)$
as the moduli spaces
$\oM_{0,4} \times \oM_{0,7}$, $\oM_{0,5} \times \oM_{0,5} \times \oM_{0,4}$,
and $\oM_{0,8}$.
The full boundary has simple normal crossings, and restricts to
the boundary on each of these spaces.
We first check that divisors $D(\Gamma)$ and $D(\Delta)$ intersect if and only if
$\Gamma$ and $\Delta$ are orthogonal or nested.
If both the divisors have type $D(A_1)$, see
Lemma~\ref{deltairr}. If $\Gamma=A_1$ and $\Delta=A_2$,
or if $\Gamma=\Delta=A_2$, this is clear from
the proof of Lemma~\ref{delta1}.
$A_3\times A_3$ divisors are pairwise disjoint by Lemma~\ref{E}.

\begin{Lemma} \label{meetingH}
If $\Delta=A_7$ then $D(\Gamma) \cap D(\Delta) \neq \emptyset$
iff $\Gamma \subset \Delta$.
Subsystems $A_1, A_2, \break A_3\times A_3\subset A_7$
correspond to boundary divisors of $\oM_{0,8}\simeq D(A_7)$
as in Lemma~\ref{RDn}.
\end{Lemma}

\begin{proof}
$A_7$ divisors are pairwise disjoint by Lemma~\ref{H}.
By the proof of Lemma ~\ref{H}, $\delta_{\irr,\alpha} \cap \cH_{Q}
\neq \emptyset$ iff $Q(\alpha) =1$.
And $Q(\alpha) =1$ iff $A_1(\alpha) \subset A_7(Q)$.
Now consider $\cH_{Q} \cap \delta_{1,W} \neq \emptyset$. One
checks using the description of $Q$ in the proof of Lemma~\ref{H}
that $Q|_W =1$, or equivalently $A_2(W) \subset A_7(Q)$.
The case $\Gamma=A_3\times A_3$ follows from Lemma~\ref{H},
where it is shown that $\cT^{(2)}$ is the transverse intersection of
two $\ocH^{(2)}_{Q}$ that correspond to two $A_7$'s intersecting in an $A_3\times A_3$.
For each type of $\Gamma$ there is only one corresponding
divisor of $\tcH$: it is clear from Lemma~\ref{tm3nc}
that for $\delta_{\irr}$
this is of type $B_2$, for $\delta_1$ of type $B_3$, and
for $\cE$ of type $B_4$.
\end{proof}

Let $\Gamma=A_1$ or $A_2$ and $\Delta=A_3\times A_3$. Let $b:\,\tcM_3^{(2)}\to\ocM_3^{(2)}$.
Then $b(D(\Delta))$ is~$\cT^{(2)}$, the self-intersection
of $\delta_{irr,\alpha}$ for $\alpha=\Delta^{\perp}$.
$b(D(\Gamma))$ either contains $\cT^{(2)}$
(in which case $\Gamma=A_1=\Delta^{\perp}$)
or intersects it along a divisor.
But then $\Gamma\subset\Delta$ as only cycles in $\Delta$ can
be vanishing cycles for degenerations of a generic curve in $\cT^{(2)}$.

Now for each divisor we consider the simplicial
complex for its boundary. We check
that these links agree with the corresponding links of~$\cR(E_7)$.
For $A_1$ this follows from Lemma~\ref{deltairr}.
For the remaining divisors, it follows from the well-known fact (Lemma~\ref{RDn})
that a subset of divisors of $\oM_{0,n}$ (and thus of products $\oM_{0,n_1}\times\ldots\times\oM_{0,n_k}$)
has a non-empty intersection iff they intersect pairwise.
\end{proof}

\section{Fan of Real Components $\cG(\Delta)$}\label{realcomp}

\begin{Review}\textsc{Motivation.}
Fix $\Delta\ne E_8$.
Before proving that $\cF(\Delta)$ is a strictly simplicial fan supported on $\cA(Y(\Delta))$, 
we introduce a remarkable fan $\cG(\Delta)\subset N_\bQ$ which contains~$\cF(\Delta)$ but is much more
symmetric: its maximal cones are permuted by~$W$.
The origins of the fan are geometric:
$Y$ is defined over $\bR$ and the set $Y(\bR)$
is divided into connected components, which we call~cells.
The Weyl group acts transitively on the set of cells\footnote{This is easy for $D_n$ and for $E_n$
this follows from the classical
result that real del Pezzo surfaces of the same topological type
are connected in the moduli space, see \cite{KollarRDP} for the modern exposition.
Points of $Y(\bR)$ correspond to real del Pezzo surfaces
that are blow-ups of $\bP^2$ in $n$ real points.}.
The closure of each cell in $\oY(\Delta)$ is
diffeomorphic (as manifold with corners) to a polytope, which we denote~$P(\Delta)$.
This is clear locally near the boundary
(since we will prove that $\oY(\Delta)$ has normal crossings),
but the additional claim is that the closure is contractible.
$P(D_n)$ is identified with the associahedron in \cite{Devadoss}, and $P(E_n)$
is described in \cite{SY,S2} (for~$E_7$ without the proof).

Faces of $P(\Delta)$ correspond to the
boundary divisors of $\oY(\Delta)$ which meet the closure of the cell
(the face being the intersection of the boundary divisor with
the closure of the cell). This gives a natural way of dividing
up the boundary divisors, and we define the fan $\cG(\Delta)$ by declaring
(the convex hull of) a collection of rays (from $\cF \subset N_\bQ$)
to be a cone if the corresponding boundary divisors all
meet a single closed cell, i.e., they correspond to faces of $P(\Delta)$.
Since the zero strata of $\oY$ are real points, each cone of
$\cF^{\partial \oY}$ is a cone of $\cG(\Delta)$.
The remarkable fact is that the
number of faces of $P(\Delta)$ is equal to the rank of $N(\Delta)$ (except in
the $E_7$ case, where the rank is one larger). This raises the
hope that such cones form a maximal dimensional, strictly
simplicial fan, or equivalently that for each pair of a closed cell and one of its
maximal faces, one can find a unit (necessarily unique up to scaling) which vanishes
to order one along the face (i.e. along the corresponding boundary divisor) and
is generically
regular and non-vanishing on the other faces. This indeed turns out to be the case, and 
in fact the units are $D_4$ units. Our initial inspiration
for this fan comes from \cite{Brown}, where the units are given in the $D_n$
case. 

We will not prove these statements and will not use this geometric description.
It will suffice for our purposes to define $\cG(\Delta)$ combinatorially.
\end{Review}

Recall that any collection of $A_1$'s in $\Delta$ gives rise to a graph,
with vertices indexed by $A_1$'s and where vertices are connected by
an edge iff the corresponding $A_1$'s
are not orthogonal. For example, a collection of simple roots gives the Dynkin diagram,
adding the lowest root gives the affine Dynkin diagram, etc.

\begin{figure}[htbp]
  \includegraphics[width=0.9\hsize]{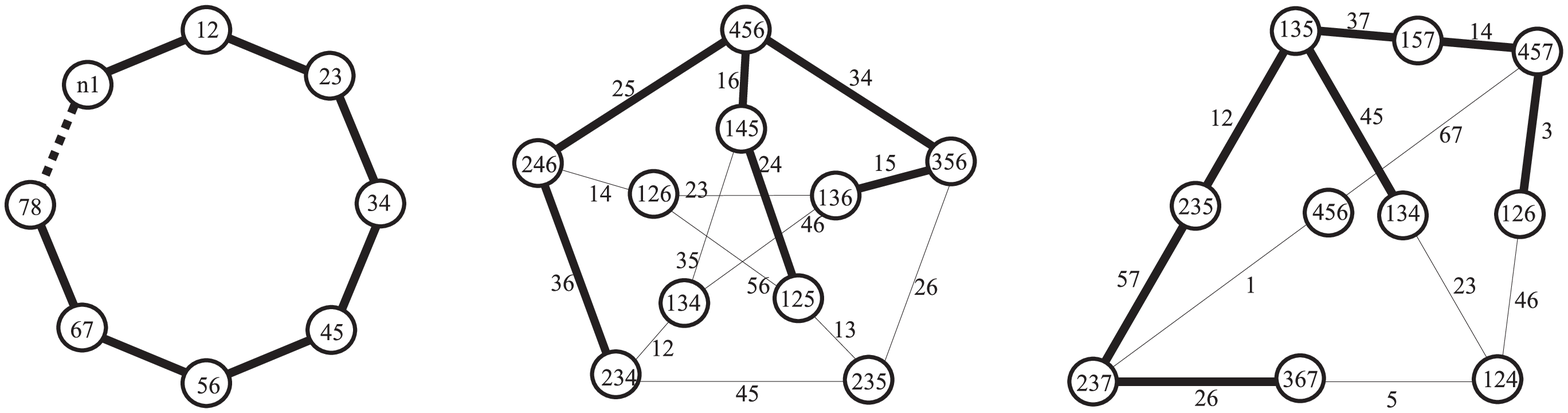}
  \caption{\small Sekiguchi's $n$-gon, pentadiagram, and tetradiagram with highlighted affine Dynkin diagrams.
  We mark each edge of the last two diagrams by a unique root that forms an $A_2$
  with vertices of the edge.}\label{pentatetra}
\end{figure}

\begin{Definition}
An {\em $n$-gon diagram}\/ is a set of $A_1$'s in $A_{n-1}$ which form
an $n$-gon.

A {\em pentadiagram}\/ is a collection of $10$ $A_1$'s in $E_6$
which form the Petersen graph.

A {\em tetradiagram}\/ is a collection of $10$ $A_1$'s
in $E_7$ which form the graph with $10$ vertices given by the $4$ vertices of a tetrahedron, together
with the $6$ midpoints of the $4$ edges, with edges the
$12$ half edges of the tetrahedron, see Figure~\ref{pentatetra}
\end{Definition}

\begin{Theorem}[Sekiguchi]\label{SekiSymmetries} For $\Delta = D_n$, $E_6$, and $E_7$,
the Weyl group acts
transitively on the set of $n$-gons, pentadiagrams,
and tetradiagrams, respectively. The normalizer of a diagram
is respectively the dihedral group $D_{2n}$,
$S_5$, and $S_4 \times (\bZ/2\bZ)$.
Its action on the diagram is  a natural one:
the action of the dihedral group on the $n$-gon,
the~action of $S_5$
on the Petersen graph $\cR(D_5)$ from Lemma~\ref{RDn},
and the action of~$S_4$ permuting vertices of the tetrahedron
($\bZ/2\bZ$ acts trivially).
\end{Theorem}

\begin{Lemma}\label{extraA1}
For each tetradiagram there is a unique $A_1$ which is
perpendicular to the $A_1^{\times 6}$ formed by the midpoints
of the edges of the tetrahedron.
\end{Lemma}

\begin{proof}
Uniqueness follows from $\rk E_7=7$.
For Figure~\ref{pentatetra}, the extra root is $\alpha_{247}$.
\end{proof}

\begin{Definition}
We define rays of $\cG$ to be the same as rays of $\cF$.
We define a maximal cone in $\cG$ for each $n$-gon, pentadiagram,
and tetradiagram, respectively.
We will check that they are indeed cones, in fact strictly simplicial ones, in Th.~\ref{dualbasis}.

\smallskip

$(D_n)$
We label an edge from $\eps_i-\eps_j$ to $\eps_j-\eps_k$ by $j$ --
so the diagram is now an $n$-gon with edges labeled by
the elements of $N:=\{1,\ldots,n\}$.
Each non-edge chord defines
a bipartition $N = I \coprod I^c$, and so a ray of $\cF(D_n)$. A collection
of rays forms a cone in $\cG(D_n)$ iff the corresponding chords do not intersect in the interior of the $n$-gon.
The number of rays is equal to the number of chords of an $n$-gon,
i.e., ${n(n-3)\over2}=\dim N(D_n)$.

\smallskip

$(E_6)$
A pentadiagram gives rise to $10$ $A_1$'s in $E_6$ (the vertices)
and $5$ subsystems of
type $A_2 \times A_2 \times A_2$ --- subsystems of
this type whose Dynkin diagram is a subdiagram of the
pentadiagram.
A collection of rays forms a cone in $\cG(E_6)$  iff they all come
from a single pentadiagram.
The number of rays is equal to $15=\dim N(E_6)$.

\smallskip

$(E_7)$
A tetradiagram gives rise to $10$ $A_1$'s in $E_7$ (the vertices), $12$ $A_2$'s (the edges),
$3+6$ subsystems of type
$A_3 \times A_3$ and $3$ $A_7$'s -- all the subsystems of these types whose
Dynkin diagrams are subdiagrams of the tetradiagram.
The number of rays of the corresponding cone in $\cG(E_7)$ is  $34$, i.e.~$\dim N(E_7)-1$.
Let $\cG'(E_7)$ be the fan obtained by adding to each maximal cone of~$\cG(E_7)$ the extra ray of Lemma~\ref{extraA1}.
\end{Definition}

\begin{Lemma} \label{subfan}
Each cone of $\cF(\Delta)$ is
a cone of $\cG(\Delta)$.
\end{Lemma}

\begin{proof}
For $D_n$: it is easy to check by induction that any set of compatible irreducible subsystems
of $A_{n-1}$ can be simultaneously tubed on the Dynkin diagram.
In particular it can be tubed on the affine Dynkin diagram, i.e.~on the $n$-gon.

For $E_n$: By Prop.~\ref{classification}, any maximal simplex of $\cR(E_n)$
can be tubed on the affine Dynkin diagram $\tilde\Gamma(E_n)$ (except $VII_9$).
But $\tilde\Gamma(E_n)$ is a subgraph of the pentadiagram (resp.~tetradiagram), see Fig.~\ref{pentatetra}.
Finally, $VII_9$ is a simplex of $6$ orthogonal $A_1$'s which can be realized by the midpoints
of the edges of the tetrahedron.
\end{proof}

\begin{Theorem} \label{dualbasis}
Let $R$ be a ray
of a maximal cone $\sigma \in \cG$ ($\cG'$ for $E_7$) defined above.
There exists a unique $D_4$ and a unique pair of $4$-tuples
$F_1 \coprod F_2 \subset \cA_1(D_4)$ as in Lemma~\ref{U(D4)}
such that $u(F_1,F_2)$ is equal to~$1$
on the first lattice point along $R$ and is equal to~$0$ on other rays of $\sigma$.
Here $u(F_1,F_2)\in M(D_4)$ is the function that
is equal to $1$ on $F_1$, $-1$ on $F_2$, and $0$ on $F_3$.
In particular, $\sigma$ is a strictly simplicial cone.
\end{Theorem}

\begin{proof}
By symmetry we can assume that $\sigma$ is given by Figure~\ref{pentatetra}.
Strict simpliciality follows from the previous statement,
since the number of rays  is equal to $\dim N$.
Uniqueness follows as well.
It clearly suffices to prove that
$u(F_1,F_2)$ is equal to~$1$ on some point of $R$, and so by Th.~\ref{PsiCalculation},
it will suffice to prove that
\begin{equation}\label{formulaforU}
u(F_1,F_2)(\Theta)=|F_1\cap\Theta|-|F_2\cap\Theta|=d,\cooltag
\end{equation}
where $\Theta\in\cR(\Delta)$ gives the ray $R$
and $d|\psi(\Theta)$ in $N(\Delta)$ ($d$ is given in Th.~\ref{PsiCalculation}).

{\em Proof for $D_n$.} Here rays are given by chords $R$ of the $n$-gon:
induced bipartitions $N=I\coprod I^c$ correspond to subsystems
$\Theta_R=D_I\times D_{I^c}\subset D_n$.
Note that $\psi(\Theta_R)\in R$ and is divisible by $4$ by Th.~\ref{PsiCalculation}.
Up to dihedral symmetry, $R$ corresponds to the bipartition $\{1,\ldots,k\}\coprod\{k+1,\ldots,n\}$.
It is illustrated in Figure~\ref{crazychordes}.
We take the $D_4=D_{k,k+1,n,1}$ and fourtuples $F_1=\{\eps_1\pm\eps_k,\ \eps_{k+1}\pm\eps_n\}$,
$F_2=\{\eps_1\pm\eps_{k+1},\ \eps_n\pm\eps_k\}$.
Note that $F_1\subset \Theta_R$ and $F_2\cap\Theta_R=\emptyset$.
So $u(F_1,F_2)$ is equal to $4$ on $\psi(\Theta_R)$.

\begin{figure}[htbp]
  \includegraphics[width=\hsize]{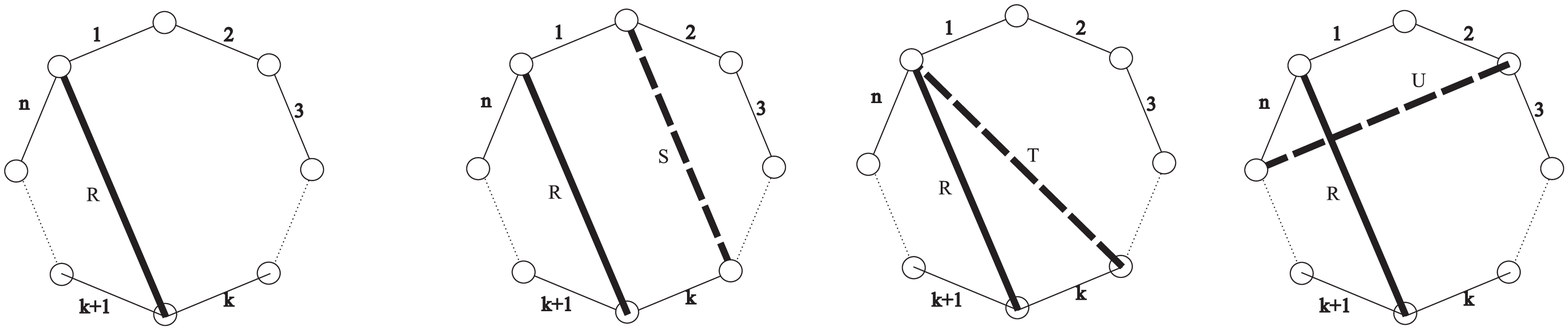}
  \caption{\small}\label{crazychordes}
\end{figure}

It remains to prove that $u(F_1,F_2)$ vanishes on other rays of $\sigma$.
There are three cases illustrated in Figure~\ref{crazychordes}.
If $S$ is a chord that does not intersect $R$
(including endpoints), then $F_1,F_2 \subset \Theta_S$.
And so $u(F_1,F_2)$ is equal to $4-4=0$ on $\psi(\Theta_S)$.
The next chord $T$ intersects $R$ at the endpoint. Then
$|F_1\cap\Theta_T|=|F_2\cap\Theta_T|=2$
and $u(F_1,F_2)$ is equal to $2-2=0$ on $\psi(\Theta_T)$.
Finally, $U$ is a chord that crosses $R$. Then
$F_1\cap\Theta_T=F_2\cap\Theta_T=\emptyset$
and $u(F_1,F_2)$ is equal to $0-0=0$ on $\psi(\Theta_U)$.

\smallskip

{\em Proof for $E_6$.} By Th.~\ref{SekiSymmetries}, the subgroup of $W$ normalizing
the maximal cone of $\cG$ is~$S_5$.
There are two $S_5$ orbits on rays, corresponding
to the type of the subsystems, either $\Theta = A_1$
or $\Theta = A_2^{\times 3}$.
Note that we can read all their roots from Figure~\ref{pentatetra}.

$\Theta = A_1$. By symmetry we can assume $\Theta =456$.
$D_4$ is determined as follows. It is clear that $\Theta\in F_1$
and no other vertex of the pentadiagram is in $F_1$ or~$F_2$.
Next, $\Theta$ is contained in three
subsystems $\Gamma$ of type $A_2^{\times 3}$ in the pentadiagram.
Each has form
$E \times A_2 \times A_2$, where
$E\simeq A_2$ corresponds to an edge in the diagram containing the
vertex $\Theta$. $\psi(E)$ is a lattice point on
the ray $\psi(\Gamma)$ by Th.~\ref{PsiCalculation}, and so $u(F_1,F_2)$ must vanish on $E$.
The only possibility is that the root on the edge $E$ in Fig.~\ref{pentatetra} is in $F_2$.
In our case we find that $456\in F_1$ and $16,34,25\in F_2$.
Inspecting the list of $D_4$'s in Remark~\ref{D4class}, we see that
the only possibility is $D(16,34,25)$ and
$$
F_1=\{145,123,246,356\},\quad F_2=\{7,16,34,25\}.
$$
It is clear that $u(F_1,F_2)$ vanishes on the two remaining $A_2^{\times3}$
because they contain only vertices of the pentadiagram other than $145$
and edges other than $16,34,25$.

$\Theta=A_2^{\times 3}$. By symmetry
$\Theta$ is given by edges $25$, $46$, and $13$ of Figure~\ref{pentatetra}.
The~$D_4$ and fourtuples are defined as follows.
No vertex of the pentadiagram is in $F_1$ or~$F_2$.
But since $u(F_1,F_2)$ has valuation $3$ on $\Theta$, all edges of $\Theta$ have to be in~$F_1$.
It follows from Remark~\ref{D4class} that
$D_4 = D(13,25,46)$ with fourtuples
$$
F_1=\{7,13,25,46\},\quad F_2=\{124,156,236,345\}.
$$
Since neither $F_1$ nor $F_2$ contain other edges,
$u(F_1,F_2)$ vanishes on the $4$ remaining~$A_2^{\times3}$.

\smallskip

{\em Proof for $\cG'(E_7)$.}
In Figure~\ref{supertetrafig} we plot all roots of all rays of~$\cG'(E_7)$.
There are three $S_4$-orbits with  $\Theta = A_1$: vertices, midpoints of edges,
and the special root $247$ of Lemma~\ref{extraA1}.
There is one orbit with $\Theta=A_2$ (half-edges).
There are two orbits with $\Theta = A_3^{\times2}$. The first type corresponds
to a pair of opposite edges, e.g.~$135,235,237$ and $124,126,457$.
We attach to each midpoint the only ``extra'' root in $A_3$ not contained in its $A_1$'s or $A_2$'s,
e.g.~we attach $137$ to $235$. The second type corresponds to a pair of vertices,
e.g.~$235,135,157$ and $367,124,126$ correspond to a pair $(135,124)$.
We attach extra roots in $A_3$'s to midpoints (and align them in the direction of an edge
connecting midpoints).
Finally there is one orbit with $\Theta=A_7$: it corresponds to a pair of opposite edges
(remove them from the tetrahedron to get the affine Dynkin diagram of $A_7$).
Each $A_7$ contains $4$ extra roots not contained in $A_1$, $A_2$, or $A_3$'s
and we plot them near the corresponding pair of opposite edges.

\begin{figure}[htbp]
  \includegraphics[width=\hsize]{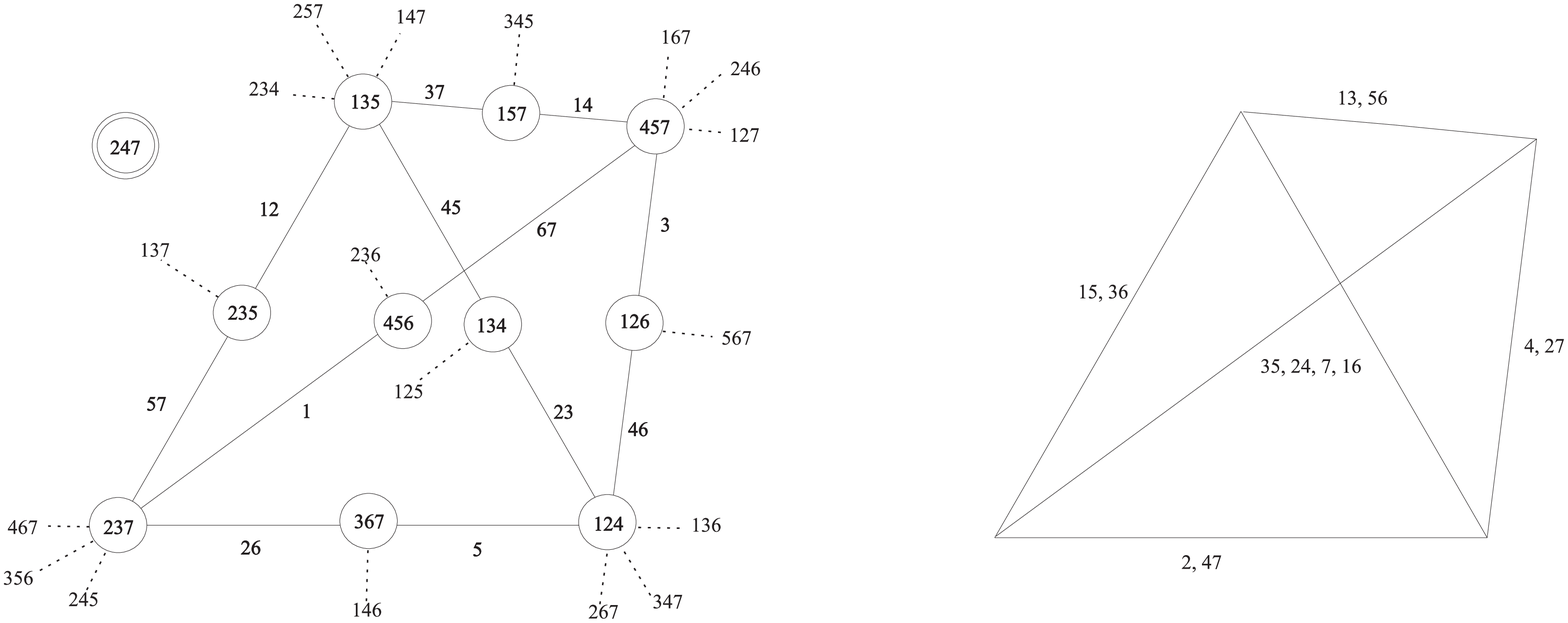}
  \caption{\small}\label{supertetrafig}
\end{figure}

In each case we list the $D_4$ and the $4$-tuples
and show that they have an intrinsic characterization, i.e., they are preserved by the stabilizer $(S_4)_\Theta$.
Thus it will suffice to prove \eqref{formulaforU} and to show
that $|F_1\cap\Gamma|=|F_2\cap\Gamma|$ for
a representative $\Gamma$ of each $(S_4)_\Theta$-orbit.
We leave this count (similar to the case of $E_6$) to the reader.

$\Theta=247$ is an $A_1$ of Lemma~\ref{extraA1}.
We take the vertices of the tetrahedron and the unique $D_4=D(17,25,34)$ that contains them with fourtuples
$$
\{247,123,357,145\},\quad\{6,17,25,34\}.
$$

$\Theta=237$ is a vertex. We take $D_4 = D(34,26,57)$
(its Dynkin diagram is formed by $\Theta$ and midpoints of adjacent edges)   with fourtuples
$$
\{237,245,467,356\},\quad \{1,34,26,57\}
$$

$\Theta=367$ is a midpoint. Take the two $A_1$ endpoints $\Gamma_1$ and $\Gamma_2$ of the edge (in our case $237$
and $124$). Take the unique $A_7\not\ni\Theta$ and consider its four extra roots.
Exactly two of them ($2$ and $56$ in our case)
are perpendicular to $\Gamma_1$ and $\Gamma_2$.
This gives $4$ roots that belong to a unique $D_4=D(14,37)$. We take the fourtuples
$$
\{6,25,146,367\},\quad\{5,26,145,357\}.
$$

$\Theta=A_2$, e.g.~given by a half-edge $26$.
$A_2$ is contained in three $A_3$'s and their three extra (i.e.~not contained in $A_2$'s)
roots ($146$, $245$, and $356$ in our case) along with the root on the half-edge $26$
form a Dynkin diagram of $D_4=D(26,34,15)$ with fourtuples
$$\{7,15,26,34\},\quad \{123,146,245,356\}.$$

$\Theta = A_3^{\times2}$ given by two vertices, e.g.
$237$ and $457$.
We consider the following $6$ roots: extra roots of $A_3$'s in  ($356$ and $246$),
extra roots of $A_7\supset\Theta$ that are not perpendicular to the two vertices
($7$ and $16$), the midpoint of an edge connecting the two vertices ($456$), and
the extra root of an $A_3$ formed by the two vertices and the midpoint above ($236$).
Then there is only one $D_4$ containing all $6$ roots, namely $D(16,25,34)$, and
we take the fourtuples
$$
\{7,16,34,25\},\quad \{123,145,356,246\}.
$$

$\Theta= A_3^{\times2}$ given by a pair of disjoint edges, e.g.,
$135-157-457$ and $237-367-124$.
There are two $A_7$'s that contain $\Theta$, and each of them has two
extra roots perpendicular to midpoints of the disjoint edges
($7$, $24$, $15$, $36$ in our case).
There is only one $D_4$ that contains them, namely $D(24, 15, 36)$, and we take the fourtuples
$$
\{123,345,146,256\},\quad \{7,24,15,36\}.
$$

$\Theta = A_7$, e.g., not containing $134$ and $456$.
Consider $4$ extra roots of $A_7$ (namely, $7$, $16$, $24$, $35$)
and the unique $D_4= D(16,24,35)$ that contains them, with fourtuples
$$
\{7,16,24,35\}, \quad \{123,145,346,256\}.
$$
\end{proof}

\section{$\oY(\Delta)$ is Smooth Log Canonical, and Tropical Compactification}\label{logcantrick}

\begin{Theorem} \label{lcm} Let $(\oY,B)$ be a pair of a
(possibly reducible) proper variety with reduced boundary, which has stable toric singularities.
Then $K_{\oY} + B$ is ample iff each irreducible open stratum is log minimal
(in the reducible case strata are defined using boundary and double locus divisors).
\end{Theorem}

\begin{proof} Since for an irreducible component $\oS \subset \oY$ we
have
$(K_{\oY} + B)|_{\oS} = K_{\oS} + B_{\oS}$, we may assume $\oY$ is
irreducible.
By adjunction, for any closed stratum $\oS$ with interior $S$, 
$(K_{\oY} + B)|_{\oS} = K_{\oS} + B_{\oS}$, where $B_{\oS} = \oS \setminus S$. 
Thus if $K_{\oY} + B$ is ample, $S$ is log minimal and $S \subset \oS$ is its log canonical
compactification.

In the other direction we induct on $\dim\oY$. To show that $K_{\oY} + B$ is ample,
we apply the Nakai--Moishezon criterion, i.e., we show that the restriction to every subvariety is big.
By induction (and adjunction) it suffices to consider a subvariety
which
meets the interior  $Y := \oY \setminus B$.
But $m$-pluri log canonical forms of $Y$ extend to any
log canonical compactification and so for some $ m > 0$, the rational
map
$$
\oY \dasharrow \bP(H^0(\oY,\cO(m(K_{\oY}+B)))^*)
$$
restricts to an immersion on $Y$, since by assumption $Y$ is log
minimal. Thus the  restriction of the map to any
subvariety meeting $Y$ is birational.
\end{proof}

\begin{Corollary}\label{n6lcm}
$\oY(E_6)$ and $\oY(D_n)$ are the log canonical compactifications.
\end{Corollary}

\begin{proof} The compactifications have simple normal crossing
boundary, so by Th.~\ref{lcm} we only need to check that open
strata are all log minimal. Let $X(r,n)$ be the moduli space of ordered $n$-tuples
of points in $\bP^{r-1}$ in linear general position. All strata are
open subsets of products of various $X(r,n)$ (indeed of $X(3,6)$
or $X(2,i)$). This is familiar for $\oM_{0,n}$. For $\oY(E_6)$ see
\cite{Naruki}. An open subvariety of
a log minimal variety is log minimal, and $X(r,n)$ is
log minimal by \cite[2.18]{ChQII}. \end{proof}

Recall that, for $\ch k \neq 2$, $\oY(E_7)$ denotes the compactification of $Y(E_7)$ given by $\tM^{(2)}_3$,
the coarse moduli space of the blowup of the moduli stack of stable curves of genus $3$ with level $2$ structure
at the $T$ locus.

\begin{Review}\textsc{Fano Simplex.}\label{Fano}
Let $\bP^2(\bF_2)=\{1\dots 7\}$. Each of the $7$ lines $\bP^1(\bF_2)\subset\bP^2(\bF_2)$
is given by a subset $\{ijk\}$ and
gives rise to the root $\alpha_{ijk}\in E_7$. These roots
are orthogonal because two lines intersect at a point.
Fig.~\ref{fanoplane} shows one possible choice:
\begin{figure}[htbp]
  \includegraphics[width=1.4in]{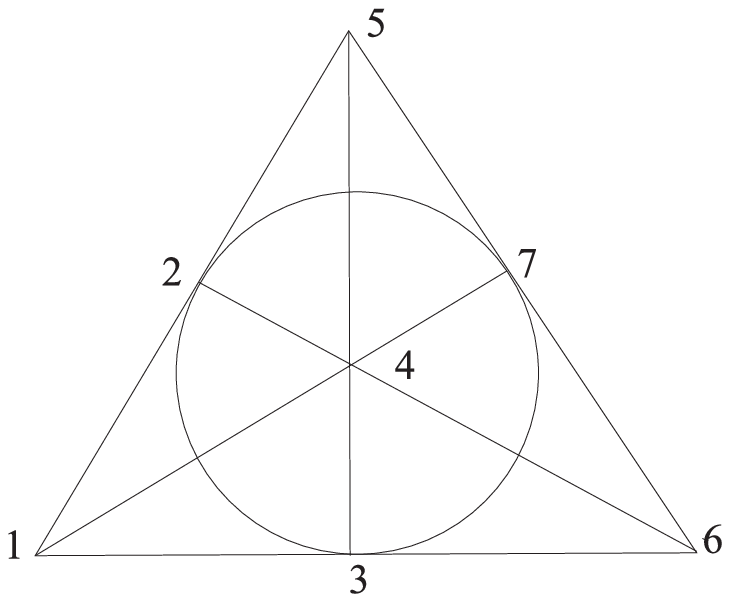}\\
  \caption{\small Fano plane}\label{fanoplane}\label{fanosimplex}
\end{figure}
\end{Review}

\begin{Lemma} \label{A1locus} Let $A_1(\oY(E_7))\subset \oY(E_7)$ be the
complement to the union of boundary divisors
not of type $A_1$. This open set, with its boundary, is smooth
with simple normal crossings. It is covered by open subsets
$(V,B)$ of the following form.
Consider the configuration of $6$ lines in $\bP^2$
as in Fig.~\ref{fanosimplex} without the (impossible in $\hbox{\rm
char}\ne2$) line $\{237\}$.
Let $\tilde V \subset (\bP^2)^7$ be the
open subset of $7$-tuples of points with no degenerations other than
those in the configuration (so any triple of points other than the $6$
triples in Fig.~\ref{fanosimplex} are not collinear, and no $6$ points lie on a conic).
$\PGL_3$ acts freely on $\tilde V$, let $V$ be the quotient.
Let $B \subset V$ be the union of the $6$ boundary divisors.
\end{Lemma}

\begin{proof}
Over $V$ there is a family of del Pezzo
surfaces with at worse ordinary $A_1$ singularities (the
anti-canonical model of the blowup of $\bP^2$ at the given $7$ points).
The anti-canonical map for each fibre is a double cover of $\bP^2$ branched over
a quartic with at worst ordinary nodes. This determines
a regular map $V \to \oM_3 \setminus (\oH \cup \delta_1)$.

\begin{Claim}
The map is surjective and quasi-finite.
\end{Claim}

\begin{proof}
Choose $[C] \in \ocM_3 \setminus (\ocH \cup \delta_1)$. $C$ embeds
in $\bP^2$ as a quartic nodal curve. Let $S' \to \bP^2$ be the
double cover branched over $C$, and $\tS \to S'$ the minimal
desingularisation. The nodes of $C$ give rise to a collection
of disjoint $(-2)$-curves on~$\tS$. By \eqref{rootsEn}, it corresponds to a
subsystem $A_1^{\times r} \subset E_7$, $r \leq 6$. We can
extend this to an $A_1^{\times 6}$. Any two such are conjugate
under $W(E_7)$ by Th.~\ref{classification}, which implies we can find
a birational morphism $\tS \to \bP^2$ realizing $\tS$ as the
blowup at $7$ distinct points, with only collinear degeneracies,
as in the statement of the theorem.
\end{proof}

Let $Z \subset \tMM \times V$ be
the closure of the graph of the rational birational map $V \to \tMM$. Then
$Z \to V$ is proper and birational, and quasi-finite.
Thus as $V$ is normal, it is an isomorphism. So
$V \to \tMM$ is a birational quasi-finite morphism, and so
an open immersion by Zariski's main theorem.
Now we take the orbit of $V \subset \tMM$
under $W(E_7)$. Each of these open sets admits an
analogous description, and by the Claim they cover
$p^{-1}(\oM_3 \setminus \oH\cup\delta_1)$, which
is precisely $A_1(\oY(E_7))$. \end{proof}

\begin{Lemma} \label{qs} Let $(\cM,\cD_1 + \cD_2 \dots + \cD_r)$
be a smooth Deligne-Mumford stack with a normal crossing divisor.
Assume
the restriction of the inertia stack to each of the substacks
$\cD_i$ contains a finite \'etale sub group-scheme, and
that these sub-group schemes together generate the automorphisms
of $\cM$. Then the coarse moduli space
$(M,D_1 + \dots D_r)$ is smooth with normal crossings.
\end{Lemma}

\begin{proof} The question is \'etale local, and
it follows from the construction of the coarse moduli space in \cite{KeelMori}
that for any geometric point $\Spec k\to \cM$, with image $x\in M$,
there exists an ´etale neighborhood $U$ of $x$ in $M$ such that $\cM\times_M U$
is isomorphic to a stack of the form $[X/G]$ for some scheme $X$,
where $G$ is the geometric stabilizer group of $x$.
Using our assumptions, we can assume further that
$X$ is a vector space with a linear action of $G$,
with a collection of $G$-invariant
hyperplanes $\cE_i$ ($[\cE_i/G] = \cD_i$), given by linearly
independent
linear functionals ($X$ is the Zariski tangent space to $\cM$ at a point).
The conditions imply that $G$ is
generated by
subgroups $G_i$ acting trivially on $\cE_i$. Clearly $G_i$ is cyclic,
generated
by a reflection in $\cE_i$. We can mod out by the intersection $\cap
\cE_i$
(on which $G$ acts trivially) to reduce to the case when the linear
functionals
give a basis. Now we can choose coordinates so that $G_i$ is generated
by
a diagonal matrix with all ones, except for a root of unity in the
$i$th position.
One checks the claim easily in this case. \end{proof}

\begin{Lemma} \label{smoothcodim2} $\oY(E_7)$ has normal
crossings generically along any codimension two boundary stratum
contained in a divisor of type $D(A_7)$.
\end{Lemma}

\begin{proof}
We show the corresponding statement for $(\tM_3,B)$ holds.
The result then follows from the proof Th.~\ref{stacksmooth} --
indeed the loops used in it are (homotopic to loops that are)
contained
in the automorphism free locus of $M_3$, and the proof of
Th.~\ref{stacksmooth}
depends only on the monodromy action of these loops.
Now we check the criteria of Lemma~\ref{qs}.
The possible strata correspond to boundary strata of
$\oH=\oM_{0,8}/S_8$, thus $B_2, B_3$ and $B_4$ as in \cite{CornalbaHarris}.

The generic point of $B_2$ corresponds to a nodal curve such that its normalization
is a smooth genus $2$ curve and a pair of points over the node
is a conjugate pair for the (unique) hyperelliptic involution.
It is clear that this nodal curve has automorphism group
$\bZ/2\bZ$, the (restriction of the) hyperelliptic involution.

The generic point of $B_3$ corresponds to a $1$-pointed elliptic curve glued
to a $1$-pointed hyperelliptic curve (of genus~$2$) at a Weierstrass
point. In this case the automorphism group is $(\bZ/2\bZ)^2$,
generated by the restrictions of the involutions at
the generic points of $\oH$ and $\Delta_1$.

The generic point of $B_4=T \subset \oM_3$ is the union of a pair
of two 2-pointed genus~$1$ curves. The automorphism group is
(generically)
$\bZ/2\bZ$ -- the involution on each component interchanging the marked
points, the restriction of the hyperelliptic involution.
The same then holds at $T = E \cap \tilde{H} \subset \tM_3$.
\end{proof}

\begin{Lemma} \label{surjectA7} The restriction map
$r:\,M^{D(A_7)}(E_7) \to \cO^*(M_{0,8})/k^* = M(D_8)$
is surjective,
where $M^{D(A_7)}$ means units with valuation~$0$ along $D(A_7)$ as in~\eqref{M^Snotation}.
\end{Lemma}

\begin{proof}
First we show that we can find a unit $u\in M^{D(A_7)}(E_7)$ which has valuation~$1$
along some boundary divisor meeting $D(A_7)$,
i.e., the image of $r$ is indivisible. Choose perpendicular
subsystems $A_1,A_1'$
of the given $A_7 \subset E_7$.
It is easy to see using the explicit chart of Lemma~\ref{A1locus}
that we can find a unit with valuation~$1$ on $A_1$, and~$0$
on $A_1'$. If it has zero valuation on $D(A_7)$ we
are done. Otherwise we can choose an element of the Weyl
group which preserves the $A_7$ and interchanges $A_1$ and $A_1'$.
Now take the difference between the original
relation, and the one obtained from it by applying the
permutation.

The map $r$ is $S_8$-equivariant and by Lemma~\ref{torsionfree} the image is irreducible.
$M(E_7)_\bQ$ contains a single copy of $M(D_8)_\bQ$
(e.g., by dimension count), and therefore it remains to show that
there exists a surjective map $r':\,M^{D(A_7)}(E_7) \to M(D_8)$ of $\bZ[S_8]$-modules:
surjectivity of $r$ then follows from Schur's lemma.
We choose the standard $A_7\subset E_7$ \eqref{standA7}.
Let $f:\,\bZ^{\cA_1(E_7)}\to\bZ^{\cA_1(A_7)}\to\bZ^{\cA_1(D_8)}$
be the composite map, where the first map is the natural restriction
and the second map sends $[\eps_i-\eps_j]$ to $[\eps_i-\eps_j]+[\eps_i+\eps_j]$.
Here $i,j\in\{1,\ldots,8\}$ and we identify the root $\beta_{i8}\in A_7\subset E_7$
with the root $\eps_i-\eps_8$ of the standard $A_7\subset D_8$ of \eqref{anrem}.
Let $r'$ be the restriction of $f$ to $M^{D(A_7)}(E_7)\subset\bZ^{\cA_1(E_7)}$.
We claim that $r'$ is surjective onto $M(D_8)\subset \bZ^{\cA_1(D_8)}$.
By Th.~\ref{messy}, any $m\in M^{D(A_7)}(E_7)$ is a linear combination of $D_4$-units,
i.e., a linear combination of differences $F_1-F_2$, where $F_1$ and $F_2$
are fourtuples of a $D_4$ of one of the three types of Remark~\ref{D4class}.
For $D_4$'s of the first and third type, $F_1-F_2\in M^{D(A_7)}(E_7)$
(i.e., the coefficients at the roots in $A_7$ add up to $0$) and $r'(F_1-F_2)$
is a $D_4$-unit of $D_8$.
For a $D_4$ of the second type, $F_1-F_2$ (restricted to $A_7$)
is either equal to $0$ or is a linear combination of $4$ orthogonal roots of $A_7$.
It is easy to see that if $m\in M^{D(A_7)}(E_7)$ is
a linear combination of such units then its restriction to $A_7$
can be rewritten using restriction to $A_7$ of $D_4$-units
of the first and second type. Since $M(D_8)$ is also generated by $D_4$ units,
we are done.
\end{proof}

\begin{Lemma}\label{paulobserv}
$D(A_3\times A_3)\subset \oY(E_7)$
is isomorphic to $\oM_{0,5} \times \oM_{0,5} \times \oM_{0,4}$,
and the pullback of one boundary point from the third factor is
$D(A_1) \cap D(A_3 \times A_3)$ where $A_1 = (A_3 \times A_3)^{\perp}$,
while the pullback of the other
two boundary points are $D(A_7) \cap D(A_3 \times A_3)$ for
the two $A_7$'s which contain $A_3 \times A_3$.
\end{Lemma}

\begin{proof}
The divisor $D(A_3 \times A_3)$ is a component of $E^{(2)} \subset \tM_3^{(2)}$,
so this follows from Lemma~\ref{E}.
\end{proof}

\begin{Theorem} \label{pzeros}\label{surjBergman}
There is a commutative diagram
$$
\begin{CD}
M(\Delta) @>{\oplus b_{\Gamma}}>> \bZ^{\cR(\Delta)} \\
@V{i}VV  @V{j}VV \\
\cO^*(Y(\Delta))/k^* @>{\val}>> \bZ^{\partial \oY(\Delta)}
\end{CD}
$$
where $\val$ sends a unit to the sum of its valuations at each
of the boundary divisors,
$i$ is the identification of Th.~\ref{pullbackofunits},
$j$ is the identification of Lemma~\ref{RDn} and Prop.~\ref{classification},
and, for each type $\Gamma$ of root system in $\cR(\Delta)$, $b_{\Gamma}:=c_\Gamma/m_\Gamma$, where
$c_\Gamma$ is the composition
$$
c_{\Gamma}:\, M(\Delta) \arrow^{\psi^\vee} \bZ^{\Delta_{+}} \arrow^A
\bZ^{\Gamma(\Delta)},\quad A([\alpha])=\sum_{\alpha\in\Theta\in\Gamma(\Delta)}[\Theta]
$$
and $m_{\Gamma}$ is the largest integer such that $c_\Gamma(M(\Delta))\subset m_\Gamma \bZ^{\Gamma(\Delta)}$.

We define a collection of cones $\cF(\partial\oY(\Delta))\subset N(\Delta)_{\bQ}$
from the combinatorics of the boundary of $\oY(\Delta)$ as in Cor.~\ref{constcasencc}.
The map $i^\vee$ identifies this collection with~$\cF(\Delta)$.
In particular, $\cF(\Delta)$ is supported on the tropicalization $\cA(Y(\Delta))$.
\end{Theorem}

\begin{proof}
The map dual to $c_\Gamma$ has the form $c_\Gamma^\vee:\, [\Theta]\mapsto\psi(\Theta)$.
It is clear from Th.~\ref{dualbasis} that $c_{\Gamma}\ne0$ (and is therefore injective)
and that $m_\Gamma$ is given by Th.~\ref{PsiCalculation}.

Let $\Gamma \in \cR$ be a root subsystem, let
$D(\Gamma) \in \partial \oY(\Delta)$ be the corresponding boundary divisor,
and let $[D(\Gamma)]\in N_{Y(\Delta)}$ be the corresponding point.
Let $\xi(\Gamma)$ be the first lattice point of $N(\Delta)$ along
the ray spanned by $\psi(\Gamma)$. Note commutativity of
the diagram is equivalent to $i^\vee([D(\Gamma)])=\xi(\Gamma)$.
Assuming this, the equality of the collections of cones
$\cF(\Delta)$ and $\cF(\partial\oY(\Delta))$ follows from Lemma~\ref{RDn} and Prop.~\ref{classification}.
Let $H\subset W(\Delta)$ be the normalizer of $\Gamma$ (and thus of $D(\Gamma)$).
Note that $\psi(\Gamma)$ and $[D(\Gamma)]$ are $H$-invariant.

\begin{Claim}
$\dim N(\Delta)^H=1$.
\end{Claim}

\begin{proof}
Let $\Delta=D_n$.
We use the equivariant surjection $\psi:\,\bZ^{\cA_1}\to N$ of Lemma~\ref{exsequn}.
$H=W(D_k)\times W(D_{n-k})$ has three orbits on $\cA_1$, namely
$\cA_1(D_k)$, $\cA_1(D_{n-k})$, and their complement.
It follows that $N(\Delta)^H$ is generated by $\psi(D_k)$, $\psi(D_{n-k})$, and
$\psi(D_n)-\psi(D_k)-\psi(D_{n-k})$. The first two vectors are equal
(see the proof of Th.~\ref{PsiCalculation}), and $\psi(D_n)=0$
(being a $W$-invariant vector in an irreducible module).

If $\Delta=E_6$ or $E_7$ then we can either make a similar calculation or note that,
by Frobenius reciprocity, $\dim N(\Delta)^H$ is equal
to the multiplicity of $N(\Delta)$ in the permutation module
$\bZ^{\Gamma(\Delta)}$. The latter is equal to $1$ by the tables of~\cite{Atlas}.
\end{proof}

It follows that vectors $\psi(\Gamma)$ and $[D(\Gamma)]$ are collinear.

\begin{Claim}
The vectors $\psi(\Gamma)$ and $[D(\Gamma)]$ are not opposite.
\end{Claim}

\begin{proof}
We know that $[D(\Gamma)]\in \cA_{Y(\Delta)}$ by \eqref{constcasencc}.
We claim that $\psi(\Gamma)$ also belongs to the tropicalization.
By Th.~\ref{pullbackofunits}, the dominant map $\cX(\Delta)\arrow^\Psi Y(\Delta)$ induces
the linear map $\bZ^{\Delta_+}\arrow^\psi N(\Delta)$
that restricts to the map of tropicalizations $\cA(\cX)\to\cA(Y)$.
The fan $\tilde\cF(\Delta)$ supported on $\cA(\cX)$ is described in~\cite{Williams}
but we only need the description of rays: they correspond to 
irreducible root subsystems $\tilde\Theta\subset\Delta$ such that $\Gamma(\tilde\Theta)\subset\Gamma(\Delta)$.
The corresponding ray of $\tilde\cF(\Delta)$ is given by 
$\sum_{\alpha\in\tilde\Theta^+}[\alpha]$.
Therefore, it suffices to find, for each $\Theta\in\cR(\Delta)$, a subsystem $\tilde\Theta\subset\Delta$
realizable on the Dynkin diagram and such that $\psi(\Theta)$ and $\psi(\tilde\Theta)$
are proportional. This is easy using the table of Theorem~\ref{PsiCalculation}.

Now we argue by contradiction. If $\psi(\Gamma)$ and $[D(\Gamma)]$ are opposite
then $\cA(Y(\Delta))$ contains a line $L$. Take any $D_4\subset\Delta$.
Then $\cA(Y(\Delta))$ maps onto $\cA(Y(D_4))$ by Theorem~\ref{D4inclpullback}.
Since $\cA(Y(D_4))=\cA(M_{0,4})$ does not contain lines
(being the union of three rays), $L$ maps to $0$ in $N(D_4)$.
But by Theorem~\ref{messy}, the map $N(\Delta)\to\bigoplus\limits_{D_4\subset\Delta}N(D_4)$ is injective.
This is a contradiction.
\end{proof}

To complete the proof it's enough to
show that $[D(\Gamma)]$ is the first lattice point along
its ray. For this its enough to find some unit with
valuation one along this divisor.
For $\oY(D_n) = \oM_{0,n}$ this is easily checked --- one can
take a classical cross-ratio of some $4$-tuple, see~\cite{Keel}.
For $\oY(E_6) = \oY^6$ (where there are only two orbits of
divisors) we can take any of Naruki's cross-ratio maps, see~\cite{Naruki}.

Now assume $\Delta = E_7$. Suppose first that $\Gamma\ne A_7$.
We choose an $A_7$ such that $D(A_7)$ intersects $D(\Gamma)$.
Note if we have a unit $u$
with valuation zero along $Y(A_7)$, then, by Lemma~\ref{smoothcodim2},
its valuation along $D(\Gamma)$ is
the same as the valuation of $u|_{D(A_7)}$ along $D(\Gamma) \cap D(A_7)$.
As above, on $D(A_7) = \oM_{0,8}$ we can find a unit with valuation one along
any given boundary divisor. So the result follows from Lemma~\ref{surjectA7}.

To finish the proof it suffices to find a unit with valuation one
along some $Y(A_7)$. Suppose it does not exist.
Then, by equivariance, valuation of any unit along any $D(A_7)$ is divisible by $d>1$.
Consider subsystems $A_3 \times A_3$ and
$A_1 := (A_3 \times A_3)^{\perp}$.
By Th.~\ref{intfan} (and the previous case $\Gamma\ne A_7$), we can find a unit $u$ that
has valuation $0$ along $D(A_3 \times A_3)$ and valuation $1$ along~$D(A_1)$.
By Lemma~\ref{paulobserv}, the valuations
of $u$ at the two $A_7$ boundary divisors that meet $D(A_3\times A_3)$ add up to $-1$, and therefore
must be coprime.
This contradicts our assumption.
\end{proof}

\begin{Corollary} \label{snc} $\oY(E_7)=(\tM^{(2)},B)$ is smooth with
simple normal crossings.
\end{Corollary}

\begin{proof}
By the strict simpliciality of the cones $\sigma\in\cF(\Delta)$ (Th.~\ref{intfan}),
each boundary divisor is Cartier: this is a local question, and by the strict
simpliciality, given a collection of boundary divisors
cutting out a stratum $S$, we can find a unit with valuation
one on one of them, and zero on the others --- and so a
boundary divisor is locally principal.
Our calculations have shown each boundary divisor is smooth.
Thus $\oY(E_7)$ is smooth.
Since each boundary divisor has simple normal crossing boundary
induced (by the strict simpliciality) by the boundary of $\oY(E_7)$,
we have simple normal crossings for $\oY(E_7)$.
\end{proof}

\begin{Theorem} \label{intfan}
Let $\cP=\prod \limits_{D_4 \subset \Delta} \cF(D_4)\subset
N'=\bigoplus \limits_{D_4 \subset \Delta} N(D_4)$
be the product fan.
$\cF(\Delta)$ is a convexly disjoint and strictly simplicial fan,
and each of its cones is a cone in the intersection
fan $N(\Delta)_{\bQ} \cap \cP$, where
we identify $N(\Delta)$ with a sublattice in $N'$ using Th.~\ref{messy}.
The map of toric varieties $X(\cF) \to X(\cP)$
is an immersion, i.e., an isomorphism of $X(\cF)$
with an open subscheme of a closed subscheme of~$X(\cP)$.
\end{Theorem}

\begin{proof}
Let $\sigma$ be a cone of $\cF(\Delta)$.
By Cor.~\ref{D4inclpullback}, each map
$N(\Delta) \to N(D_4)$ is induced by the
morphism of very affine varieties $Y(\Delta) \to Y(D_4)$.
So by Th.~\ref{dominanta} we
have a surjective map of tropical sets
$|\cF(\Delta)| \to |\cF(D_4)|$.
Since $\cF(D_4)$ is one-dimensional, the image of any face of $\sigma$ is a cone of $\cF(D_4)$.
By Lemma~\ref{subfan}, $\sigma$ is a cone of $\cG(\Delta)$.
It follows (by Th.~\ref{dualbasis})
that $\Ann\langle\sigma\rangle\subset M(\Delta)_\bQ$
is generated by the restriction of elements from the various $M(D_4)$,
where $\langle\sigma\rangle\subset N(\Delta)_{\bQ}$ means the linear
span and $\Ann\langle\sigma\rangle\subset M(\Delta)_\bQ$
is the annihilator $\Ker[M(\Delta)_\bQ\to\langle\sigma\rangle^\vee]$.
Now it follows from Lemma~\ref{conecond} below that $\sigma$ is a cone of the intersection fan.

It follows that $\cF(\Delta)$ is a fan. It is strictly simplicial by Th.~\ref{dualbasis} and Lemma~\ref{subfan}.
To show that it is convexly disjoint we note that: (1) $\cF(D_4)$ is convexly disjoint;
(2)~the product of convexly disjoint collections is convexly disjoint;
(3) the intersection of a convexly disjoint collection with a linear subspace is convexly disjoint.

Let $U$ be the $T_N$-toric variety given by the intersection fan $N(\Delta)_{\bQ} \cap \cP$.
It contains $X(\cF)$ as a toric open subset.
Since $N$ is saturated in $N'$ by Th.~\ref{messy}, $U$ is the normalisation of its image under the 
the natural map $\nu:\,U\to X(\cP)$, 
and $\nu$ induces a bijection of the sets of $T_N$ orbits of $U$ and $\nu(U)$ \cite[1.4]{Oda}.
Therefore it suffices to prove that $\nu(U(\sigma))$ is normal
for any affine open toric subset $U(\sigma)$ given by a cone $\sigma\in\cF(\Delta)$.
Let $\sigma'=\prod_i \pi_i(\sigma)\subset N'_\bQ$, where $\pi_i:\,N\to N(D_4)$ are the projections.
We claim that $U(\sigma)\to U(\sigma')$ is a closed embedding.
It suffices to prove that the map of semigroups $(\sigma')^\vee\cap M'\to\sigma^\vee\cap M$ is surjective.
But this follows from Th.~\ref{dualbasis}.
\end{proof}

\begin{Remark} Th. \ref{intfan} holds if we
replace $\cF$ by $\cG$ everywhere (in particular, $\cG$ is itself
a fan). We don't need this here.
The same
proof holds if one can establish that under $N(\Delta) \to N(D_4)$
each cone of $\cG(\Delta)$ maps onto a cone of $\cG(D_4)$. Since
this is true of the one dimensional cones (which are part of $\cF$)
it's enough to show this is a map of fans. This follows
from Sekiguchi's descriptions of the real cells (which appear
without proof in the $E_7$ case) -- as it's
clear (by continuity) that $Y(\Delta)(\bR) \to Y(D_4)(\bR) =
M_{0,4}(\bR)$
sends a cell into a cell.
\end{Remark}

\begin{Lemma} \label{conecond} Let $N \subset \bigoplus_i N_i$
be a sub-lattice of a finite direct sum of lattices. Let
$\sigma \subset N_{\bQ}$ be a rational polyhedral cone. Assume
that for each face $\gamma \subset \sigma$ and each projection
$\pi_i: N \to N_i$, $\pi_i(\gamma)$ is a face of $\pi_i(\sigma)$.
Then $\sigma =N_{\bQ}\cap\prod_i \pi_i(\sigma)$
if for each face $\gamma \subset \sigma$
the natural map $\bigoplus_{i} \Ann\langle\pi_i(\gamma)\rangle\to\Ann\langle\gamma\rangle$
is surjective.
\end{Lemma}

\begin{proof} For any face $\gamma \subset \sigma$ let
$P(\gamma)= \prod_i \pi_i(\gamma)$ and $I(\gamma)= P(\gamma) \cap N_\bQ$.
By assumption, $P(\gamma)$ is a face of $P(\sigma)$ and
so $I(\gamma)$ is a face of $I(\sigma)$.
$\gamma \subset I(\gamma)$ and,
if we can show that they have the same dimension,
then every face of $\sigma$ lies in a face of $I(\sigma)$ of the same dimension
and therefore $\sigma=I(\sigma)$.
It suffices to check that $\langle\gamma\rangle=\langle I(\gamma)\rangle$,
i.e., $\langle\gamma\rangle=N_\bQ\cap\bigoplus_i \langle\pi_i(\gamma)\rangle$.
But this clearly follows from the surjectivity of
$\bigoplus_{i} \Ann\langle\pi_i(\gamma)\rangle\to\Ann\langle\gamma\rangle$.
\end{proof}

\begin{Proposition}\label{boundarycorrespondence}
Let $\Gamma \in \cR(\Delta)$. Th.~\ref{pzeros} induces
the identification of boundary divisors of $D(\Gamma)$ with $\link_{\Gamma}(\cR(\Delta))$.
It has the following explicit form:

\underline{$A_1 \subset E_6$.} $\link_{\Gamma}(\Delta)$ is given by $A_1$ and $A_2^{\times2}$
subsystems of $\Gamma^{\perp}=A_5$. In the notation of Lemma~\ref{anrem},
these correspond to
$(2,4)$ and $(3,3)$ bipartitions of $N_6$ respectively, which
correspond to boundary divisors of $D(\Gamma) =\oM_{0,N_6}$.

\underline{$A_2^{\times3}\subset E_6$.}
The link
is given by $A_1$ subsystems of $\Gamma$. Identifying each factor
with
a single standard $A_2$, and fixing $N_3 \subset N_4$, each $A_1$
gives a two element subset of one of three copies of $N_3 \subset N_4$,
and thus a boundary
divisor of $D(\Gamma) = \oM_{0,N_4}^{\times 3}$.

\underline{$A_7 \subset E_7$.} The link consists of $A_1$,
$A_2$, and $A_3^{\times2}\subset A_7$. These correspond to $(2,6)$, $(3,5)$ or $(4,4)$
bipartitions of $N_8$, and thus boundary divisors of $D(\Gamma)=\oM_{0,8}$.

\underline{$A_2 \subset E_7$.} Identify
$\Gamma$ and $\Gamma^{\perp}=A_5$ with the standard $A_2$ and $A_5$.
Fix $N_3 \subset N_4$ and $N_6 \subset N_7$. The link
consists of $A_1 \subset \Gamma$
(boundary divisors of $\oM_{0,4}$) and of
$A_1 \subset \Gamma^{\perp}$, $A_2 \subset \Gamma^{\perp}$, $A_3^{\times2}\supset\Gamma$, or $A_7\supset\Gamma$.
The last two possibilities are
equivalent to $A_3:=A_3^{\times2}\cap\Gamma^\perp \subset \Gamma^{\perp}$
and $A_4:=A_7\cap\Gamma^\perp \subset\Gamma^{\perp}$, respectively.
These subsystems correspond to $k$-element subsets of $N_6$, $2\le k\le 5$, and thus to boundary divisors of
the second factor of $D(\Gamma) = \oM_{0,4} \times \oM_{0,7}$.

\underline{$A_3^{\times2}\subset E_7$.} The link consist of $A_1$ or
$A_2$ subsystems of the two factors, $A_7\supset \Gamma$ or $A_1 =
\Gamma^{\perp}$.
The last two possibilities correspond to the boundary divisors from
the third factor of $D(\Gamma) = \oM_{0,5} \times \oM_{0,5} \times
\oM_{0,4}$
and the other possibilities correspond (as above) to boundary divisors
of the first two factors.
\end{Proposition}

\begin{proof}
It is easy to check using Table~\ref{perpsE7} that the given map is a $W(\Delta)_{\Gamma}$-equivariant bijection.
Now we claim that $\Aut_{W_\Gamma}(\link_\Gamma\cR)$ is trivial
except for $A_3^{\times2}\subset E_7$, in which case
it is equal to $\bZ/2\bZ$ induced by the involution of $D(\Gamma)$
given by the reflection in $A_1:=(A_3^{\times2})^\perp$.
Indeed, let $G:=W(\Gamma\times\Gamma^\perp)$ Then $G\subset W_\Gamma$ and therefore
if $f\in\Aut_{W_\Gamma}(\link_\Gamma\cR)$ then $G_\Theta=G_{f\Theta}$ for any $\Theta$ in the link.
Since $\Gamma\times\Gamma^\perp$ has at most one direct summand~$A_1$,
a~subsystem $\Theta\subset\Gamma\times\Gamma^\perp$ of type $A_1$ is uniquely determined by $G_\Theta$.
Since any $A_1\subset\Gamma\times\Gamma^\perp$ belongs to the link,
$f$ must preserve all of them.
It follows that $f$ preserves all subsystems in the link that belong to $\Gamma\times\Gamma^\perp$
because they are uniquely determined by incident $A_1$'s in the~link.
A simple calculation shows that in fact $f$ preserves all subsystems in the~link, except
when $\Gamma=A_3^{\times2}$: in this case it is possible to permute the two $A_7$'s that contain~$\Gamma$.
This permutation is induced by the reflection in $A_1:=(A_3^{\times2})^\perp$.
\end{proof}

\begin{Theorem}\label{monhubsch} Let $\Delta=D_n$
or $E_n$, $n \leq 7$. $Y(\Delta)$ is H\"ubsch,
$\oY(\Delta)$ is the log canonical
compactification, and $\cF(\Delta)$ is the log canonical
fan. The rational map
\begin{equation}\label{sekiembedding}
\oY(\Delta) \to P:=\prod \limits_{D_4 \subset \Delta} \oY(D_4)
\cooltag\end{equation}
is a closed embedding. Each open (resp. closed) boundary stratum
is the scheme-theoretic inverse image of the open (resp. closed)
stratum of $P$ which contains~it.
\end{Theorem}

\begin{Remark}
Th.~\ref{monhubsch} and Cor.~\ref{snc} imply Sekiguchi's Conjecture~\ref{sekconj}.
\end{Remark}

\begin{proof}
All statements about the embedding of $\oY(E_7)$ into
the product of $\bP^1$'s follow from Th.~\ref{intfan}
once we prove that $\cF(\Delta)$ is the log canonical fan.

$Y(\Delta) \subset \oY(\Delta)$ is a simple
normal crossing compactification by Cor.~\ref{snc}.
The open strata are $Y(\Delta)$, $D^0(A_1)\subset\oY(E_7)$,
or products of $M_{0,i}$'s. They are log minimal
as open subsets of the complements to connected hyperplane arrangements
(for $D^0(A_1)$ this follows from Lemma~\ref{A1locus}).
It follows that $\oY(\Delta)$ is the log canonical compactification (Th.~\ref{lcm}).
We prove $Y(\Delta)$ is H\"ubsch by showing
for each stratum the map
$S \to T^S_Y$ of Th.~\ref{hubsch} is an immersion. For $Y(\Delta)$
this is clear, as $Y(\Delta)$ is very affine. For purely
$A_1$ strata of $\oY(E_7)$
this is easily checked using the
explicit charts of Lemma~\ref{A1locus}. All other strata
are products of $M_{0,i}$, which are very affine. So
it's enough by Th.~\ref{hubsch} to show
that
$M_{Y(\Delta)}^S \to M_S$ is surjective for such strata.

Proof for $D_n$: by Cor.~\ref{D4inclpullback} and
Th.~\ref{messy}, $M_S$ is generated
by the pullback of units from the canonical
cross-ratio maps $S \to M_{0,4}$.
Thus it is enough to show that each such
map is the restriction of a map
$M_{0,n}\to M_{0,4}$. This is well known.

Proof for $E_6$ and $E_7$: all boundary divisors (other than $D(A_1)\subset\oY(E_7)$) are products of $\oM_{0,i}$, and
so, by the $D_n$ case, to prove
$M_{Y(\Delta)}^S \to M_S$ is surjective
it's
enough to consider the case where $S$ is codimension
one. Let $r:=\rk M_S$ (given in~Lemma~\ref{torsionfree}).
It's enough to find
$r$ boundary divisors $D_1,\dots, D_r$, all incident
to $\oS$, and units $u_1,\dots,u_r$ so that
$\val_{D_i}(u_j) = \delta_{ij}$, and $\val_S(u_i) = 0$,
for then $u_i \in M_Y^S$, and their images in $M_S$
will give a basis of the lattice. We choose a maximal cone
$\sigma \in \cG(\Delta)$ which has a ray corresponding
to the boundary divisor $S$.
By Th.~\ref{dualbasis}, it suffices to prove the following claim:

\begin{Claim} The cone $\sigma$
contains $r$ rays corresponding to
divisors (other than $S$) which
have non-empty intersection with $\oS$.
\end{Claim}

\begin{Remark}
This follows  from
Th.~\ref{dualbasis} (applied to the boundary divisor)
if we know the rays in a cone indeed correspond to the faces of a cell, as
the intersection of a cell with a boundary divisor will be a cell for
the corresponding compactification.
\end{Remark}

We check the claim for each $W(\Delta)$ orbit of boundary
divisor, using Fig.~\ref{pentatetra}.
Proof for $E_6$: an $A_2^{\times3}$ subdiagram of the pentadiagram
contains $6$ $A_1$'s which yield $u_1,\dots,u_6$ for the boundary divisor
$\oM_{0,4}^{\times 3} = D(A_2^{\times3}) \subset \oY(E_6)$.
An $A_1$ subdiagram is contained in $3$ $A_2^{\times3}$ subdiagrams and is perpendicular
to $6$ $A_1$'s. This gives $9$ units for the boundary divisor
$\oM_{0,6} = D(A_1) \subset \oY(E_6)$.
Proof for $E_7$: an $A_7$ of the tetra-diagram
contains $8$ $A_1$'s, $8$ $A_2$'s, and $4$ $A_3^{\times2}$'s,
which yield $u_1,\dots,u_{20}$ for the boundary divisor
$\oM_{0,8} = D(A_7) \subset \oY(E_7)$.
Consider an $A_3^{\times2}$ subdiagram that corresponds to the pair of opposite edges.
It contains $6$ $A_1$'s and $4$ $A_2$'s, and is contained in $2$ $A_7$'s.
Another type of an $A_3^{\times2}$ subdiagram contains $6$ $A_1$'s and $4$ $A_2$'s, is contained in an $A_7$,
and is perpendicular to an $A_1$ from the tetradiagram.
In both cases this gives $12$ units of $M_{0,5}\times M_{0,5}\times M_{0,4}$.
Finally, consider an $A_2$ subdiagram. It contains $2$ $A_1$'s,
is perpendicular to $5$ $A_1$'s and $4$ $A_2$'s, and is contained in $3$ $A_3^{\times2}$'s
and $2$ $A_7$'s.
This gives $16$ units of $M_{0,4}\times M_{0,7}$.
\end{proof}

\section{Moduli of Stable Del Pezzo Surfaces}\label{sds}

\begin{Proposition} \label{pimapoffans}
Let $\Delta' \subset \Delta$ be either $D_n \subset D_{n+1}$, or
$E_n \subset E_{n+1}$ for $n \leq 6$.
The~map $\pi:\,N(\Delta) \to N(\Delta')$ of Lemma~\ref{subs} induces a map of fans $\cF(\Delta) \to \cF(\Delta')$.
Let $p:\,Y(\Delta) \to Y(\Delta')$ be the following map:
For $D_n$, $M_{0,n+1}\to M_{0,n}$ is the forgetful map.
For $E_n$, $Y^{n+1}\to Y^{n}$ is \eqref{mapP}.
Then $\pi^\vee$ is the pullback of units $p^*:\,\cO(Y(\Delta'))/k^*\to\cO(Y(\Delta))/k^*$.
We have a commutative diagram of morphisms
$$
\begin{CD}
\cX(\Delta) @>{\Psi}>> Y(\Delta) @>>> X(\cF(\Delta)) \\
@V{\pr}VV   @V{p}VV @V{\pi}VV \\
\cX(\Delta') @>{\Psi'}>> Y(\Delta') @>>> X(\cF(\Delta'))
\end{CD}
$$
where $\pr$ is the projection dual to the embedding of root lattices $\Lambda(\Delta')\hookrightarrow \Lambda(\Delta)$.
\end{Proposition}

\begin{proof}
Under the identifications of Th.~\ref{pullbackofunits} (for $\Delta$ and $\Delta'$),
the pullback of units $p^*$ is identified with $\pi^\vee$ by commutativity of~\eqref{basicdiagram}.
By Th.~\ref{dominanta}, $p$ induces a surjective map of tropicalizations
$\cA_{Y(\Delta)}\to\cA_{Y(\Delta')}$. This map is automatically a map of fans
$\cF(\Delta) \to \cF(\Delta')$ by Th.~\ref{minfan}.
It induces the map of toric varieties and therefore (by~Th.~\ref{monhubsch})
the regular map $\oY(\Delta) \to \oY(\Delta')$.
\end{proof}

Next we describe $\pi:\,\cF(\Delta) \to \cF(\Delta')$ explicitly.
For any $[\Gamma]\in\cR$, we denote by $\zeta(\Gamma)$ the first lattice
point along the corresponding ray of $\cF$.

\begin{Proposition} \label{piforDn}\label{piforE6}
\underline{$D_n\subset D_{n+1}$.}
Rays of $\cF(\Delta)$ correspond to bipartitions of $N_{n+1}$,
\begin{equation}\label{pondivofm0n}
\pi(\zeta(I\sqcup I^c)) =
\begin{cases}
\zeta((I\cap N_n)\sqcup(I^c\cap N_n)) & \hbox{\rm if}\ |I\cap N_n|,\ |I^c\cap N_n|>1\cr
0                   & \hbox{\rm otherwise}\cr
\end{cases}\cooltag
\end{equation}

\underline{$E_5:=D_5 \subset E_6$.}
For $A_1 \in \cR(E_6)$,
\begin{equation}\label{pondivofnar}
\pi(\zeta(A_1)) =
\begin{cases}
\zeta(D_2\subset D_5) & \hbox{\rm if}\ A_1\subset D_5\cr
0                   & \hbox{\rm otherwise (there are $16$ of those)}\cr
\end{cases}\cooltag
\end{equation}
We use that, in $D_5$,  $A_1\times A_1^\perp=(A_1\times A_1')\times D_3=D_2\times D_3$.
For $A_2^{\times 3}\in\cR(E_6)$,
\begin{equation}\label{pondivofnar2}
\pi(\zeta(A_2^{\times3})) = \zeta(D_2\subset D_5),
\ \hbox{\rm where}\ A_2^{\times 3}\cap D_5 = A_2 \times A_1^{\times 2}=A_2\times D_2.
\cooltag
\end{equation}
In both cases the toric map $X(\cF(\Delta)) \to X(\cF(\Delta'))$
is flat, with reduced fibres.
\end{Proposition}

\begin{proof}
Rays of $\cF(D_n)$ correspond to bipartitions by~Lemma~\ref{anrem}.
The first lattice points are found in Th.~\ref{PsiCalculation}
(see also the proof of Th.~\ref{dualbasis}).
The identities (\ref{pondivofm0n}--\ref{pondivofnar2}) follow from Lemma~\ref{subs} and Th.~\ref{PsiCalculation}.
The identity $A_2^{\times3} \cap D_5 = A_2 \times A_1^{\times2}$ is easy to check.
The last statement follows from Lemma~\ref{toricflat}.
\end{proof}

\begin{Lemma}\label{toricflat}
Let $N \to N'$ be a surjection of
lattices giving a map of fans $\cF \to \cF'$.
Assume $\cF'$ is strictly simplicial and any ray of $\cF$ maps to a ray of $\cF'$.
Then the toric map $X(\cF)\arrow^\pi X(\cF')$ is flat.
If, in addition, the first lattice point of any ray of
$\cF$ maps to the first lattice point of a ray
in $\cF'$ (or to $0$) then $\pi$ has reduced fibres.
\end{Lemma}

\begin{proof}
Since $\cF'$ is strictly simplicial, every cone in $\cF$ maps onto a cone
of $\cF'$, which implies that $\pi$ is equidimensional.
Since $X(\cF')$ is regular and $X(\cF)$ is Cohen--Macaulay \cite{Da},
this implies that $\pi$ is flat \cite[6.1.5]{EGA4}.
The condition on lattice points implies that the scheme-theoretic inverse image
of any toric stratum in $X(\cF')$ is reduced.
This implies that all fibers are reduced by equivariance.
\end{proof}

\begin{Definition}
Let $\tcF(E_6)$ be the refinement of $\cF(E_6)$
obtained by performing the barycentric subdivision
of cones of type $\psi\{A_1,\ldots, A_1\}$
and taking the corresponding minimal subdivision of
cones of type $\psi\{A_2^{\times 3}, A_1, \ldots, A_1\}$
(meaning we barycentrically subdivide the codimension one face spanned
by the $A_1$ rays, and then add the final $A_2^{\times3}$
ray to each cone in this subdivision).
\end{Definition}

\begin{Proposition}\label{piE7toE6} For $E_6 \subset E_7$, $\pi$
has the following effect on the rays of $\cF(E_7)$:
\begin{equation}\label{pieffectseki1}
\pi(\zeta(A_1)) =
\begin{cases}
\zeta(A_1\subset E_6) & \hbox{\rm if}\ A_1\subset E_6\cr
0                   & \hbox{\rm otherwise (there are $27$ of those);}\cr
\end{cases}
\cooltag\end{equation}
\begin{equation}\label{pieffectseki2}
\pi(\zeta(A_2)) =
\begin{cases}
\zeta(A_1\subset E_6) & \hbox{\rm if}\ A_2\cap E_6=A_1\cr
\zeta(A_2^{\times3}) & \hbox{\rm if}\ A_2\subset E_6,
\ \hbox{\rm where}\ A_2^{\times3}=A_2\times A_2^\perp;\cr
\end{cases}
\cooltag\end{equation}
\begin{equation}\label{pieffectseki3}
\pi(\zeta(A_3^{\times2})) =
\begin{cases}
\zeta(A_2^{\times3})   & \hbox{\rm if}\ A_3^{\times2}\cap E_6=A_2^{\times2},\
\hbox{\rm here}\ A_2^{\times3}:=A_2^{\times2}\times(A_2^{\times2})^\perp\cr
\zeta(A_1)+\zeta(A_1') & \hbox{\rm if}\ A_3^{\times2}\cap E_6=A_3\times A_1\times A_1';\cr
\end{cases}
\cooltag\end{equation}
\begin{equation}\label{pieffectseki4}
\pi(\zeta(A_7))=\zeta(A_1\subset E_6),
\ \hbox{\rm where}\ A_7\cap E_6=A_1\times A_5.
\cooltag\end{equation}

The map $X(\cF(E_7)) \to X(\cF(E_6))$ has reduced fibres.

The fan $\tcF(E_6)$ is strictly simplicial. It is the unique
minimal refinement of $\cF(E_6)$ so that $\pi(\sigma)$ is a union
of cones for each $\sigma \in \cF(E_7)$.
The collection of cones
$$\tcF(E_7)=\{\pi^{-1}(\gamma) \cap \sigma\ |\ \gamma \in \tcF(E_6),\ \sigma \in \cF(E_7)\}$$
is a fan.
The map $X(\tcF(E_7)) \to X(\tcF(E_6))$ is flat, with reduced fibres.
\end{Proposition}

\begin{proof}
The first lattice points are found in Th.~\ref{PsiCalculation}
(see also the proof of Th.~\ref{dualbasis}). Formulas for the image of $\pi$
follow from Lemma~\ref{subs} and Th.~\ref{PsiCalculation}
up to (easily checked) formulas for the intersection of root subsystems of $E_7$
with $E_6$.

The only case when the image of a ray in $\cF(E_7)$
is contained in the relative interior of a cone in $\cF(E_6)$
is $\psi\{A_3^{\times2}\}$ mapped onto the bisector
of~$\psi\{A_1,A_1'\}$.
Fig.~\ref{fignonflata3} shows a cone $\sigma \in \cF(E_7)$ such that
$\pi(\sigma)$ is a $4$-dimensional cone with rays
$\psi\{A_1'\}, \psi\{A_1''\}, \psi\{A_1'''\}$, and the bisector of $\psi\{A_1,A_1'\}$.
It easily follows that $\tcF(E_6)$ is the minimal
refinement such that $\pi(\sigma)$ is a union
of cones for each $\sigma \in \cF(E_7)$.
\begin{figure}[htbp]
  \includegraphics[width=4in]{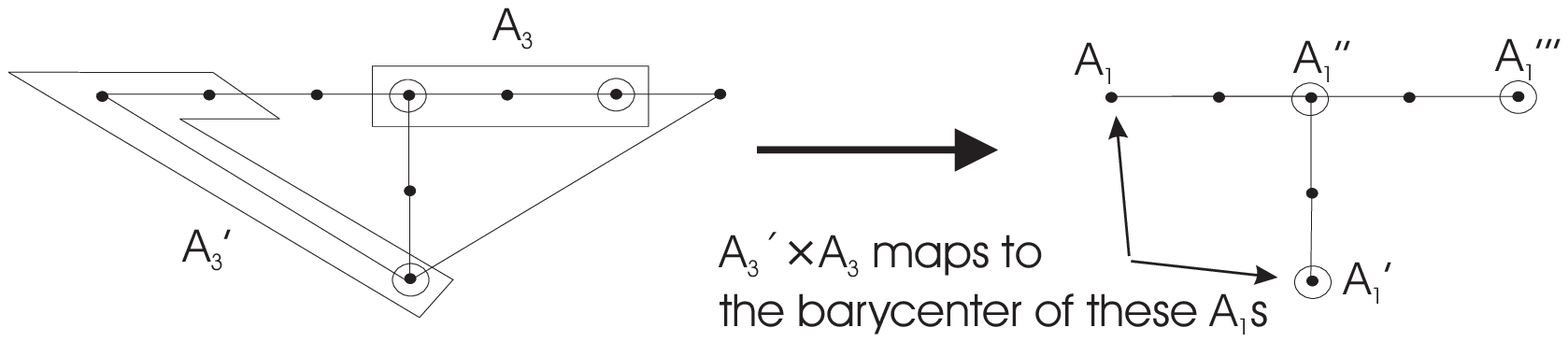}\\
  \caption{\small}\label{fignonflata3}
\end{figure}

It is clear $\tcF(E_7)$ is a fan, and that
a cone of $\tcF(E_7)$ (as a subset of the lattice!)
maps onto a cone of $\tcF(E_6)$.
As the latter is strictly simplicial, the map of toric varieties is
flat, with
reduced fibres, by Lemma \ref{toricflat}.
\end{proof}

\begin{Definition}\label{rigid} A very
affine variety $Y$ is called {\em rigid in its intrinsic torus}~$T_Y$ if
$$
\{t \in T_Y\,|\, t \cdot Y \subset Y\} = \{e\}.
$$
\end{Definition}

\begin{Lemma} \label{hcomprigid}
The complement to a connected hyperplane arrangement is rigid in its intrinsic torus.
Any open boundary stratum of $\oY(E_6)$ is rigid in its intrinsic torus.\end{Lemma}

\begin{proof} For the first statement, see \cite[\S4]{Tevelev}.
The boundary strata of $\oY(E_6)$ are all
products of $M_{0,k}$ and thus complements to connected hyperplane arrangements.
\end{proof}

\begin{Lemma}\label{brigid} Let $\oZ \arrow^f \oY$ be
a surjective morphism of normal proper varieties with divisorial
boundary
such that $f^{-1}(\partial \oY) \subset \partial \oZ$.
The pullback of boundary divisors $\bZ^{\partial \oY} \arrow^{f^*} \bZ^{\partial \oZ}$
determines $f$ if $Y := \oY \setminus \partial \oY$ is
very affine and  rigid in
its intrinsic torus.
\end{Lemma}

\begin{proof}
The map $\cO^*(Y)/k^* \to \cO^*(Z)/k^*$ is determined by the pullback of boundary divisors
(by pulling back the divisor of zeroes and poles of $u\in\cO^*(Y)$).
This determines the composition $Z \to Y \to T_Y$
up to a translation in $T_Y$ preserving~$Y$. \end{proof}

\begin{Proposition} \label{restrictions}
In the notation of Prop.~\ref{piforE6} and Prop.~\ref{piE7toE6},
we describe the restriction $q:=p|_{D(\Gamma)}$ of $\oY(\Delta)\arrow^p\oY(\Delta')$ to
$D(\Gamma)\subset\oY(\Delta)$ for $\Gamma\in\cR(\Delta)$.
Below all maps $\oM_{0,i} \to \oM_{0,j}$, $j\le i$ are canonical fibrations given by dropping points.
$$\hbox{\rm for \eqref{pondivofnar}:}\quad
q=
\begin{cases}
\oM_{0,6} \to \oM_{0,4}=D(D_2\subset D_5)\cr
\oM_{0,6} \to \oM_{0,5}=\oY(D_5);\cr
\end{cases}$$
$$\hbox{\rm for \eqref{pondivofnar2}:}\quad
q:\ \oM_{0,4} \times \oM_{0,4} \times\oM_{0,4} \arrow^{\pr_1}  \oM_{0,4} = D(D_2\subset D_5);$$
$$\hbox{\rm for \eqref{pieffectseki2}:}\quad
q=
\begin{cases}
\oM_{0,4} \times \oM_{0,7} \to \{pt\}\times\oM_{0,6} = D(A_1)\cr
\oM_{0,4} \times \oM_{0,7} \to \oM_{0,4} \times [\oM_{0,4}
\times \oM_{0,4}]= D(A_2),\cr
\end{cases}$$
where in the second case the map is induced (with appropriate
ordering) from
\begin{equation}\label{funnyM07}
\begin{CD}
\oM_{0,7} @>{\pi_{1,2,3,4} \times \pi_{4,5,6,7}}>> \oM_{0,4} \times
\oM_{0,4};\end{CD}\cooltag
\end{equation}
$$\hbox{\rm for \eqref{pieffectseki3}:}\quad
q=
\begin{cases}
\oM_{0,5} \times \oM_{0,5} \times \oM_{0,4} \to \oM_{0,4} \times
\oM_{0,4} \times \oM_{0,4}=D(A_2^{\times3})\cr
\oM_{0,5} \times \oM_{0,5} \times \oM_{0,4} \arrow^{\pr_1} \oM_{0,5}=D(A_1)\cap D(A_1').\cr
\end{cases}$$
$$\hbox{\rm for \eqref{pieffectseki4}:}\quad
q:\ \oM_{0,8} \to \oM_{0,6} = D(A_1 \subset E_6).$$
\end{Proposition}

\begin{Definition}
We refer to the second possibility for $A_2 \subset E_7$
(resp.~the second possibility for $A_3^{\times2}\subset E_7$)
as a {\em non-flat $D(A_2)$} (resp.~a {\em non-flat $D(A_3 \times A_3)$}).
\end{Definition}

\begin{proof}
By Th.~\ref{monhubsch} the map of fans determines the image of each boundary
stratum, which in turn determines the pullback map
on Weil divisors for the map
of a stratum to its image (since by Th.~\ref{monhubsch}, Prop.~\ref{piforE6}, and Prop.~\ref{piE7toE6}
the pullback of a boundary divisor for this map will be reduced).
By Lemma~\ref{brigid} it is enough to check that this pull-back agrees with the pull-back
for the maps in the theorem (which we denote by $t$).
We use Prop.~\ref{boundarycorrespondence}, \ref{piforE6}, and~\ref{piE7toE6}
and suggest that the reader draws Dynkin diagrams.
We use the standard $D_5=E_5\subset E_6\subset E_7$ (as in Rk.~\ref{rootsEn}) and
$A_2^{\times3}\subset E_6$, $A_3^{\times2},A_7\subset E_7$ (as in Fig.~\ref{conesofF(Delta)}) unless otherwise noted.
Each paragraph has its own temporary notations for root subsystems.

\underline{\eqref{pondivofnar}, $A_1\subset D_5$.}
We choose $A_1=7$ in $D_5$ with simple roots $7,123,23,34,45$.
We identify $D(A_1)=M_{0,\{1,\ldots,6\}}$ and $t=\pi_{2,3,4,5}$.
Up to conjugation, any boundary divisor $D$ of $p(D(A_1))$ is its intersection
with $D(D_2'=A_1''\times A_1'''\subset D_5)$, where $A_1''=23$ and $A_1'''=45$.
The only components of $p^{-1}(D(D_2'))$ intersecting $D(A_1)$ are
$D(A_1'')$, $D(A_1''')$, $D(A_2^{\times 3})$, and $D(A_2'^{\times 3})$,
where $A_2=\{123,456,7\}$, and $A_2'=\{145,236,7\}$.
Thus
$q^{-1}(D)=\delta_{23}\cup\delta_{45}\cup\delta_{123}\cup\delta_{145}=t^{-1}(\delta_{23})$.

\underline{\eqref{pondivofnar}, $A_1\not\subset D_5$.}
We choose $A_1=7$. We have $t=\pi_{1,2,3,4,5}$. Up to conjugation,
a boundary divisor $D$ of $\oY(D_5)$ is $D(D_2=A_1'\times A_1'')\subset\oY(D_5)$,
where  $A_1'=45$, and $A_1''=123$.
The only components of $p^{-1}(D)$ intersecting $D(A_1)$ are
$D(A_1')$ and $D(A_2^{\times3})$, where $A_2=\{12,13,23\}$.
Thus $q^{-1}(D(D_2))=\delta_{123}\cup\delta_{45}=t^{-1}(\delta_{45})$.

\underline{\eqref{pondivofnar2}.}
Up to conjugation any boundary divisor $D$ of $p(D(A_2^{\times3}))$ is its intersection
with $D(D_2=A_1\times A_1'\subset D_5)$, where $A_1=12$ and $A_1'=345$.
The only component of $p^{-1}(D(D_2))$ intersecting $D(A_2^{\times 3})$ is $D(A_1)$.
It follows that $q^{-1}(D)=\delta_{12}\times \oM_{0,4}\times\oM_{0,4}=t^{-1}(\delta_{12})$.

\underline{\eqref{pieffectseki2}, $A_2\cap E_6=A_1$.} Choose $A_2=\{56,57,67\}$
and identify the $\oM_{0,7}$ component of $D(A_2)$ with $\oM_{0,\{0,\ldots,4,x,y\}}$,
$t=\pi_{1,\ldots,4,x,y}$.
$p(D(A_2))=D(A_1\subset E_6)=\oM_{0,\{1,2,3,4,x,y\}}$ has $3$ conjugacy classes of boundary divisors:

(a) the intersection with $D(A_1'\subset E_6)$, where $A_1'=7$.
The only components of $p^{-1}(D(A_1'))$ that intersect $D(A_2)$ are $D(A_7)$
and $D(A_3^{\times2})$, where $A_3$ has simple roots $12,23,34$.
It follows that $q^{-1}(D)=\oM_{0,4}\times[\delta_{01234}\cup\delta_{1234}]=t^{-1}(\delta_{xy})$.

(b) the intersection with $D(A_1'\subset E_6)$, where $A_1'=123$.
The only two components of $p^{-1}(D(A_1'))$ that intersect $D(A_2\subset E_7)$ are $D(A_1')$
and $D(A_2)$, where $A_2=\{123,4,567\}$.
It follows that $q^{-1}(D)=\oM_{0,4}\times[\delta_{4x}\cup\delta_{04x}]=t^{-1}(\delta_{4x})$.

(c) the intersection with $D(A_2'^{\times3}\subset E_6)$, where $A_2'=\{123,7,456\}$.
The only two components of $p^{-1}(D(A_2'^{\times3}))$ that intersect $D(A_2\subset E_7)$
are $D(A_2'')$ and $D(A_3^{\times2})$, where $A_2''=\{12,13,23\}$.
It follows that $q^{-1}(D)=\oM_{0,4}\times[\delta_{0123}\cup\delta_{123}]=t^{-1}(\delta_{456})$.

\underline{\eqref{pieffectseki2}, $A_2\subset E_6$.}
Choose $A_2=\{123,456,7\}$. Then $A_5=A_2^\perp$ has simple roots $12$, $23$, $347$, $45$, $56$.
We have $D(A_2)=\oM_{0,\{a,b,c,d\}}\times\oM_{0,\{1,2,3,x,5,6,y\}}$,
$t=id\times[\pi_{123y}\times\pi_{x56y}]$.
Up to symmetries, $p(D(A_2))=\oM_{0,4}^{\times3}$ has $2$ types of boundary divisors:

(a) the intersection with $D(A_1)$ for $A_1=7$.
The only component of $p^{-1}(D(A_1))$ that intersects $D(A_2)$ is
$D(A_1)$. Thus $q^{-1}(D)=\delta_{c,d}\times\oM_{0,7}=t^{-1}(\delta_{c,d}\times\oM_{0,4}\times\oM_{0,4})$.

(b) the intersection with $D(A_1)$ for $A_1=12$.
Components of $p^{-1}(D(A_1))$ that intersect $D(A_2)$ are
$D(A_1)$, $D(A_2^{(i)})$, $i=4,5,6$,
$D({}^{(i)}\!A_3^{\times2})$, $i=4,5,6$, and~$D(A_7)$,
where $A_2^{(i)}=\{12,1i7,2i7\}$, ${}^{(i)}\!A_3$ has simple roots $7,123,3i$,
and $A_7$ has simple roots $12$, $1$, $7$, $123$, $34$, $45$, $56$.
The corresponding decomposition of $q^{-1}(D)$ is
$$\oM_{0,4}\times[\delta_{12}\cup\delta_{12x}\cup\delta_{125}\cup\delta_{126}
\cup\delta_{1256}\cup\delta_{12x6}\cup\delta_{12x5}\cup\delta_{12x56}]=t^{-1}(\oM_{0,4}\times\delta_{12}\times\oM_{0,4}).
$$

\underline{\eqref{pieffectseki3}, $A_3^{\times2}\cap E_6=A_2^{\times2}$}.
$D(A_3^{\times2})=\oM_{0,\{0,1,2,3,x\}}\times\oM_{0,\{4,5,6,7,y\}}\times\oM_{0,\{a,b,c,d\}}$.
Then $t=\pi_{1,2,3,x}\times\pi_{4,5,6,y}\times id$.
Up to symmetries, $p(D(A_3^{\times2}))=\oM_{0,4}^{\times3}$ has the following types of boundary divisors:

(a) the intersection with $D(A_1)$ for $A_1=7$.
The only component of $p^{-1}(D(A_1))$ that intersects $D(A_3^{\times2})$
is $D(A_7)$.
Thus $q^{-1}(D)=\oM_{0,5}^{\times2}\times\delta_{c,d}=t^{-1}(\oM_{0,4}^{\times2}\times\delta_{c,d})$.

(b) the intersection with $D(A_1)$ for $A_1=123$.
The only component of $p^{-1}(D(A_1))$ that intersects $D(A_3^{\times2})$
is $D(A_1)$. Thus $q^{-1}(D)=\oM_{0,5}^{\times2}\times\delta_{b,d}=t^{-1}(\oM_{0,4}^{\times2}\times\delta_{b,d})$.

(c) the intersection with $D(A_1)$ for $A_1=12$.
Components of $p^{-1}(D(A_1))$ that intersect $D(A_3^{\times2})$ are
$D(A_1)$ and $D(A_2)$ with $A_2=\{1,12,2\}$.
It follows that $q^{-1}(D)=[\delta_{12}\cup\delta_{012}]\times\oM_{0,5}\times\oM_{0,4}=t^{-1}(\delta_{12}\times\oM_{0,5}\times\oM_{0,4})$.

\underline{\eqref{pieffectseki3}, $A_3^{\times2}\cap E_6=A_3\times A_1\times A_1'$}.
We choose $A_3^{\times2}$ containing $A_3$ with simple roots $12,23,34$.
The $\oM_{0,5}$ component of $D(A_3^{\times2})$ that corresponds to $A_3$
is identified with $\oM_{0,\{1,2,3,4,x\}}$ and its boundary divisors are identified with
$A_1\subset A_3$ ($ij\mapsto\delta_{ij}$) and $A_2\subset A_3$ ($\{ij,jk,ik\}\mapsto\delta_{ijk}$).
Since the link of $\sigma=\{A_1,A_1'\}\subset\cR(E_6)$ is also equivariantly identified
with $A_1,A_2\subset A_3$ (with $A_2$ corresponding to $A_2\times A_2^\perp$),
the boundary divisors of $\oM_{0,5}=D(A_1)\cap D(A_1)\subset\oY(E_6)$
correspond to $A_1,A_2\subset A_3$ in the same way as above
(by the argument used in Prop.~\ref{boundarycorrespondence}).
It remains to note that the only component of $p^{-1}(D(A_1))$ (resp.~$p^{-1}(D(A_2^{\times3}))$)
for $A_1\subset A_3$ (resp.~$A_2\subset A_3$)
that intersects $D(A_3^{\times2})$ is $D(A_1)$ (resp.~$D(A_2)$).
Therefore, $q^{-1}(D)=t^{-1}(D)$.

\underline{\eqref{pieffectseki4}.}
There are $2$ types of boundary divisors $D$ of $p(D(A_7))=\oM_{0,6}$.

(a) the intersection with $D(A_1)$ for $A_1=12$.
The components of $p^{-1}(D(A_1))$ that intersect $D(A_7)$ are $D(A_1)$, $D(A_2)$ with $A_2=\{1,12,2\}$,
$D(A_2')$ with $A_2'=\{12,17,27\}$, and $D(A_3^{\times2})$, where $A_3$ has simple roots
$1,12,27$. It follows that
$q^{-1}(D)=\delta_{12}\cup\delta_{012}\cup\delta_{127}\cup\delta_{0127}=t^{-1}(\delta_{12})$.

(b) the intersection with $D(A_2^{\times3})$ for $A_2=\{12,23,13\}$.
The components of $p^{-1}(D(A_2^{\times3}))$ that intersect $D(A_7)$ are $D(A_2)$, $D(A_3^{\times2})$,
where $A_3$ has simple roots $1,12,23$,
$D(A_3'^{\times2})$, where $A_3'$ has simple roots $12,23,37$,
and $D(A_2')$ with $A_2'=\{45,46,56\}$.
It follows that
$q^{-1}(D)=\delta_{123}\cup\delta_{0123}\cup\delta_{1237}\cup\delta_{456}=t^{-1}(\delta_{123})$.
\end{proof}

\begin{Lemma}\label{normalbundleZ}
Suppose $A_1, A_1'\subset E_6$ are orthogonal.
The codimension $2$ stratum $Z:=D(A_1)\cap D(A_1')\subset\oY(E_6)$ has trivial projectivized normal bundle.
\end{Lemma}

\begin{proof}
It suffices to prove that normal bundles to $Z$ in $D(A_1)$ and $D(A_1')$ are isomorphic.
It is well-known that, for the embedding $\oM_{0,n}=\delta_{n,n+1}\subset \oM_{0,n+1}$,
$\cO(\delta_{n,n+1})|_{\delta_{n,n+1}}=\psi_n^*$.
Therefore it suffices to check that $Z$ embeds both in $D(A_1)$ and $D(A_1')$ as a section $\delta_{xy}$,
where we identify both $D(A_1)$ and $D(A_1')$ with $\oM_{0,\{1,2,3,4,x,y\}}$
and $Z$ with $\oM_{0,\{1,2,3,4,x\}}$ as in the previous proof
(the case \eqref{pieffectseki3}, $A_3^{\times2}\cap E_6=A_3\times A_1\times A_1'$).
But this follows from the identification of the link of $Z$ in $\cR(E_6)$ with
boundary divisors of $\oM_{0,5}$ because other $\delta_{ij}$ divisors of $D(A_1)$ and $D(A_1')$
either properly intersect $Z$ or are disjoint from it.
\end{proof}

\begin{Lemma}\label{flatcriterion} Let $\pi:\,Y \to Y'$ be
a dominant morphism from an integral scheme to a normal
scheme, with reduced fibres of constant dimension. Then $\pi$ is flat.
\end{Lemma}

\begin{proof} Use \cite[15.2.3]{EGA4} and \cite[14.4.4]{EGA4}.
\end{proof}

\begin{Review}
For basic properties of log structures we refer to
\cite{Kato} and \cite{Ol}.  Any log structure
we use in this paper will be toric, i.e., the space will
come with an evident map to a toric variety, and we endow
the space with the pullback of the toric log structure on
the toric variety. In fact, we do not make any use of
the log structure itself, only the bundles of log (and
relative log) differentials as in \cite[2.19]{ChQII}.
\end{Review}

If $X$ is a toric variety with torus $T$ and $Y \subset X$ is a closed subvariety, we refer to the multiplication map
$T \times Y \rightarrow X$, $(t,y) \mapsto t \cdot y$
as the \emph{structure map}.

\begin{Proposition}\label{logsmooth}
Let $X \arrow^\pi X'$ be a dominant toric map of toric varieties with reduced fibres of constant dimension.
Let $Y \subset X$, $Y' \subset X'$ be closed subvarieties with smooth
structure maps. Let $W \subset Y$, $W' \subset Y'$ be irreducible strata.
Suppose $\pi$ induces dominant maps $Y \arrow^p Y'$ and $W\arrow^q W'$.
If\/ $q$ is log smooth at $z \in W$ then $p$ is log smooth at~$z$.
\end{Proposition}

\begin{proof}
Suppose first that $\codim_Y W=\codim_{Y'}W'=1$.
We have a commutative diagram of vector bundles on $W$
$$
\begin{CD}
0 @>>> \Omega^1_W(\log) @>>> \Omega^1_Y(\log)|_W @>{r}>> \cO_W @>>> 0 \\
@.    @A{q^*}AA          @A{p^*}AA             @|  @. \\
0 @>>> q^*(\Omega^1_{W'}(\log)) @>>> q^*(\Omega^1_{Y'}(\log)|_{W'}) @>{r}>>
q^*(\cO_{W'}) @>>> 0
\end{CD}
$$
Here $\Omega^1_Y(\log)=\Omega^1_Y(\log B_Y)$ and $r$ is the residue map.
The first vertical arrow is a subbundle embedding by assumption.
Thus the second arrow is a subbundle embedding, i.e., $Y \rightarrow Y'$ is log smooth at
$z \in W$.

Suppose now that $\codim_Y W=1$ and $W'=Y'$. Then the map
$p^*\Omega^1_{Y'}(\log B_{Y'}) \to \Omega^1_Y(\log B_Y)|_W$ factors through
$\Omega^1_W(\log B_W)$ and so is a subbundle embedding.

The general case follows by induction:
If $W'=Y'$ find a stratum $W_0\supset W$ such that $\codim_{Y}W_0=1$.
Then $W_0 \rightarrow Y'$ is log smooth, and we are done by induction on
$\codim_YW$.
If $W' \neq Y'$, choose a component $W_0$ of $p^{-1}(W')$ containing~$W$.
Then $W_0 \rightarrow W'$ is log smooth.
Choose a stratum $W'_0\supset W'$ such that $\codim_{W'_0}W'=1$,
and let $W_1$ be a component of $p^{-1}(W'_0)$ such that $W_0 \subset W_1$
and $W_1 \rightarrow W'_0$ is dominant.
Then $\codim_{W_1}W_0=1$,
so $W_1 \rightarrow W'_0$ is log smooth and we are done by induction on $\codim_{Y'}W'$.
\end{proof}

\begin{Proposition} \label{ls} Let $X \arrow^\pi X'$ be a flat map of toric varieties
with reduced fibers extending a surjective homomorphism $T \to T'$ of algebraic tori, with kernel~$T''$.
Assume $X'$ is smooth.
Let $Y \subset X$, $Y' \subset X'$ be subvarieties with smooth
structure maps and assume $\pi|_Y$ is a dominant map $Y \arrow^p Y'$,
log smooth at $y \in Y$.
Then the induced map $f:\,Y \times T'' \to F:=Y' \times_{X'} X$ is smooth at $y$.
\end{Proposition}

\begin{proof}
Define $Q$ (locally near $y$) by the exact sequence
$$
0 \to \Omega^1_{Y' \times T'/X'}|_{Y \times T''} \to  \Omega^1_{Y
\times T/X}|_{Y \times T''}
\to  Q \to 0
$$

The first term is the pullback of the log cotangent bundle
of $Y'$
and the second of the log cotangent bundle of $Y$ \cite[2.19]{ChQII}. Thus $Q$ is a vector bundle
of rank
$$
r = \dim Y - \dim Y' = \dim (Y \times T'') - \dim F.
$$
A simple computation shows that $Q = \Omega^1_{Y \times T''/F}$.
Therefore $f$ is equidimensional and has smooth fibers (locally near $y$).
It suffices to prove that $F$ is normal: then $f$ is flat by Lemma \ref{flatcriterion},
and so smooth \cite[6.8.6]{EGA4} (locally near $y$).

Observe $F \to Y'$ is flat, has reduced Cohen--Macaulay fibres, and normal generic fibre.
Since $Y'$ is smooth, $F$ is Cohen--Macaulay and has $R_1$.
Thus $F$ is normal.
\end{proof}

We recall the definition of the Koll\'ar--Shepherd-Barron--Alexeev moduli stack $\ocM$ of stable surfaces with boundary.
Assume $\ch k = 0$.
For a scheme $T/k$, the objects of $\ocM(T)$ are families $(\cS,\cB= \cB_1+\cdots+\cB_n)/T$, where
$\cS/T$ is a flat family of surfaces and each $\cB_i/T$ is a flat family of codimension one subschemes,
such that each closed fibre $(S,B)$ is a pair with semi log canonical singularities and $\omega_S(B)$ ample.
If $\omega_S(B)$ is invertible, then no further conditions are required (this is the case in our example).
In general, one considers the covering stack of $(S,B)$ defined by the $\bQ$-line bundle $\omega_S(B)$ and
requires that the deformation $(\cS,\cB)/T$ is induced by a deformation of this stack, cf. \cite[\S 3]{H}.
The stack $\ocM$ is a proper Deligne--Mumford stack. Let $\oM$ denote its coarse moduli space.

In our example we consider smooth del Pezzo surfaces $S$ of degree $9-n$, $4 \le n \le 8$, with boundary $B$ the sum of the $(-1)$-curves.
Then $K_S+B$ is ample and the pair $(S,B)$ is log canonical if $B$ is a normal crossing divisor. A marking of $S$ corresponds to
a labelling of the $(-1)$-curves.

\begin{Definition}
Let
$B_p \subset \oY(\Delta)$   
be the union of the boundary
divisors which surject onto $\oY(\Delta')$ --- thus for
$D_n$ these are the divisors $D(D_2)$ for $D_2 \not \subset D_{n-1}$,
and for $E_n$ the divisors $D(A_1)$ for $A_1 \not \subset E_{n-1}$.
We call them {\em horizontal divisors}.
\end{Definition}

\begin{Theorem}\label{MainThDegree5}
For $\Delta=D_n$ or $E_n$, $n \leq 6$,
$p: (\oY(\Delta),B_p) \to \oY(\Delta')$ is a flat log smooth
family of stable pairs. It induces an isomorphism
$\oY(\Delta')\stackrel{\sim}{\to} \oM_{0,n-1}$ for $D_n$ and a closed embedding $\oY(\Delta') \subset \oM$ 
for $E_n$ (assuming $\ch k = 0$).
\end{Theorem}

\begin{proof}
This is well-known for $D_n$: $p$ is
the universal family of stable rational curves. For
$\Delta= E_5$, $\oY(\Delta')=\{pt\}$ and the claim is easy. So
it is enough to consider $\Delta = E_6$.

The claims are clear over the interior
$Y(E_5)$ --- it is well known that the $-1$ curves on a del
Pezzo of degree at least $4$ have normal crossings.
For each boundary divisor $D \subset \oY(E_6)$, the
restriction $p|_D$ is described by Prop.~\ref{restrictions}. It
has reduced fibers and is flat of relative dimension $\le2$ and log smooth (instances of the $D_n$ case).
Flatness of $p$ follows from Lemma~\ref{flatcriterion} and
log smoothness from Prop.~\ref{logsmooth}.

Let $T_{\pi}$ be the kernel of $\pi:\,T_{\oY(E_6)} \to T_{\oY(E_5)}$.
Let $F:=\oY(E_5)\times_{X(\cF(E_5))}X(\cF(E_6))$.
By Prop.~\ref{ls}, $T_{\pi} \times \oY(E_6) \to F$ is smooth over $\oY(E_5)$.
It follows that fibers of~$p$ have stable toric singularities.
By Cor.~\ref{n6lcm}, $K_{\oY(E_6)} + B$ is ample.
Since $B = p^*(B) + B_p$,  $K_{\oY(E_6)} + B$ restricts
to the log canonical bundle of each fibre, and thus $p$ is
a family of stable pairs. Consider the corresponding map $\oY(E_5) \to \oM$.
The closed embedding \eqref{sekiembedding} factors through $\oM$ by Lemma~\ref{CRsonKSBAspace}.
Thus $\oY(E_5) \to \oM$
is a closed embedding. \end{proof}

\begin{Lemma}\label{CRsonKSBAspace}

Assume $\ch k = 0$ and $n \le 6$. The KSBA cross-ratios are regular in a neighbourhood of the closure of the locus of smooth marked del Pezzo surfaces of degree $(9-n)$ in the moduli stack $\ocM$.
\end{Lemma}
\begin{proof}
Recall that, for a smooth del Pezzo surface with normal crossing boundary, a KSBA cross-ratio is the cross-ratio of 4 points on a $(-1)$-curve cut out by 4 other
$(-1)$-curves. Consider the universal family $(\bS,\bB) \rightarrow \ocM$. Let $(S,B)$ be a fibre over the closure of the locus of smooth marked del Pezzo surfaces.
We claim that $(B_i, \sum_{j \neq i} B_j|_{B_i})$ is a stable pointed curve of genus $0$ for each $i$.
First, $B_i$ is a nodal curve and the marked points are smooth and distinct
because $(S,B)$ has semi log canonical singularities. Second, $(S,B)$ has stable toric singularities by the first clause of Thms.~\ref{MainThDegree5} and ~\ref{MainThDegree6}, so, in particular, $\omega_S(B)$ is invertible.
It follows that the adjunction formula $\omega_{B_i}(\sum_{j \neq i} B_j|_{B_i})=\omega_S(B)|_{B_i}$ holds \cite[16.4.2]{FA}.
So $\omega_{B_i}(\sum_{j \neq i} B_j|_{B_i})$ is ample and $(B_i,\sum_{j \neq i} B_j|_{B_i})$ is stable, as required.
Because stable curves deform to stable curves, the same is true for nearby fibres. We obtain a map from a neighbourhood of the closure of the locus of smooth del Pezzo
surfaces in $\ocM$ to $\oM_{0,m}$, where $m$ is the number of marked points on $B_i$. Forgetting all but $4$ of the marked points we obtain maps to $\oM_{0,4}$ extending
the KSBA cross-ratios on the locus of smooth del Pezzo surfaces.
\end{proof}

\begin{Definition} \label{eckpoints}
Let $\tY(E_n) \subset X(\tcF(E_n))$ be the closure of $Y(E_n)$, $n=6$ or $7$.
Let $\tp: \tY(E_7) \to \tY(E_6)$ denote the restriction of $\tpi: X(\tcF(E_7)) \to X(\tcF(E_6))$.
An {\em Eckhart point}\/ of $\oY(E_7)$ (resp.~$\tY(E_7)$) is
a point on the intersection of three horizontal $D(A_1)$ divisors
(resp.~their strict transforms).
\end{Definition}

\def\nN{\oN}
\def\nS{\oS}

\begin{Proposition} \label{flatoutsideT}
$(\oY(E_7),B_p) \arrow^p \oY(E_6)$ has the following properties: it is flat
with reduced fibers outside of the union $D$ of
non-flat $D(A_3^{\times2})$ divisors; log smooth outside of $D$ and Eckhart points;
at~an Eckhart point $z$ away from $D$, $p$ is smooth and
$B_p$ restricts on the fibre to $3$
pairwise transversal curves intersecting at~$z$.
\end{Proposition}

\begin{proof}
$W(E_6)$ acts transitively on Fano simplices of $\cR(E_7)$
(since $E_6$ contains at most $4$ orthogonal $A_1$'s, $E_7\setminus E_6$ at most $3$,
and $4$ orthogonal $A_1$'s are contained in a unique Fano simplex).
This implies the result on the purely $A_1$ locus in $\oY(E_7)$,
as each fibre $(S,B)$ of $p|_V$, for a chart $V$ of Lemma~\ref{A1locus}, is an open subset of $\bP^2$ of Fig.~\ref{fanosimplex}
with the  lines through one of the points as the boundary.

By Prop.~\ref{restrictions}, $p|_{D(\Gamma)}$
is flat of relative dimension $\le2$, has reduced fibers, and is log smooth for $\Gamma=A_7$, a flat $A_3 \times A_3$,  or a flat $A_2$.
By Prop.~\ref{logsmooth}, it suffices to prove that the map $p|_{D(A_2)}$ has properties of Prop.~\ref{flatoutsideT}
for any $A_2 \subset E_6 \subset E_7$.
By Prop.~\ref{restrictions}, it suffices to prove this
for the map \eqref{funnyM07}, which we denote by $q$.


$q|_{M_{0,7}}$ is smooth,
and thus by Prop.~\ref{logsmooth}
we just have to consider how $q$ restricts to the possible boundary divisors of $\oM_{0,7}$.
One checks (by running through the list of
possibilities) that $q$ restricts to a flat log smooth map onto a
stratum of $\oM_{0,4}^{\times2}$ except for the boundary divisors
given by subsets of form $\{1,5,4\}$, or $\{1,5\}$
(for example, $\{1,4\}$ gives the boundary divisor $\oM_{0,\{2,\ldots,7\}}$
and $q$ restricts to the canonical projection $\pi_{\{4,5,6,7\}}$).
It is easy to check using Prop.~\ref{boundarycorrespondence} that divisors of type $\{1,5,4\}$
(resp.~of type $\{1,5\}$) are precisely restrictions of non-flat $D(A_3 \times A_3)$ divisors
(resp.~vertical $D(A_1)$ divisors).
Consider the open set $U \subset \oM_{0,7}$,
the complement to all the boundary divisors not of the type $\{1,5\}$.
It is enough to prove the claim for $q|_U$. Note $U$ has $9$ boundary divisors, $\{i,j\}$ with
$i<4$, $j > 4$, and that two meet iff the subsets are disjoint.
On their intersection, say $\delta_{1,5} \cap \delta_{2,6}$,
$$
\begin{CD}
q|_{\delta_{1,5} \cap \delta_{2,6}}:\,\oM_{0,\{1,2,3,4,7\}} @>{\pi_{1,2,3,4} \times \pi_{1,2,4,7}}>>
\oM_{0,4} \times
\oM_{0,4}
\end{CD}
$$
is the blowup of $\bP^1 \times \bP^1$ at three points along
the diagonal. The exceptional divisors are outside of $U$. So the
map is an isomorphism, i.e., the intersection of any two
boundary divisors gives a section. Thus $q$ is smooth in a neighborhood
of any codimension two stratum of $U$. The restriction to the
interior of any boundary divisor is easily seen to be smooth.
Further, each codimension three stratum is carried isomorphically
to the diagonal of $\oM_{0,4} \times \oM_{0,4}$ (which is not
a stratum). It follows that $q|_U$ is flat, and log smooth
outside of the codimension three strata.
\end{proof}

\begin{Proposition} \label{flatandlogsmooth} $\tp:\, \tY(E_7) \to \tY(E_6)$ is flat,
has reduced fibers, and is log smooth outside of Eckhart points. At Eckhart points
the map is smooth, and $B_p$ restricts on the fibre to three
pairwise transversal curves intersecting at a point.
\end{Proposition}

\begin{proof}
By Prop.~\ref{flatoutsideT} $p$ is flat away from the
non-flat $D(A_3 \times A_3)$ divisors. It follows that on
this locus, $\tp$ is just the pullback \cite[2.9]{Tevelev}. Thus (as log smoothness and reduced fibers
are preserved by pullback)
Prop.~\ref{flatoutsideT} implies that $\tp$ has the required properties
outside the inverse image of non-flat $D(A_3 \times A_3)$ divisors.

Fix a non-flat divisor $D(\Theta)$, where $\Theta=A_3 \times A_3'$, $A_3\subset E_6$.
Let $Z\subset\oY(E_6)$ be its image, $Z=D(A_1^{(1)})\cap D(A_1^{(2)})\simeq\oM_{0,5}$.
By Prop.~\ref{restrictions},
$D(\Theta) = \oM_{0,5}(A_3) \times \oM_{0,5}(A_3') \times \oM_{0,4}$
and $p$ restricts to the projection onto the first factor.

\begin{Lemma}\label{thetapicture}\label{Claim}
$p^{-1}D(A_1^{(i)})$ near $D(\Theta)$ has
irreducible components $D(\Theta)$, $D(A_1^{(i)})$, $D(A_7^{(i)})$, $D(A_2^{(i)})$, $D(\tA_2^{(i)})$:
take all $A_7$'s containing~$\Theta$ and all $A_2$'s contained in~$A_3'$.
For any of these divisors $D\ne D(\Theta)$, the map $p|_D$ is log smooth and flat near~$D(\Theta)$.
The intersection of $D(\Theta)$ with $D(A_7^{(i)})$
(resp.~with one of $D(A_1^{(i)})$, $D(A_2^{(i)})$, $D(\tA_2^{(i)})$) is pulled back from a point
of the $\oM_{0,4}$ component
(resp.~from a divisor of the $\oM_{0,5}(A_3')$ component) of~$D(\Theta)$
according to the following picture, where
we realize $\oM_{0,5}$ as the blow up of $\bP^2$ in $4$ points (corresponding to $A_2$ divisors).
Let $Q$ be the union of codimension~$2$ strata in $p^{-1}(Z)$ (near $D(\Theta)$)
not contained in $D(\Theta)$. Dotted lines illustrate $D(\Theta)\cap Q$.
\begin{figure}[htbp]
  \includegraphics[width=4.5in]{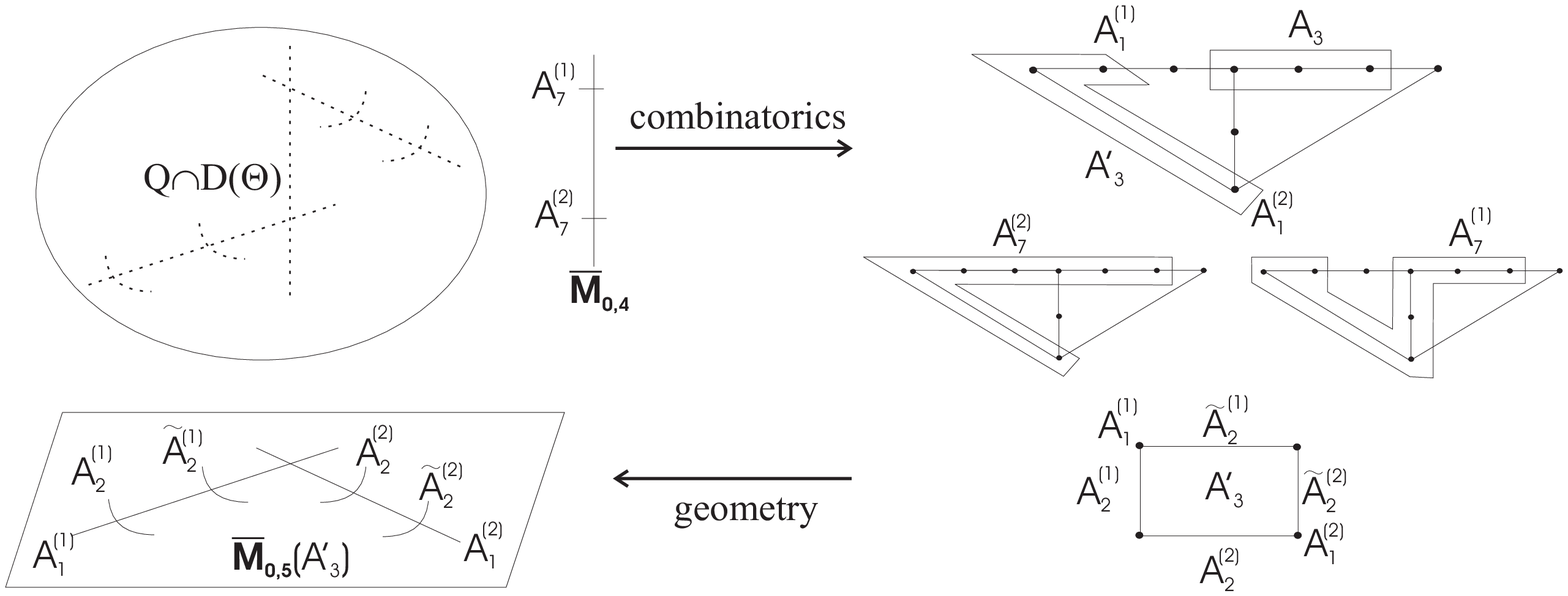}
  \caption{\small}\label{prettypicture}
\end{figure}

\noindent
\end{Lemma}

\begin{proof}
The description of $p^{-1}(D(A_1^{(i)}))$ follows from Prop.~\ref{piE7toE6}.
Intersections of divisors are computed in Prop.~\ref{boundarycorrespondence}.
By Prop.~\ref{restrictions} we only have to prove that $p|_{D(A_1^{(i)})}$
is log smooth and flat.
But $D(A_1^{(i)}) \cap D(\Theta)\simeq\oM_{0,5} \times \oM_{0,4} \times \oM_{0,4}$
and $p|_{D(A_1^{(i)}) \cap D(\Theta)}\simeq\pr_1$. Clearly log smooth and flat
of relative dimension~$2$.
Thus $p|_{D(A^{(i)}_1)}$ is log smooth by Prop.~\ref{logsmooth}
and flat by Lemma~\ref{flatcriterion} (since $p$ is flat outside $D(\Theta)$).
\end{proof}

By Prop.~\ref{piE7toE6} and Th.~\ref{monhubsch} we have a commutative diagram
$$\begin{CD}
\tY(E_7) @>>> \hY(E_7) @>{b_7}>> \oY(E_7) \cr
@V{\tp}VV @V{\hp}VV @VpVV \cr
\tY(E_6) @>>> \hY(E_6) @>{b_6}>> \oY(E_6) \cr
\end{CD}$$
where $b_6$ (resp.~$b_7$) is the blowup of $Z$ (resp.~of $W:=p^{-1}(Z)$).
For any divisor $D$ on $\oY(E_6)$ or $\oY(E_7)$, we denote by $\hD$ its strict transform.
Prop.~\ref{flatandlogsmooth} follows from an analogous statement for $\hp$,
which is proved in Lemma~\ref{Eckhart} below.
\end{proof}

\begin{Remark}
It is clear from Fig.~\ref{prettypicture} and transversality that
$I_W = I_{D(\Theta)} \cdot I_Q$ and $Q$ is regularly embedded, near $D(\Theta)$.
Thus blowing up $W$ is the same as blowing up~$Q$
and $\hD(\Theta)=b_7^{-1}(D(\Theta))$ as $D(\Theta)$ does not contain any component of $Q$.
Note that, at the generic point of intersection of irreducible components of~$Q$,
the blow-up is locally a product of $\bA^3$ with $S$,
the blow up of the union of two lines in $\bA^3$ through the origin.
Note that $S$ has an ordinary double point.
\end{Remark}

\begin{Lemma} Suppose that $D \subset\partial \oY(E_7)$ surjects onto $D(A_1^{(i)}) \subset \oY(E_6)$.
Then, near $D(\Theta)$, $\hD \simeq D$, $\hD(A_1^{(i)})\simeq D(A_1^{(i)})$,
and $\hD \to \hD(A_1^{(i)})$ is flat and log smooth.
\end{Lemma}

\begin{proof}
We blow up a divisor of $D(A_1^{(i)})$,
so $\hD(A_1^{(i)})\to D(A_1^{(i)})$ is an isomorphism.
It~is clear from Fig.~\ref{prettypicture} that
the scheme-theoretic intersection $Q\cap D$ is a divisor of~$D$,
so $\hD \to D$ is an isomorphism. The last remark now follows from Lemma \ref{Claim}.
\end{proof}

$A_3'$ contains $6$ $A_1$'s, $A_1^{(1)},\ldots,A_1^{(6)}$, divided into three orthogonal pairs.
One pair is $A_1^{(1)},A_1^{(2)}$. Taking the other two pairs, and adding $\Theta^{\perp}$ to each, gives
two triples of orthogonal $A_1$'s in $E_7\setminus E_6$, i.e., Eckhart triples.
Note $\hD(A_1) = b_7^{-1}(D(A_1))$ for any $A_1 \not \subset E_6$,
since $D(A_1)$ does not contain any component of $Q$ 
(because $D(A_1)\rightarrow \oY(E_6)$ is flat away from non-flat $D(A_3^{\times 2})$'s by Prop.~\ref{flatoutsideT}).

\begin{Lemma} \label{Eckhart}
Along $\hD(\Theta)$, $\hp$ is flat,
and log smooth away from the strata $\hD(A_1)\cap\hD(A_1')\cap\hD(\Theta^{\perp})$
for the two Eckhart triples above. Near an Eckhart triple
$\hp$ is smooth, and each of the three $\hD(A_1)$'s from the triple meets
the fibre in a smooth curve.
These curves are pairwise transversal and intersect at a point.
\end{Lemma}

\begin{proof}
By Lemma~\ref{normalbundleZ}, the exceptional divisors are products $E_{b_6}=Z\times\bP^1$ and $E_{b_7}=Q\times\bP^1$.
By Lemma \ref{Claim}, for a
boundary divisor $D$ mapping onto $D(A_1^{(i)})$, the restriction $p|_D$ is flat and log smooth near
$D(\Theta)$. In
particular each component of $Q$ is then flat and
log smooth over its image. And thus the same
holds for irreducible components of $E_{b_7}$, since
$E_{b_7} \to E_{b_6}$ is a pullback of $Q \to Z$ (indeed just
product with $\bP^1$). Now by
Prop.~\ref{logsmooth}, log smoothness holds around $E_{b_7}$.
Moreover, by the previous Lemma, we have log smoothness
around any $\hD$.

So it's enough to work outside of $E_{b_7}$, where
$b_7$ is the identity, and away from any of the
divisors $D$ of Lemma \ref{thetapicture}.
On this open subset $D(\Theta)$ is isomorphic to
$Z\times U \times V$, where $U\subset\bP^2$ is the complement of lines $L_i$, $i=1,2$
(that correspond to $A_1^{(i)}$ on Fig.~\ref{prettypicture})
and $V=\bP^1\setminus\{P_1,P_2\}$ (where $P_i$ corresponds to $A_7^{(i)}$).

Now we can describe the map
$\tp: D(\Theta) \to E_{b_6}=Z\times\bP^1$ (on the open set on which we are working).
As it is the same over each point of $Z = \oM_{0,5}$, we describe it as a rational map
$f:\,\bP^2\times\bP^1\to \bP^1$.
Since the projective normal bundle of~$Z$ is trivial,
it is the rational
map given by the linear system spanned by the two
divisors $p^{-1}(D(A_1^{(i)})) - D(\Theta)$, $i=1,2$, restricted on $D(\Theta)$. On
$\bP^2 \times\bP^1$ these divisors are $(L_i \times \bP^1)\cup(\bP^2 \times P_i)$.
Thus in coordinates $f$ is given by
$$
f:\,((X:Y:Z),\ (A:B) ) \to (XA:YB)
$$
where the lines $L_i$ are given by $X=0$, $Y=0$, and
the points $P_i$ are are given by $A=0$, $B=0$.
On $U\times V$ this map is clearly
smooth, of relative dimension~$2$. This proves flatness.
For log smoothness, we consider how the boundary
divisors meet the fibres. There are $5$ boundary divisors that meet $U\times V$,
namely $D(A_1^{(i)})$ for $i=3,4,5,6$, and $D(\Theta^{\perp})$.
On Fig.~\ref{prettypicture}, the first $4$ divisors project on lines in~$\bP^2$
pairwise connecting $4$ blown up points, and the last is the inverse
image of a point from $\bP^1$, corresponding to $D(\Theta^{\perp})$.
Clearly the picture on each
fibre of $f$ is of smooth divisors on a smooth surface that either intersect transversally
or three of them meet at a point and are pairwise transversal.
The triple points come from the orthogonal pairs of $A_1^{(i)}$, intersecting
with $D(\Theta^{\perp})$, which are exactly the
Eckhart triples mentioned above.
\end{proof}

\begin{Lemma} \label{subprod} Let $P$ be product
of smooth (not necessarily proper) curves of log
general type. Let
$S \to P$ be a map with $S$ smooth, which has a generically
injective derivative. Then $S$ is of log general type.
\end{Lemma}

\begin{proof} We may replace $P$ by $\dim S$ factors and
assume the map is generically~\'etale.
If $C$ is a curve of log general type then $\Omega^1_C$ is globally generated by log $1$-forms.
It~follows that $\omega_P$ is globally generated by log canonical forms (by wedging),
so the result now follows by pulling back forms.
\end{proof}

\begin{Proposition} \label{A_1^2CR}
Let $Z=D(A_1)\cap D(A_1')\subset\oY(E_6)$ be a codimension~$2$ stratum.
Then there is a KSBA cross-ratio $g$ such that near $Z$ we have
$(g)=D(A_1)-D(A_2)$.
\end{Proposition}

\begin{proof}
Let $M$ denote the lattice spanned by $D_4$ and let
$\lambda,\rho,\mu,\nu$ be the characters corresponding to simple roots
$\eps_1-\eps_2$, $\eps_2-\eps_3$, $\eps_3-\eps_4$, $\eps_3+\eps_4$.
Let $T=\Hom(M,\bG_m)$.
For $\alpha \in D_4^+$ let $H_{\alpha}$ denote the hypertorus
$(\chi^{\alpha}=1) \subset T$.
Cayley \cite{C} found an explicit family of smooth cubic
surfaces over $T \setminus \bigcup H_{\alpha}$, see also \cite[5.1]{Naruki}:
\begin{equation}\label{family}
\begin{array}{c}
\rho W(\lambda X^2 + \mu Y^2 + \nu Z^2 +
(\rho-1)^2(\lambda\mu\nu\rho-1)^2W^2\\
+(\mu\nu+1)YZ+(\lambda\nu+1)XZ
+(\lambda\mu+1)XY
\\-(\rho-1)(\lambda\mu\nu\rho-1)W((\lambda+1)X+(\mu+1)Y+(\nu+1)Z)) +XYZ = 0.
\end{array}\cooltag
\end{equation}
The marking of the family \eqref{family} is given by labelling the tritangent planes
of the fibres \cite[Table~1, p.~10]{Naruki}.
By \cite[Prop.~7.1]{Naruki}, the morphism $T \setminus \bigcup H_{\alpha} \rightarrow Y(E_6)$ defined
by the family (\ref{family}) is an open embedding.

\begin{Lemma} \label{chars}
The character $\chi^\alpha$ for $\alpha\in D_4$
is a KSBA cross-ratio.
\end{Lemma}

\begin{proof}
We use the equations of the tritangent planes for the family
(\ref{family}) given in \cite[Table~1, p.~10]{Naruki}.
A fibre of the family contains the line $l=(W=X=0)$, and has tritangent
planes $(W=0)$, $(X=0)$ and $(X+\rho(\mu-1)(\nu-1)W=0)$
through~$l$. These give line pairs intersecting $l$ in the points
$Z/Y=\{0,\infty\}$, $\{-\mu,\frac{-1}{\nu}\}$, and
$\{-1,\frac{-\mu}{\nu}\}$ respectively. Thus $\mu=\eps_3-\eps_4$ is a
KSBA cross-ratio (the cross-ratio of $0,-1,-\mu, \infty$).
By $W(D_4)$-equivariance, the same is true for any $\alpha\in D_4$.
\end{proof}

We recall the construction of $\oY(E_6)$ from \cite[p.~21--24]{Naruki}.
Let~$\Sigma$ be the fan of Weyl chambers in $N_{\bR}$ given by the
Coxeter arrangement of root hyperplanes in $D_4$.
Let $X=X(\Sigma)$ be the corresponding $T$-toric variety.
Let $\oH_{\alpha} \subset X$ be the closure of~$H_\alpha$.
$X$ is smooth and the arrangement of $\oH_\alpha$'s at the identity $e \in T$
is locally isomorphic to the Coxeter arrangement of hyperplanes in $N \otimes k$.
Let $\tilde{X} \rightarrow X$ be the blowup of this arrangement
which corresponds to the ``wonderful blowup'' of the Coxeter arrangement.
That is, we blowup the strata corresponding to irreducible root
subsystems $\Gamma\subset D_4$ in order of increasing dimension.
Explicitly, we blowup the identity $e \in T$, 12 curves, and 16
surfaces, which correspond to
sub root systems of type $D_4$, $A_3$, and $A_2$ respectively.
Following Naruki, we call these strata
$A_3$-curves and $A_2$-surfaces. Then the KSBA  cross-ratios of type $I$ define a
morphism $\tilde{X} \rightarrow \oY(E_6)$ which contracts
the strict transforms of the exceptional divisors over the $A_3$-curves
to $A_1^2$-strata in $\oY(E_6)$ and is an isomorphism elsewhere.

Let $A_3=\{\eps_i-\eps_j\}\subset D_4$ and let $\Gamma \subset X$ be the corresponding $A_3$-curve.
$\Gamma$~meets the toric boundary in two points, which are interior
points of the toric boundary divisors
$\Delta_1$ and $\Delta_2$ corresponding to the positive and
negative rays of the line $A_3^{\perp}=\bR \cdot (1,1,1,1)$
in the fan of $X$. The primitive generators of these rays in
$N=M^{\vee}$ are $\pm \frac{1}{2}(1,1,1,1)$.
So, if $\alpha =\eps_1+\eps_2 \in D_4$ then the divisor
$(\chi^\alpha)$ equals $\Delta_1-\Delta_2$ plus some other toric boundary divisors.
On the Naruki space $\Delta_1$ and $\Delta_2$ become $A_1$-divisors
$Y_1$ and $Y_2$ meeting in a $A_1^2$ stratum $Z$,
and $g=\chi^{\alpha}$ is a cross-ratio with the desired properties.
\end{proof}

\begin{Theorem} \label{MainThDegree6}
Let $\eY(E_7) \to \tY(E_7)$ be the blowup of the union
of Eckhart points. The map
$(\eY(E_7),B_{\ep}) \arrow^{\ep} \tY(E_6)$ is a flat
family of stable pairs with stable toric singularities.
For $\ch k = 0$, the induced map $\tY(E_6)\to\oM$
is a closed embedding, with image the closure of the locus of pairs
consisting of a smooth cubic surface without Eckhart points with boundary its 27 lines.
The product of all KSBA cross-ratio maps embeds $\tY(E_6)$ in the product of~$\bP^1$'s.
\end{Theorem}

\begin{proof}
By Prop.~\ref{flatandlogsmooth}, the union
of Eckhart points is a disjoint union of smooth codimension three strata contained in the smooth locus of $\tp$.
Fibres of $\ep$ have stable toric singularities by
Prop.~\ref{flatandlogsmooth}, Th.~\ref{monhubsch} and Prop.~\ref{ls}.
Let $(F,B)$ be a fibre, and $(\tF,B)$ the corresponding
fibre of $\tp$. The inverse image of an Eckhart point of $(\tF,B)$ is $\bP^2$
with $3$ general lines on it, glued to the rest of $F$ along the fourth general line.
Thus $K_F + B$ restricts to an ample divisor on the exceptional locus for $F\to\tF$.
To prove $K_F+ B$ is ample, by Th.~\ref{lcm} it is enough to
show that any open stratum $S$ of $(\tF,B)$ is log minimal.
Let $\tG\subset X(\tcF(E_7))$ be
the fibre of $\tilde\pi$ containing $\tilde F$. The
structure map $T_{\pi} \times \tF \to \tG$ is smooth outside of Eckhart
points,
and the stratification of $\tF$ is induced from the stratification of
$\tG$ by restriction. Strata of $\tG$ are orbits for $T_\pi$, thus
in particular, $S$ is smooth. Since $X(\tcF(E_7))$ embeds
in the normalization of $X(\tcF(E_6)) \times_{X(\cF(E_6))}
X(\cF(E_7))$,
$\tF$~maps finitely to a fibre of $p$, and $S$ has a quasi-finite map
to an open stratum of~$\oY(E_7)$. Thus $S$ has a quasi-finite map
to a product of $M_{0,4}$'s by Th.~\ref{monhubsch}, and so is
log minimal by Lemma~\ref{subprod}. Thus $\ep$ is a flat family of
stable pairs. Let $\tY(E_6)\to\oM$ be the induced map.

By Lemma~\ref{CRsonKSBAspace}, the KSBA cross-ratio maps are
regular on $\ocM$, and therefore on $\tY(E_6)$.
Thus it suffices
to prove that the product
of KSBA cross-ratio maps embeds $\tY(E_6)$.
By Th.~\ref{monhubsch} the cross-ratios of type I define an embedding of $\oY(E_6)$.
So to prove that all cross-ratios embed $\tY(E_6)$ we can work locally over $\oY(E_6)$.
By Th.~\ref{monhubsch}, and because an analogous statement holds for the map
of toric varieties $X(\tF(E_6))\to X(\cF(E_6))$,
$\tY(E_6)\to\oY(E_6)$ is the blow up of all purely $A_1$
strata in order of increasing dimension.
Let $P \in \oY(E_6)$ be a point lying on exactly $m$ $D(A_1)$-divisors and
write locally
$(P \in \oY(E_6))=(0 \in \bA^4_{x_1,\ldots,x_4})$, where
$(x_1=0),\ldots,(x_m=0)$ are the
$D(A_1)$-divisors through $P$.
An easy toric calculation shows that the
rational functions $x_i/x_j$, $1 \le i<j \le m$,
embed $\tY(E_6)$ in $\oY(E_6)\times (\bP^1)^m$ (locally near $P$)).
On the other hand, Prop.~\ref{A_1^2CR} shows that these rational functions
correspond (locally) to KSBA cross-ratio maps on $\oY(E_6)$.
\end{proof}

\medskip
\noindent
Paul Hacking,  Department of Mathematics, University of Washington, Box 354350, Seattle, WA~98195. \texttt{hacking@math.washington.edu} \\
\\
Sean Keel, Department of Mathematics, University of Texas at Austin, Austin, TX~78712. \texttt{keel@math.utexas.edu} \\
\\
Jenia Tevelev, Department of Mathematics and Statistics, University of Massachusetts Amherst, Amherst, MA~01003. \texttt{tevelev@math.umass.edu}\\

\end{document}